\input amstex
\documentstyle{amams} 
\document
\def\ritem#1{\item{{\rm #1}}}
\input boxedeps.tex 
\SetepsfEPSFSpecial 
\HideDisplacementBoxes
\def\figin#1#2{
$$
 {\BoxedEPSF{#1.ps scaled
#2}%
}%
$$
\noindent}

\annalsline{156}{2002}
\received{June 23, 1998}
\startingpage{713}
\font\tenrm=cmr10

\catcode`\@=11
\font\twelvemsb=msbm10 scaled 1100

\font\ninemsb=msbm10 scaled 800
\newfam\msbfam
\textfont\msbfam=\twelvemsb  \scriptfont\msbfam=\ninemsb
  \scriptscriptfont\msbfam=\ninemsb
\def\msb@{\hexnumber@\msbfam}
\def\Bbb{\relax\ifmmode\let\next\Bbb@\else
 \def\next{\errmessage{Use \string\Bbb\space only in math
mode}}\fi\next}
\def\Bbb@#1{{\Bbb@@{#1}}}
\def\Bbb@@#1{\fam\msbfam#1}
\catcode`\@=12

 \catcode`\@=11
\font\twelveeuf=eufm10 scaled 1100
\font\teneuf=eufm10
\font\nineeuf=eufm7 scaled 1100
\newfam\euffam
\textfont\euffam=\twelveeuf  \scriptfont\euffam=\teneuf
  \scriptscriptfont\euffam=\nineeuf
\def\euf@{\hexnumber@\euffam}
\def\frak{\relax\ifmmode\let\next\frak@\else
 \def\next{\errmessage{Use \string\frak\space only in math
mode}}\fi\next}
\def\frak@#1{{\frak@@{#1}}}
\def\frak@@#1{\fam\euffam#1}
\catcode`\@=12


\define\ext{\oper{Ext}}

\define\imm{\oper{Im}}

\define\ree{\oper{Re}}

\define\tec{Teich\-m\"ul\-ler\ }
\define\wei{Weierstrass\ }

%

%
%

\define\a{\alpha}
\predefine\barunder{\b}
\redefine\b{\beta}

\define\cd{\cdot}
\predefine\dotunder{\d}
\redefine\d{\delta}
\define\Del{\Delta}
\define\e{\varepsilon}

\define\f{\frac}
\define\g{\gamma}
\define\G{\Gamma}

\define\lb\{{\left\{}

\define\lra{\longrightarrow}

\define\om{\omega}

\define\oper{\operatorname}

\define\ov{\overline}

\define\p{\partial}
\define\rb\}{\right\}}

\define\vp{\varphi}
\define\wh{\widehat}

\define\x{\times}
\define\z{\zeta}
\define\({\left(}
\define\){\right)}
\define\[{\left[}
\define\]{\right]}
\define\<{\left<}
\define\>{\right>}
\def\slantline#1#2#3#4#5{\hbox to 0pt{
}}


\def\SC{{\Cal C}}

\def\SF{{\Cal F}}

\def\SH{{\Cal H}}

\def\SK{{\Cal K}}

\def\SM{{\Cal M}}

\def\SR{{\Cal R}}

\def\SY{{\Cal Y}}



\def\bbc{{\Bbb C}}

\def\bbe{{\Bbb E}}

\def\bbr{{\Bbb R}}

\define\per{\operatorname{Per}}
\define\height{\operatorname{\SH}}

\def\sgn{\oper{sgn}}

\def\supp{\oper{supp}}

\define\id{\oper{id}}

\define\im{\oper{Im}}

\define\ogup{{\Omega_{Gdh}}}
\define\ogdn{{\Omega_{G^{-1}dh}}}

\title{Teichm\"uller theory and handle addition\\ for minimal surfaces} 
\shorttitle{Teichm\"uller theory}  
\acknowledgements{The second author was partially supported by NSF grant number DMS-9626565
and the SFB.}
 \twoauthors{Matthias Weber}{Michael Wolf}
  \institutions{University of Indiana,
Bloomington, IN\\
{\eightpoint {\it E-mail address\/}: matweber@indiana.edu}
\\
\vglue6pt
Rice University, Houston, TX\\
{\eightpoint {\it E-mail address\/}: mwolf@rice.edu
}}

 \def\sni#1{\smallbreak\noindent{#1}. }
\def\ssni#1{\vglue-1pt\noindent\hskip18pt {#1}}

\sni{1} Introduction
\ssni{1.1.} The surfaces
\ssni{1.2.} The proof
\sni{2} Background, notation and a sketch of the argument
\ssni{2.1.} Minimal surfaces
\ssni{2.2.} \tec theory
\ssni{2.3.} A brief outline of the proof
\sni{3} The geometry of orthodisks
\ssni{3.1.} Orthodisks
\ssni{3.2.} From orthodisks to Riemann surfaces
\ssni{3.3.} From orthodisks to \wei data
\ssni{3.4.} Geometric significance of the formal \wei data
\ssni{3.5.} Examples of simple orthodisks
\ssni{3.6.} Orthodisks for the  Costa towers
\ssni{3.7.}  More orthodisks by drilling holes
\sni{4} The space of orthodisks
\ssni{4.1.} Introduction
\ssni{4.2.} Geometric coordinates for the ${\rm DH}_{1,1}$-surface
\ssni{4.3.} Height function for the ${\rm DH}_{1,1}$-surface
\ssni{4.4.} Geometric coordinates for the ${\rm DH}_{m,n}$-surfaces
\ssni{4.5.} Height functions for the ${\rm DH}_{m,n}$-surfaces
\ssni{4.6.} Properness of the height functions for the ${\rm DH}_{m,n}$-surfaces
\ssni{4.7.} A monodromy argument 
\sni{5} The gradient flow
\ssni{5.1.} Overall strategy
\ssni{5.2.} Deformations of ${\rm DH}_{1,1}$
\ssni{5.3.} Infinitesimal pushes
\sni{6} Regeneration
\sni{7} Nonexistence of the ${\rm DH}_{m,n}$-surfaces with $n<m$
\sni{8} Extensions and generalizations
\ssni{8.1.} Higher symmetry
\ssni{8.2.} Deformations with more catenoidal ends
\ssni{8.3.} Embeddedness aspects of ${\rm DH}_{m,n}$
\sni{9} References

\section{Introduction}

In this paper, we develop \tec theoretical methods to construct new
minimal surfaces in $\bbe^3$ by adding handles and planar ends to
existing minimal surfaces in $\bbe^3$. We exhibit this method on an
interesting class of minimal surfaces which are likely to be embedded,
and have a low degree Gau\ss map for their genus.
In particular, we exhibit a two-parameter family of complete minimal
surfaces in the Euclidean three-space $\bbe^3$; these surfaces are
embedded (at least) outside a compact set and are indexed (roughly)
by the number of ends they have and their genus. They have at most
eight self-symmetries despite being of arbitrarily large genus, and
are interesting for a number of reasons. Moreover, our methods
also extend to prove that some natural candidate classes of surfaces
cannot be realized as minimal surfaces in $\bbe^3$. As a result of both
aspects of this work, we obtain a classification of a family 
of surfaces as either realizable or unrealizable as minimal surfaces.

This paper is a continuation of the study we initiated in [WW]; in a
strong sense it is an extension of that paper, as the essential
organization of the proof, together with many details, have been
retained. Indeed, part of our goal in writing this paper was a
demonstration of the robustness of the methods of [WW], in that here
we produce minimal surfaces of a very different character than those
produced in [WW], yet the proof changes only in a few quite
technical ways. (In particular, the present proof handles the
previous case of Chen-Gackstatter surfaces of high genus as an
elementary case.) Indeed in the intervening years between our initial
preparation of this manuscript and its final revision for 
publication, this method has been applied to produce other families
of surfaces of substantively different characteristics or to
prove their nonexistence (\cite{WW2},
\cite{MW}).

 \vglue12pt
1.1. {\it  The surfaces}. Hoffman and Meeks (see \cite{Ho-Me}) have
conjectured that any complete embedded minimal surface in space has
genus at least $r-2$, where $r$ denotes the number of ends of the
surface. In this paper, we provide significant evidence for this
conjecture in the situation where the surfaces have eight
symmetries. This is an important case for two reasons: first, it is
presently unknown whether there are any complete embedded minimal
surfaces which have no symmetries
\footnote{$^{1}$}{Added in proof. M.\ Traizet [Tr] has announced the proof of the 
existence of a complete embedded minimal surface with no 
symmetries.}, and second, there are very few
families of examples known where there are more than four ends.
(Indeed, the only such constructions available are from the recent
work of Kapouleas [Kap], where the genus is both high and
inestimable.)

In particular, we consider two families of surfaces, with the first
included in the second. The first case consists of surfaces $CT_g$
which generalize Costa's example [Cos]. We prove

\nonumproclaim{Theorem A} For all odd genera $g${\rm ,} there is a complete
minimal surface $CT_g\subset\bbe^3$ which is embedded outside a
compact surface with boundary of genus~$g${\rm ,} with $g$ parallel {\rm (}\/horizontal\/{\rm )} planar ends
and two catenoid ends. The symmetry group of $CT_g$ is generated
by reflective symmetries about a pair of orthogonal vertical planes
and a rotational symmetry about a horizontal line.
\endproclaim

These surfaces represent the borderline case for the conjecture.
(The even genus cases have substantially different combinatorics, and
require a different treatment.) Consider the Riemann surface
underlying such an example: it is a fundamental theorem  of
Osserman \cite{Oss1} that such a surface is conformally a compact surface
of genus~$g$, punctured at points corresponding to the ends. Let $Z$
denote the vertical coordinate of such a minimal surface: clearly,
$Z$ is critical at the $g$ points corresponding to the planar ends,
the two points corresponding to the catenoid ends, and $g$ interior
points where the two reflective planes meet the surface.

We generalize these surfaces as follows, imagining {\bf D}rilling
additional\break {\bf H}oles to obtain surfaces ${\rm DH}_{m,n}$ (see \S3.7).

\nonumproclaim{Theorem B} {\rm(i)} For every pair of integers $n\ge
m\ge1${\rm ,} there exists a complete minimal surface ${\rm DH}_{m,n}\subset\bbe^3$
of genus $m+n+1$ which is embedded outside a compact set with the
following properties\/{\rm :} it has $2n+1$ vertical normals{\rm ,} $2m+1$ planar
ends{\rm ,} and two catenoid ends. The symmetry group is as in Theorem~{\rm A.}

{\rm(ii)} For $n<m${\rm ,} there is no complete minimal surface with those
symmetries of the type ${\rm DH}_{m,n}$ 
{\rm (}\/and $2n+1$ vertical normals{\rm ,} $2m+1$ planar ends{\rm ,} and two
catenoid ends\/{\rm ).}
\endproclaim

In the second statement, the surfaces for which we prove nonexistence
are in precise analogy with the surfaces for which we prove existence.
There are many configurations of surfaces which have the given
eight symmetries and $2n+1$ vertical normals, $2m+1$ planar
ends, and two catenoid ends, and we will indicate the range of
possible choices in \S3: in the nonexistence section, we concentrate
{\it only} on the candidates which the rest of the paper indicates are
most likely to exist.  In particular, we do not prove a more general
statement ruling out all surfaces of the rough description of
having $2n+1$ vertical normals, $2m+1$ planar
ends, and two catenoid ends (although many of the possible configurations
we do not treat would also have no minimal representatives, with 
the proofs of nonexistence being precisely analogous to the
proof we give in \S7). We give a precise formulation of the statement
of Theorem~B(ii) at the outset of \S7.  

Theorem~A follows from Theorem~B (i) by setting $n=m=\f12(g-1)$. The
weak embeddedness statement in Theorem~B is strengthened somewhat in
\S9; we conjecture (supported by some numerical evidence) that these
surfaces are, in fact, embedded. The restriction to planar ends is
unnecessary: in \S8 we show that these surfaces are deformable to
having catenoid ends. Theorem~B is displayed in tabular form at the
end of \S3. (In that table, we also add some information about the
case $m=0$, which was excluded from the statement of Theorem B.)

In summary, for the case of ``essentially embedded''
surfaces (i.e. those surfaces which are embedded outside
of a compact set) with eight
symmetries and odd ends, the conjecture is robustly true: no
counterexamples (of the type ${\rm DH}_{m,n}$)
may exist for $g<r-2$, and any pair $(r,g)$
describes an example when $g\ge r-2$.

Below are two pictures of the surface ${\rm DH}_{1,2}$, one showing it
completely, the other exhibiting only the central planar end:
\vfill
\centerline{
\BoxedEPSF{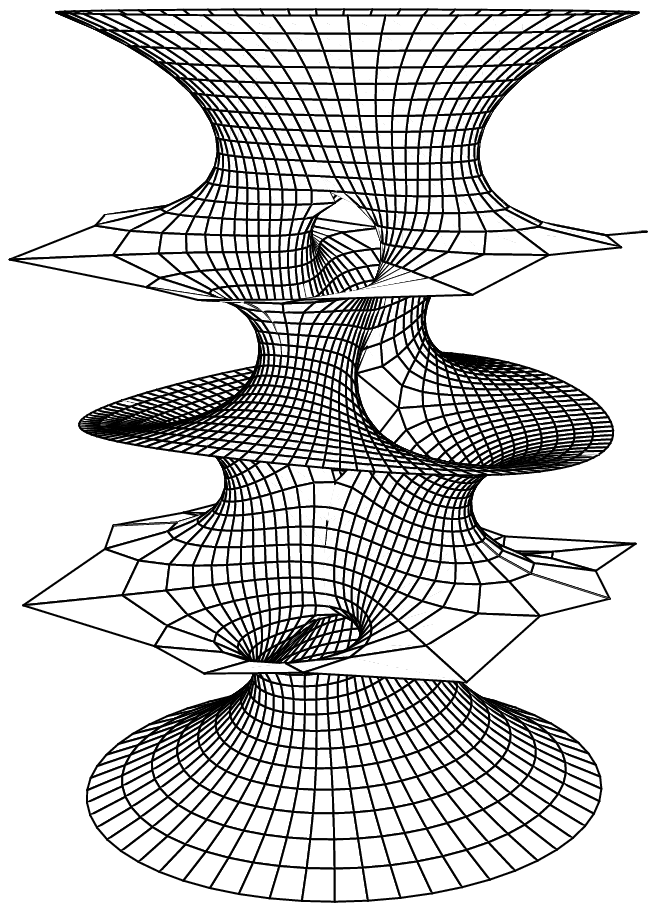 scaled 650}\hglue.25in
\BoxedEPSF{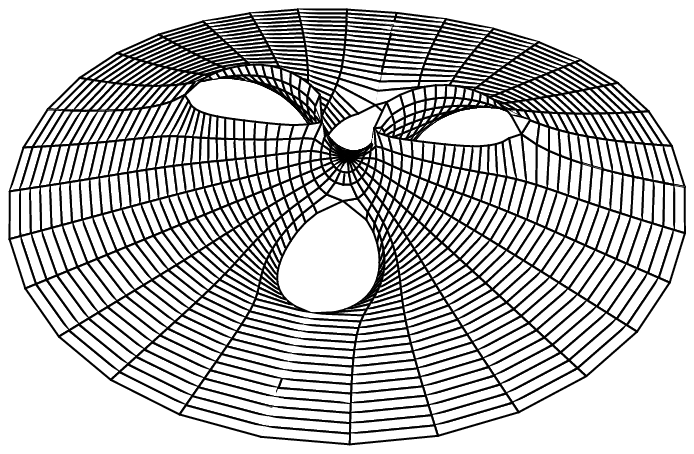 scaled 650}
}
\vfill
\centerline{${\rm DH}_{1,2}$ surface and middle end}
  \eject
1.2. {\it The proof}. As one of our principal goals is the description
of this method, we now give an overview; 
in \S2.3, we give a fairly detailed sketch. In particular, our
discussion in that subsection is quite general, as we aim to 
outline an approach to proving existence results for wide classes
of minimal surfaces, even if we only carry out that plan for
the specific classes described above.

Minimal surfaces in space can be defined in terms of their
Weierstrass data: a pair of meromorphic forms on the underlying
Riemann surface. Conversely, a pair of meromorphic forms on a simply
connected domain naturally determines a (local) minimal surface, up
to some mild compatibility requirement on the divisors. Defining a
minimal surface with some homology is substantially more difficult
however, as this requires a compatibility of periods of three
one-forms (say ($\a_1$, $\a_2$, $\a_3$)) defined via the data. 
A common approach to this global ``period problem'' is to prescribe
sufficient symmetry of the minimal surface so that the conformal
structure of the minimal surface is apparent; one then searches
for appropriate meromorphic 
one-forms $\a_1$, $\a_2$, $\a_3$ on that Riemann
surface.  As our minimal surfaces will have only a few symmetries,
we are unable to determine the conformal structure of the surface
{\it a priori\/}, or even to restrict it sufficiently well so that we might
adapt that common approach; we develop a different approach.
We handle this ``period problem'' by instead studying the
(developed) flat structures associated to the $\a_i$: in this
formulation, the periods are identifiable as vectors in $\bbc$, and
the period problem is soluble when those vectors are compatible (see
\S\S3.2, 3.3 for the precise relationships). In fact, for many
surfaces with interesting shapes, there are large moduli spaces of
triples of flat structures whose geometries are compatible in the
sense that if the underlying Riemann surfaces were conformally
identical, the flat structures would correspond to Weierstrass data
with a ``solved'' period problem.

Thus, we have translated the problem of producing a minimal surface
in $\bbe^3$ with prescribed shape into a problem in \tec theory: in a
moduli space $\Delta$ of compatible triples of flat structures, find a
triple whose underlying Riemann surfaces coincide. This we solve
nonconstructively by introducing a nonnegative height function
$\SH:\Delta\to\bbr^+$ on $\Delta$, which has the features of being proper,
and whose only critical point is at a solution to our problem. The
bulk of the paper is a description of this height function (\S4.3),
a proof of its properness (\S4.6), and a proof that its only
critical points are at solutions\break (\S5, \S6). An interesting
feature of this proof is that it is inductive: the triples
of flat structures for a slightly
less complicated minimal surface lie
on a boundary face of the compactified moduli space 
$\overline\Delta$ of 
compatible triples of flat structures for more complicated
surfaces.  We consider this solution (on $\overline\Delta$)
to the less complicated problem
as the point-at-infinity 
of a particularly good locus in the moduli space
$\Delta$ on which to restrict the height and look for a solution.
This bootstrapping from slightly less complicated solutions to
solutions is the `handle addition' referred to in the title of the
paper.

Along the way, we learn quite a bit about moduli spaces of pairs of
flat structures. The bulk of \S4 is devoted to describing asymptotic
relationships between the underlying conformal structures of flat
disks with alternating horizontal and vertical sides of
corresponding lengths. Most of \S5 concerns estimates on
infinitesimal changes in the conformal structure of a flat structure
given a prescribed infinitesimal change in its Euclidean geometry;
\S6 concerns the effects on a flat structure of the opening
of a node in a Riemann surface.

Our pace in this exposition is occasionally casual, as we have tried
to include a number of motivating examples along the way. Indeed, we
feel these illustrate the power of our approach. For instance, it is
occasionally immediately clear (and even a straightforward
calculation) that the moduli space $\Delta$ for a shape is empty --- we
then conclude that there is no minimal surface with that shape (see
\S3.5.2, \S7). There are also instances when the moduli space is an
identifiable singleton, from which we may conclude immediately that
such a shape is realizable minimally; in \S3.5.1, we find that this
is true for Costa's surface, obviating the need for the analysis of
elliptic functions in the proof of the existence of this
particular surface, and yielding a conceptually new proof of 
the existence  of this 
surface. 

Here is a detailed discussion of the contents of this paper: we
begin in \S2 with some background in \tec theory and minimal
surfaces, the two subjects which we relate in this paper. 
The proof of Theorem~B occupies the bulk of the paper 
(\S3 through \S7, \S9). While the details are occasionally quite technical
and the arguments require some space to present completely, the
basic plan and ideas are rather straightforward; we give a step-by-step
summary in \S2.3. Our
objects of study are minimal surfaces which have sufficient
symmetries so that the fundamental domain for the action is a disk.
An equivariant form on that surface induces a flat structure on the
disk with straight boundaries, and we study such domains
(``orthodisks'') and their moduli spaces in \S3. Also in \S3, we
meet our main new technical obstacle: we allow our domains to
develop into $\bbe^2$ as branched covers of the plane. This introduces many
complications into the analysis in terms of allowing many different
types of geometrically defined motions and degenerations within the
moduli spaces, as well as some difficulties in defining the
frontiers of the moduli space. We deal with these difficulties in
\S4, where we introduce our height function and prove its
properness. In \S6 we find a good locus within moduli space on which
to flow to a solution, and in \S5 we prove that we may flow along
that locus to a solution. In \S7, we prove the nonexistence portion
of the classification Theorem~B, and we conclude in \S8 by
extending some of our results: we extend to cases of higher
dihedral symmetry in \S8.1, and to nonplanar ends in \S8.2. Finally,
in \S8.3, we end 
by offering evidence that our surfaces are embedded; in particular,
we state that the surfaces are regularly homotopic to an embedding,
by a compactly supported regular homotopy. (Space considerations
force us to defer the proofs of the statements in \S8.2 and \S8.3
to the book \cite{WW2}.)

\demo{Acknowledgments} It is once again a pleasure to acknowledge
our indebtedness to Hermann Karcher for many conversations related to
this work, and for a very careful reading which led to a much
improved exposition.  We also appreciate the helpful comments
of the referee. 
\enddemo

\section{Background, notation, and a sketch of the argument}

 2.1. {\it Minimal surfaces}.

\demo{{\rm 2.1.1.} The Weierstrass representation} Any 
complete minimal surface $M$ of finite total curvature in $\bbe^3$ can be
defined
by
$$
z \mapsto \ree \int_\cd^z
(\omega_1, \omega_2, \omega_3)
$$
where $\omega_i$ are three meromorphic $1$-forms on
a compact Riemann surface $\SR$  such that
$$\omega_1^2 + \omega_2^2 + \omega_3^2  \equiv 0.$$

This last condition is usually eliminated by writing
$$\omega_1=\f12(G- G^{-1}) dh, \omega_2=\frac i 2(G+G^{-1})
dh,\omega_3=dh$$
where $G$ is the Gauss map and $dh$ the height differential.
Here $dh$ is a closed (but not necessarily exact) differential
on the underlying Riemann surface $\SR$.

The pair $G$ and $dh$ are called the \wei data for the minimal
surface~$M$.

Significant geometric data attached to such a surface are
the total absolute curvature
$$
\SK := \int_R \vert K \vert dA = 4\pi \cd \hbox{degree of the
Gauss map}
$$
and its Riemannian metric
$$
ds = \(\vert G\vert + \frac 1 {\vert G\vert}\) \vert dh\vert.
$$

At an end, the value of the Gauss map gives the normal to the
asymptotic plane. The metric
becomes infinite by completeness, and the order of decay or
degeneration describes the
type of end: The
only possible embedded ends are Catenoid ends and planar ends for which
$$ds \sim \infty^2.$$
They are distinguished by the Gauss map: for a catenoid end, it is
single-valued and
for a flat end, the Gauss map has higher multiplicity.
{\it For the planar ends we will construct{\rm ,} we will assume that the
Gauss map at a flat end has local degree three.}

Given a Riemann surface $\SR$, a meromorphic function $G$ on $\SR$ and
a meromorphic form $dh$ on $\SR$ for which the metric $ds$ above
is regular except at an acceptable collection of distinguished points
on $\SR$, we can attempt to use the constructed forms $\om_i$ 
and the formula (called the \wei representation)

$$
z \mapsto \ree \int_\cd^z
(\f12(G-G^{-1}) dh,\frac i 2(G+G^{-1})
dh , dh)
$$
to define
a minimal surface in space.  This procedure works locally, but
the surface is only well-defined globally if the periods
$$
\ree \int_\g
(\f12(G-G^{-1}) dh,\frac i 2(G+G^{-1})
dh , dh)
$$
vanish for every cycle $\g \subset \SR$.  The problem of finding
compatible meromorphic data $(G, dh)$ which satisfy the above 
conditions on the periods of $\om_i$ is known as `the period 
problem for the \wei representation'.

We will find it convenient to use the following psychologically different
well-known 
version of the above period conditions and \wei data: we will 
attempt to specify compatible forms $Gdh$ and $G^{-1} dh$ so that
the periods satisfy
$$
\ree \int_\g dh = 0 \tag2.1a
$$ 
\noindent and

$$
\int_\g Gdh = \overline{ \int_\g G^{-1} dh}. \tag2.1b
$$

The equivalence to the original period problem is elementary. 

The forms $\a_1=Gdh, \a_2=G^{-1}dh$, and $\a_3=dh$ 
lead to singular flat structures on the underlying Riemann
surfaces, defined via the line elements $ds_{\a_i} = \vert \a_i
\vert$.  These singular metrics are flat away from the support of
the divisor of $\a_i$; on elements $p$ of that divisor, the metrics 
have cone points with
angles equal to $2\pi({\rm ord}_{\a_i}(p) + 1)$.  More importantly, the
periods of the forms are given by the Euclidean geometry of the 
developed image of the metric $ds_{\a_i}$ -- a period of a cycle $\g$ is the
distance in $\bbc$ between consecutive images of a distinguished point
in $\g$.  We reverse this procedure in \S3: we use putative
developed images of the one-forms $Gdh$, $G^{-1} dh$, and $dh$ to 
solve formally the period problem for some formal \wei data.
\enddemo

 2.2. {\it Teichm{\rm \"{\it u}}ller  theory}.
For $M$ a smooth surface, let Teich~$(M)$
denote the \tec space of all conformal structures on $M$ under the
equivalence relation given by pullback by diffeomorphisms isotopic
to the identity map id:~$M\lra M$. Then it is well-known that
Teich~$(M)$ is a smooth finite dimensional manifold if $M$ is a
closed surface.

There are two spaces of tensors on a Riemann surface $\SR$ that are
important for the \tec theory. The first is the space QD$(\SR)$ of
holomorphic quadratic differentials, i.e., tensors which have the
local form $\Phi=\vp(z)dz^2$ where $\vp(z)$ is holomorphic. The
second is the space of Beltrami differentials Belt$(\SR)$, i.e.,
tensors which have the local form $\mu=\mu(z)d\bar z/dz$.

The cotangent space $T^*_{[\SR]}$(Teich~$(M)$) is canonically
isomorphic to QD$(\SR)$, and the tangent space is given by
equivalence classes of (infinitesimal) Beltrami differentials, where
$\mu_1$ is equivalent to $\mu_2$ if
$$
\int_\SR\Phi(\mu_1-\mu_2) = 0\qquad\text{for every }\
\Phi\in\text{QD}(\SR).
$$

If $f:\bbc \to \bbc$ is a diffeomorphism, then the Beltrami
differential associated to the pullback conformal structure
is $\nu = \frac{f_{\bar z}}{f_z} \frac{d\bar z}{dz}$. If $f_{\e}$ is a
family of such diffeomorphisms with $f_0= \text{id}$, then the
infinitesimal Beltrami differential is given by
$\frac d{d\e}\bigm|_{\e=0} \nu_{f_\e} 
=\(\frac d{d\e}\bigm|_{\e=0} f_\e\)_{\bar z}$.  
We will
carry out an example of this computation in \S5.2.

A holomorphic quadratic differential comes with a picture that is
a useful aid to one's intuition about them. The picture is that of
a pair of transverse measured foliations, whose properties we sketch
briefly (see \cite{FLP} for more details).

A $C^k$ measured foliation on $\SR$ with singularities
$z_1,\dots,z_l$ of order\break $k_1,\dots,k_l$ (respectively) is given by
an open covering $\{U_i\}$ of $\SR-\{z_1,\dots,z_l\}$ and open sets
$V_1,\dots,V_l$ around $z_1,\dots,z_l$ (respectively) along with
 real valued  $C^k$ functions $v_i$ defined on $U_i$ such that
\vglue4pt
\item{(i)} $|dv_i|=|dv_j|$ on $U_i\cap U_j,$
\item{(ii)} $|dv_i|=|\imm(z-z_j)^{k_j/2}dz|$ on $U_i\cap V_ji.$
\vglue4pt

Evidently, the kernels $\ker dv_i$ define a $C^{k-1}$ line field on
$\SR$ which integrates to give a foliation $\SF$ on
$\SR-\{z_1,\dots,z_l\}$, with a $k_j+2$-pronged singularity at $z_j$.
Moreover, given an arc $A\subset\SR$, we have a well-defined measure
$\mu(A)$ given by
$$
\mu(A) = \int_A |dv|
$$
where $|dv|$ is defined by $|dv|_{U_i}=|dv_i|$. An important feature
 required
of this measure is its ``translation invariance''. That is, suppose
$A_0\subset\SR$ is an arc transverse to the foliation $\SF$,
with $\p A_0$ a pair of points, one on the leaf $l$ and one on the
leaf $l'$; then, if
we deform $A_0$ to $A_1$ via an isotopy through arcs $A_t$
that maintains the
transversality of the image of $A_0$ at every time, and also
keeps the endpoints of the arcs $A_t$ fixed on the leaves
$l$ and $l'$, respectively, then we require that
$\mu(A_0)=\mu(A_1)$.

Now a holomorphic quadratic differential $\Phi$ defines a measured
foliation in the following way. The zeros $\Phi^{-1}(0)$ of $\Phi$
are well-defined; away from these zeros, we can choose a canonical
conformal coordinate $\z(z)=\int^z\sqrt\Phi$ so that $\Phi=d\z^2$.
The local measured foliations ($\{\ree\z=\oper{const}\}$,
$|d\ree\z|$) then piece together to form a measured foliation known
as the vertical measured foliation of $\Phi$, with the
translation invariance of this measured foliation of $\Phi$
following from Cauchy's theorem.

Work of Hubbard and Masur (\cite{HM}) (see also alternate proofs
in [Ke], [Gar] and [Wo]), following Jenkins (\cite{J})
and Strebel (\cite{Str}), showed that given a measured foliation
$(\SF,\mu)$ and a Riemann surface $\SR$, there is a unique holomorphic
quadratic differential $\Phi_\mu$ on $\SR$ so that the horizontal
measured foliation of $\Phi_\mu$ is equivalent to  $(\SF,\mu)$.
\vglue12pt

{\it Extremal length\/}. The extremal length $\ext_\SR([\g])$ of a class of
arcs $\G$ on a Riemann surface $\SR$ is defined to be the conformal
invariant
$$
\ext_\SR([\g])=\sup_\rho\f{\ell^2_\rho(\G)}{\text{Area}(\rho)}
$$
where $\rho$ ranges over all conformal metrics on $\SR$ with areas
$0<\text{Area}(\rho)<\infty$ and $\ell_\rho(\G)$ denotes the infimum
of $\rho$-lengths of curves $\g\in\G$. Here $\G$ may consist of all
curves freely homotopic to a given curve, a union of free homotopy
classes, a family of arcs with endpoints in a pair of given
boundaries, or even a more general class.  Kerckhoff (\cite{Ke})
showed that this definition of extremal lengths of curves extended
naturally to a definition of extremal lengths of measured foliations.

For a class $\G$ consisting of all curves freely homotopic to a
single curve $\g\subset~M$, (or more generally, a measured foliation
$(\SF, \mu)$) we see that $\ext_{(\cd)}(\G)$ (or $\ext_{(\cd)}(\mu)$)
can be constructed as a real-valued function $\ext_{(\cd)}(\G)$:
Teich$(M)\break\lra\bbr$. Gardiner (\cite{Gar}) showed that
$\ext_{(\cd)}(\mu)$ is differentiable and Gardiner and Masur ([GM])
showed
that $\ext_{(\cd)}(\mu)\in C^1$ (Teich$(M)$).
(In our particular applications, the extremal length
functions on our moduli spaces will be real analytic; this
will be explained in \S4.5.)
Moreover Gardiner
computed that
$$
d\ext_{(\cd)}(\mu)\bigm|_{[\SR]} = 2\Phi_{\mu}
$$
so that
$$
\(d\ext_{(\cd)}(\mu)\bigm|_{[\SR]}\)[\nu] =
4\ree\int_\SR\Phi_\mu\nu.\tag2.2
$$
This formula will be crucial to our discussion in \S5.3.

\demo{{\rm 2.3.} A brief sketch of the proof}
While the details of the arguments are sometimes quite involved,
the basic logic of the approach and the ideas of the proofs are
quite simple.  Moreover, it is quite likely that they extend
to large classes of other minimal surfaces (see \cite{MW},
\cite{WW2}), so in this subsection,
we sketch the approach, as a step-by-step recipe.

\vglue8pt {\it Step\/} 1. {\it Draw the surface\/}. The first step in proving the 
existence of a minimal surface is to work out a detailed proposal.
This can either be done numerically, as in the work of
Thayer ([Th]) for the Chen-Gackstatter surfaces and Boix
and Wohlgemuth ([Bo], [Woh2]) for the low genus surfaces we treat here,
or it can be schematic, showing how various portions of the surface
might fit together.  We follow the latter approach here, which
requires quite a bit of terminology -- there are many ways to 
orient and order the handles and ends in a surface of genus five
with five planar ends and two catenoid ends, even under an
additional restriction on symmetry. To narrow the list of
possibilities, one applies either some numerical work, or some
intuition, or one attempts to continue with the outline of the
proof. We develop the terminology for this process and define
the models for some of our candidates in \S3.6, and for the 
rest of them in \S3.7.
\vglue8pt
{\it Step\/} 2. {\it Compute the divisors for the forms $Gdh$
and $G^{-1}dh$\/}. From the model that we drew in Step 1, we can
compute the divisors for the \wei data, which we just defined
to be the Gauss map $G$ and the `height' form $dh$. (Note here
how important it is that the \wei representation be given in 
terms of geometrically defined quantities --- for us, this gives the
passage between the extrinsic geometry of the minimal surface 
as defined in Step 1 and the conformal geometry and \tec
theory of the later steps.) Thus we can also compute the
divisors for the meromorphic forms $Gdh$ and $G^{-1}dh$ on the
Riemann surface (so far undetermined, but assumed to exist)
underlying the minimal surface. Of course the divisors
for a form determine the form up to a constant, so the divisor
information nearly determines the \wei data for our surface.
These \wei data are computed in \S3.6  (and \S3.7, for an
extended class) as well.
\vglue8pt
{\it Step\/} 3. {\it Compute the flat structures for the forms $Gdh$
and $G^{-1}dh$ required by the period conditions\/}. 
A meromorphic form on a Riemann surface defines
a flat singular (conformal) metric on that surface; for example,
from the form $Gdh$ on our putative Riemann surface, we determine
a line element $ds_{Gdh}= |Gdh|$. This metric is locally Euclidean
away from the support of the divisor of the form and has a 
complete  Euclidean cone structure in a neighborhood of a zero or pole of the
form.  Thus we can develop the universal cover of the surface
into the Euclidean plane.

The flat structures for the forms $Gdh$ and $G^{-1}dh$ are
not completely arbitrary; because the periods for the pair of
forms must be conjugate (formula (2.1b)), the flat structures must
develop into domains which have a particular
Euclidean geometric relationship to one another. This relationship is
crucial to our approach, so we will dwell on it somewhat. If the
map $D:\Omega \lra \bbe^2$ is the map which develops the 
flat structure of a form, say $\alpha$, on a domain 
$\Omega$ into $\bbe^2$, then the map $D$ pulls back the canonical
form $dz$ on $\bbc \cong \bbe^2$ to the form $\alpha$ on $\Omega$.
Thus the periods of $\alpha$ on the Riemann surface are given
by integrals of $dz$ along the developed image of paths in $\bbc$,
i.e. by differences of the complex numbers representing 
endpoints of those paths in $\bbc$.

We construe all of this as requiring that the flat structures 
develop into domains that are ``conjugate'': if we collect all
of the differences in positions of parallel sides for the developed image
of the form $Gdh$ into a large complex-valued 
$n$-tuple $V_{Gdh}$, and we collect
all of the differences in positions of 
corresponding parallel sides for the developed image
of the form $G^{-1}dh$ into a large 
complex-valued n-tuple $V_{G^{-1}dh}$, then these
two complex-valued vectors  $V_{Gdh}$ and $V_{G^{-1}dh}$ should
be conjugate. 
This is the flat structure implication of the period condition
formula (2.1b), here using that our situation allows that
the periods of all cycles can be found from differences
of positions of parallel sides in a flat structure.
Thus, we translate the ``period problem'' into
a statement about the Euclidean geometry of the developed
flat structures. This is done at the end of \S3.6
(and again in \S3.7, for the extended class of surfaces).

The period problem (2.1a) for the form $dh$ will be trivially
solved for the surfaces we treat here.
 
\vglue6pt
{\it Step\/} 4. {\it Define the moduli space of pairs of conjugate flat
domains\/}.  Now we work backwards.  We know the general form of 
the developed images (called $\ogup$ and $\ogdn$, respectively)
of flat structures associated to the forms
$Gdh$ and $G^{-1}dh$, but in general, there are quite a few 
parameters of the flat structures left undetermined, even after
we have assumed symmetries, determined the \wei divisor data 
for the models and used the period conditions (2.1b) to 
restrict the relative Euclidean geometries of the pair
$\ogup$ and $\ogdn$. Thus, there is a moduli space $\Delta$ of 
possible candidates of pairs $\ogup$ and $\ogdn$; our 
period problem (condition (2.1b)) is now a conformal problem
of finding such  pairs which are conformally equivalent
by a map which preserves the corresponding cone points.
(Solving this problem means that there is a well-defined
Riemann surface which can be developed into $\bbe^2$ in
two ways, so that the pair of pullbacks of the form $dz$ gives
forms $Gdh$ and $G^{-1}dh$ with conjugate periods.)

The condition of conjugacy of the domains $\ogup$ and $\ogdn$
often dictates some restrictions on the moduli space, and even
a collection of geometrically defined coordinates.  We work these
out in the \S4.2 and \S4.4. 

When these moduli spaces are empty,  we have a proof of
nonexistence: this is the case for the surfaces studied
in part (ii) of Theorem B -- see \S7.

\vglue12pt
{\it Step\/} 5. {\it Solve the conformal problem using Teichm{\rm \'{\it u}}ller theory\/}.
At this juncture, our minimal surface problem has become a 
problem in finding a special point in a product of moduli
spaces of complex domains; we will have no further references
to minimal surface theory. The plan is straightforward: we will
define a height function $\height: \Delta \lra \bbr$
with the properties:
\vglue4pt
\item{1}  (Reflexivity) The height $\height$ equals $0$ only at a solution
to the conformal problem. 
\item{2}  (Properness) The height $\height$ is proper on $\Delta$.
This ensures the existence of a critical point.
\item{3}  (Noncritical Flow) If the height $\height$ at a pair $(\ogup, \ogdn)$
does not vanish, then the height $\height$ is not critical at 
that pair.
\vglue4pt

This is clearly enough to solve the problem: we now sketch the
proofs of these steps.
\vglue12pt
{\it Step\/} 5a. {\it Reflexivity\/}. We need conformal invariants of a 
domain that provide a complete set of invariants (for Reflexivity),
have estimable asymptotics (for Properness), and computable 
first derivatives in moduli space (for the Noncritical Flow Property).
One obvious choice is a set of functions of 
extremal lengths for a good choice
of curve systems, say $\Gamma = \{\gamma_1, \dots, \gamma_K\}$ 
on the domains; it is sufficient for Reflexivity that the 
extremal lengths of these curves provide coordinates for $\Delta$.
Suitable curve systems
are defined for our examples in \S4.2 and \S4.4.
We then define a height function $\height$ which vanishes only when there
is agreement between all of the extremal lengths
$\text{Ext}_{\ogup}(\gamma_i)= \text{Ext}_{\ogdn}(\gamma_i)$ and which 
blows up when $\text{Ext}_{\ogup}(\gamma_i)$ and $\text{Ext}_{\ogdn}(\gamma_i)$
either decay or blow up at different rates. See for example
Definition 4.3.1 and Lemma~4.5.5.
\vglue12pt
{\it Step\/}  5b. {\it Properness\/}. Our height function will measure 
differences in the extremal lengths $\text{Ext}_{\ogup}(\gamma_i)$ and
$\text{Ext}_{\ogdn}(\gamma_i)$. Often, but not always, a geometric
degeneration of the flat structure of either $\ogup$ or
$\ogdn$ will force one of the extremal lengths $\text{Ext}_{\cdot}(\gamma_i)$
to tend to zero or infinity, while the other 
extremal length stays finite and
bounded away from zero.  This is a straightforward situation where
it will be obvious that the height function will blow up.
A more subtle case arises when a geometric degeneration of the
flat structure forces {\it both} of the extremal lengths
$\text{Ext}_{\ogup}(\gamma_i)$ and $\text{Ext}_{\ogdn}(\gamma_i)$ to
decay {\it simultaneously\/}  (or explode). In that case, we begin
by observing that there is a natural map between the
vector $\langle\text{Ext}_{\ogup}(\gamma_i)\rangle$ and the vector
$\langle\text{Ext}_{\ogdn}(\gamma_i)\rangle$. This pair of vectors is 
reminiscent of pairs of 
solutions to a hypergeometric differential
equation, and we show, by a monodromy
argument analogous to that used in the study of those equations, that
it is not possible for corresponding components of that vector
to vanish or blow up at identical rates. In particular, we show that the
logarithmic terms in the asymptotic expansion of the extremal lengths
near zero have a different sign, and this sign difference forces a
difference in the rates of decay that is detected by the height
function, forcing it to blow up in this case.  The monodromy argument
is given in \S4.7, and the properness discussion consumes \S4.

\vglue12pt
{\it Step\/} 5c. {\it Noncritical flow\/}. The domains $\ogup$ and $\ogdn$
have a remarkable property: if
$\text{Ext}_{\ogup}(\gamma_i) > \text{Ext}_{\ogdn}(\gamma_i)$, then 
there are always deformations so that when 
we deform $\ogup$ so as to  decrease $\text{Ext}_{\ogup}(\gamma_i)$,
the conjugacy condition forces us to deform $\ogdn$ so as
to increase $\text{Ext}_{\ogdn}(\gamma_i)$.  We can thus always deform
$\ogup$ and $\ogdn$ in order  to reduce one term of the height
function $\height$. We develop this step in \S5.
\vglue12pt
{\it Step\/} 5d.\ {\it  Regeneration\/}. In the process described in 
the previous step, an issue arises: we might be able to 
reduce one term of the height function via a deformation, but
this might affect the other terms, so as not to  provide
an overall decrease in height.  We thus seek a locus $\SY$ in
our moduli space where the height function has but a single
nonvanishing term, and all the other terms vanish to at least
second order.  If we can find such a locus $\SY$, we can flow along
that locus to a solution.  To begin our search for such a locus,
we observe which flat domains arise as limits of our
domains $\ogup$ and $\ogdn$: commonly, the degenerate domains 
are the flat domains for a similar minimal surface problem,
maybe of slightly lower genus or a few fewer ends.

We find our desired locus
by considering the boundary $\p\overline{\Delta}$ of the closure 
$\overline{\Delta}$ of the moduli
space $\Delta$; this boundary has strata of moduli spaces $\Delta'$
for minimal surface problems of lower complexity.  By induction,
there is a solution $X'$ of those problems represented on such
a boundary strata $\Delta'$ (with all of the
corresponding extremal lengths in agreement), 
and we prove that there is a 
locus $\SY \subset \Delta$ inside the larger moduli space $\Delta$ 
(with $\SY$ limiting on $X'$) which 
has the analogues of those same extremal lengths in agreement.
As a corollary of that condition, the height function on $\SY$
has the desired simple properties.
This is developed in \S6.

The proof may be summarized as follows: we restrict to a locus
$\SY$ on which we have the Noncritical Flow Property. By step 5b,
the height function $\SH$ is proper on $\SY$; thus there is 
a critical point $X$ on $\SY$ for $\SH$.  The Noncritical Flow
then forces $\SH(X)=0$, so by Reflexivity, the surface 
represented by $X$ is a solution to our conformal problem, and
hence also defines a solution to the minimal surface problem.
\enddemo

\section{The geometry of orthodisks}

In this section we introduce the notion of orthodisks.

\demo{{\rm 3.1.} Orthodisks}
Consider the upper half-plane and $n\ge3$ distinguished points
$t_i$ on the real line. The point $t_\infty=\infty$ will also
be a distinguished point. We will refer to the upper half-plane together
with these data as a {\it conformal polygon}
and to the distinguished points as
{\it vertices}. Two conformal polygons are {\it conformally equivalent} if
there is
a biholomorphic map between the disks carrying vertices to vertices,
and fixing $\infty$.

Let
$a_i$ be some odd integers such that
$$a_\infty=-4-\sum_i a_i. \tag 3.1$$

By a Schwarz-Christoffel
map we mean the map
$$F:z \mapsto \int_i^z (t-t_1)^{a_1/2}\cdot\ldots\cdot(t-t_n)^{a_n/2}
dt. 
\tag 3.2$$
A point $t_i$ with $a_i>-2$ is called {\it finite}, otherwise {\it infinite}.
By (3.1), there is at least one finite vertex (possibly at $t_{\infty}$).

\demo{Definition {\rm 3.1.1}} Let $a_i$ be odd integers.
The pull-back of the flat metric on $\bbc$ by $F$ defines a
complete flat metric with boundary on $\Bbb H \cup \Bbb R$  without the
infinite vertices. We call such a metric an {\it orthodisk}. The $a_i$ are
called the {\it vertex data} of the orthodisk. The {\it edges} of an orthodisk
are the boundary segments between vertices; they come in a natural order.
Consecutive edges meet orthogonally at the finite vertices.
Every other edge is parallel
under the parallelism induced by the flat metric of the
orthodisk. Oriented distances between parallel edges are called {\it periods}.
The periods can have
four  different signs:
$+1,-1,+i,-i$.
\enddemo

{\it Remark}. The integer $a_i$ corresponds to an angle $(a_i+2)\pi/2$
of the orthodisk. Negative angles are meaningful because 
a vertex (with a negative angle $-\theta$) lies at infinity
and is the intersection of a pair of lines which also
intersect at a finite point, where they make a {\it positive} angle of
$+\theta$.
\vglue12pt

In all the drawings of the orthodisks to follow, we mean the domain to be
to the left
of the boundary.

\vglue8pt{\it Example {\rm 3.1.2}}. This is conceivably the simplest orthodisk, bounding the
second quadrant in $\Bbb R^2$:
\enddemo
\figin{Zwedge}{500}

\centerline{Simple orthodisk}

\vglue12pt {\it Example} 3.1.3. Here is an orthodisk with a branch point which is
drawn fat. The disk consists of the region northwest of the 
larger $5$-vertex boundary
{\it and} the region northwest of the fat vertex boundary.

\figin{Zbranched}{500}

\centerline{Branched orthodisk $\ldots$}

\vglue6pt

To get a clear picture of this orthodisk, we can glue the following
two domains together along the fat diagonal boundary line, in
exactly the way that a Riemann surface is assembled from several
sheets and a branch cut.

\figin{Zbranched2}{500}

\centerline{$\ldots$ decomposed into two pieces}
\vglue12pt

Denote by $\gamma_i$ an  oriented curve
connecting edge
$t_i t_{i+1}$ with edge $t_{i+2} t_{i+3}$. There are $n-1$ such
curves.  We will denote by the same name, also, their homotopy classes. It is
well-known that the extremal lengths
$\ext (\gamma_i)$ determine the conformal structure of the conformal
polygon. This follows because the extremal lengths of a conformal
quadrilateral determine the cross-ratios of the vertices of the 
quadrilateral, and the cross-ratios then determine the vertices
of the polygon, where we view the vertices as distinguished points
on the boundary of, say, the disk.

Each orthodisk has a natural conformal structure and hence determines a
conformal polygon. Vice versa, given a conformal polygon and applying
the Schwarz-Christoffel map to it, one obtains an orthodisk with certain
periods. By scaling, one can arrange for one period to be (say) $1$.
We then call this map a normalized period map.

Applying the Schwarz-Christoffel discussion above, and then
using the implicit function theorem, one can prove the following proposition
(which we will never use, so we omit a detailed proof):

\nonumproclaim{Lemma 3.1.4} The normalized  map from the space of conformal
polygons to periods of orthodisks of given vertex data is locally
injective.
\endproclaim

We will restrict our attention to a subclass of orthodisks which have a
real symmetry.

\demo{Definition {\rm 3.1.5}} An orthodisk is called {\it symmetric} if it has
a reflectional symmetry which fixes two vertices.
\enddemo

All of the orthodisks under 
discussion have angles which are odd multiples
of $\pi/2$, and, with but two exceptional cases of the vertices
representing the catenoid ends,  the angles alternate between being
congruent to $\pi/2 $ and congruent to $3\pi/2$ modulo $2\pi$; the 
result of this combinatorics of angles is that any symmetric orthodisk
must fix a pair of vertices that are on `opposite' sides of the 
orthodisks, each halfway between the vertices representing the
catenoid ends.

The above Examples 3.1.2 and 
3.1.3 are symmetric. As a convention we will draw all
orthodisks such that the symmetry line is the diagonal $y=-x$.

\demo{{\rm 3.2.}  From orthodisks to Riemann surfaces}
We begin with an orthodisk $X$ and we
describe in this section a method
to canonically construct a hyperelliptic Riemann
surface $Z$ with a meromorphic $1$-form $\om_X$ from X.

Let $X$ be a conformal polygon. First double $X$ along its boundary;
i.e., take $X$ and a copy $\overline{X}$ (with the opposite orientation)
and glue $X$ and $\overline{X}$ together
along their commonly labelled edges. 
The resulting complex space $Y$ is then topologically
a sphere with distinguished points which we 
also call $t_i$;
moreover, the sphere $Y$ (punctured at the infinite vertices)
inherits a complete singular flat structure from 
$X$ and $\overline{X}$. To see this last statement,
begin by observing that the flat structure 
on $(\text{int}\ X) \cup \text{int}(\overline{X})$ extends across the interior of the 
boundary segments of $\p X$ and $\p\overline{X}$ because those
segments are straight. At the finite vertices, the flat 
structure has cone-like singularities, while the 
infinite vertices are at infinite distance from any point
in a neighborhood of them. More precisely,
we compute that at the finite vertices,
the cone angle at the vertex of the sphere $Y$ corresponding to $t_i$
will be $\alpha_i=2\pi (a_i/2+1)$.
We use formally the same formula to attach
(negative) cone angles to the vertices at infinity. 
In fact the end corresponding
to a vertex $t_i$ at infinity in $\bbc$
is conformally the same as the infinite end of a
cone with cone angle $-\alpha_i$.

Next, construct simple closed curves 
$c_i=\gamma_i \cup \overline{\gamma_i}$ 
on $Y$ by connecting $\gamma_i$ on $X$
with the oppositely oriented arc $\overline{\gamma_i}$
on $\overline{X}$. Now, as $Y$ is a double of $X$,
the extremal lengths of the
$c_i$ determine the conformal structure of the punctured sphere $Y$.
Finally, 
construct the hyperelliptic double cover $Z$ over $Y$ branched over the
vertices $t_i$. The conformal structure of this covering is independent of the
choice of branching  slits.
However, this is not true for its \tec class:
the choice of slits affects the marking of the surface, i.e., the
choice of a basis for $\pi_1Z$.
Because we soon want to measure
distances in
\tec space between points constructed this way, 
we need to choose the slits consistently:
we take as
branch  slits the odd numbered edges.
It is then easy to check that
$c_i$ has two lifts to
$Z$ whose sum is null homologous. To pick a specific lift, we make, once
and for all, a choice: choose, for each $i$, a lift $\tilde c_i$ of
$c_i$. Subsequent constructions will depend on that choice, so we will have
to ensure
that the final statements are independent of this choice.

Metrically, on the double cover $Z$ 
we obtain a lifted flat cone metric with cone angles
$\tilde\alpha_i=2\pi (a_i+2)$ at the lifts of 
$t_i$: this angle is an odd  multiple of $2\pi$.
Furthermore, because only parallel edges 
of $X$ and $\overline{X}$ are glued together in the
whole construction, the cone metric on $Z$ has no linear holonomy.
Hence the exterior derivative of the (locally defined) developing map of
the flat structure
defines a global holomorphic $1$-form $\omega_X$ on $Z$.
(In slightly different language: we develop the 
flat structure on $Z$ into the Euclidean plane $\bbe^2$, on which
there is defined the one form $dz$, in the usual coordinates. Because
the flat structures on the orthodisks $X$ and $\overline{X}$
are bounded by horizontal and vertical lines, and our identifications
take horizontal edges to horizontal edges and vertical edges to
vertical edges, the transition maps between
the distinct developed images of, say, $X$, are given by
translations.  Thus the
form $dz$ pulls back to a well-defined one-form $\om_X$ on $Z$.)
This form $\om_X$ has  
zeroes at the lifts of $t_i$ of order $a_i+1$; here,
negative orders correspond to poles. Because the developing map is only
defined up to a complex linear transformation, the $1$-form $\omega_X$ will
only
be defined up to homothety.

\nonumproclaim{Lemma 3.2.1} Up to a  factor independent of $i${\rm ,} the
period of $\omega_X$ along
$\tilde c_i$ is  the period of $\gamma_i$.
\endproclaim

\demo{Proof} Note that the developing map of a flat metric is (locally)
only well-defined up to post-composition
with a complex linear transformation, so that the $1$-form $\omega$ is
well-defined up to multiplication 
by some nonvanishing complex number. After having
made that choice, a period of $\omega$ along $\tilde c_i$
is by construction the difference between the image 
(under the developing map defining
this choice of $\omega$) 
of the endpoint of $c_i$
and the image of its initial point;
this is because integrals of $\om_X$ along paths in $Z$ 
push down via the developing map to integrals of $dz$ along
pushed-forward paths in $\bbe^2 \cong \bbc$. 
\enddemo

 3.3. {\it From orthodisks to \wei data}.
We describe how pairs of orthodisks can be used to write down {\it formal}
\wei data
for minimal surfaces.  This procedure reverses the construction of 
\S2.1.1 where we find the flat structures associated to a 
minimal surface; here we will begin to solve the period problem
for a minimal surface by first specifying pairs of relevant
flat structures (formally $\ogup$ and $\ogdn$)
whose geometry represents compatible periods, in the sense of 
\S2.1.1.

\demo{Definition {\rm  3.3.1}}
Two  orthodisks $X_1$ and $X_2$ on the same underlying conformal polygon
but with different angles (exponents) are called {\it conformal}. By the
construction above,
they give rise to two
distinct
 meromorphic $1$-forms $\omega_1, \omega_2$ on the same Riemann surface $Z$.
\enddemo

{\it Definition} 3.3.2.
Two  orthodisks $X_1$ and $X_2$ with different vertex data are called
{\it conjugate} if 
there is a line $l \subset \bbc$ so that 
corresponding periods are symmetric with respect to that line.

\demo{Definition {\rm 3.3.3}}
Two  orthodisks $X_1$ and $X_2$ are called {\it reflexive} if they are
conformal and
conjugate.
\enddemo

Now suppose that the pair of orthodisks 
$X_1$ and $X_2$ are defined on the same conformal
polygon. Let the corresponding
$1$-forms
on $Z$ be denoted $\omega_1$ and $\omega_2$. 
We want to find a meromorphic function $G$
(the Gauss map) and a meromorphic $1$-form $dh$ (the height differential)
such that
$$\align
\omega_1 &= G dh, \\
\omega_2 &= G^{-1} dh
\endalign $$
and such that $G$ and $dh$ are \wei data of a minimal surface on $Z$.

To solve the above equations, denote the vertex data of $X_1$ by $a_i$ and the
vertex data of $X_2$ by $b_i$ (recall that these are odd numbers). Then at
a point $t_i$ on the putative underlying Riemann surface, the meromorphic
quadratic differential $\omega_1\omega_2$ has a zero of order $a_i+b_i+2$
which is an even number; thus, we first ask for $dh$ with zero of order
$\f12(a_i+b_i)+1$.

Here is a simply defined but quite general 
case where one can find such a $dh$:

\nonumproclaim{Lemma 3.3.4} Suppose that{\rm ,} for each index $i${\rm ,} 
the sum $a_i+b_i\equiv 0 \pmod 2$. Then
there is a meromorphic $1$-form $dh$ on $Z$ as above with all
periods purely imaginary.
\endproclaim

\demo{{P}roof}
Consider, on the sphere $Y$, the meromorphic $1$-form
$$\omega=\prod_{i=1}^n(t-t_i)^{(a_i+b_i)/2} dt .$$
Lift it via the hyperelliptic covering projection to $Z$ and call the
lift $dh$. Clearly, this meromorphic $1$-form has the required
zeroes (and poles) at the preimages of the $t_i$.

The periods of $dh$ are all computable as periods
of  the form $\omega$ on $Y$ which are all computable as $2\pi i$
multiples of purely real residues.  
\enddemo

This lemma will apply to all cases needed.
We deduce:

\nonumproclaim{Theorem 3.3.5} Let $X_1$ and $X_2$ be reflexive
orthodisks with exponents $a_i$ and $b_i$. Suppose that
$a_i+b_i\equiv 0 \pmod 2$. Then the above \wei data define a minimal
surface.
\endproclaim

\demo{Proof} It remains to show that the \wei data thus found have purely
imaginary periods, or equivalently (see formula (2.2)) that the
forms $Gdh$ and $G^{-1} dh$ have conjugate periods.  Yet the forms
$Gdh$ and $G^{-1} dh$ on the Riemann surface $Z$ are lifted from the
canonical forms on the orthodisks $X_1$ and $X_2$, respectively;
the conjugacy of the orthodisks $X_1$ and $X_2$ then forces the
periods of $Gdh$ and $G^{-1} dh$ on the Riemann surface $Z$ to be
conjugate.
\enddemo

{\it Definition} 3.3.6. We call the pair of vertex data $a_i,
b_i$ of a pair of orthodisks of the same genus {\it formal \wei
data}.

\demo{{\rm 3.4.} Geometric significance of the formal \wei data}
In this section we assume that we have a pair of reflexive orthodisks 
(i.e., conjugate orthodisks defined over the same underlying
conformal ploygon) and we then 
determine the geometric data
of the resulting minimal surface. These data will be given in terms of
the formal \wei   divisor data
 $a_j$ and $b_j$.

First, by construction,
the Riemann surface resulting  from a pair of conformal orthodisks with
$n=2g+1$
 vertices will have genus
$g$.

Recall (\S2.2) that
the Riemannian metric of the associated  minimal surface with
\wei data $G$ and $dh$
is given by
$$ds=(\vert G \vert +\frac1{\vert G\vert})\vert dh\vert.$$

Singularities of this metric can only occur at the $2g+2$ vertices of
the conformal polygon, where the branching (of the surface $Z$ over the
double $Y$) occurs. By the divisor discussion, the one-form
$dh$ will have a zero of order
$1+\frac{a_i+b_i}2$ and $G$ a zero of order $\frac{a_i-b_i}2$. Hence
$\vert G \vert +\frac1{\vert G\vert}$ will have a pole of order
$\frac{\vert a_i-b_i \vert}2$ and thus
to have a complete metric it is necessary and
sufficient that
$$2+a_i+b_i \le \vert a_i-b_i \vert$$
for all $i$ (including the vertex at infinity!).

It is also easy to compute the total curvature or equivalently the degree
of the
Gauss map as
$$\deg G = \frac12\sum_{a_i>b_i}a_i-b_i$$
where the index $i=\infty$ is possibly included.

In the following table
we compile a list of typical special points on minimal surfaces
and the resulting vertex data. Formal \wei data are then given 
as ordered lists of points of these types, and these types only.
In particular, 
these special points appear as vertices on the conformal
polygons, with vertex data given as fractions of cone angles
for the flat structures of the forms $(Gdh, G^{-1}dh, dh)$.
The fractions arise because a neighborhood of a special point
on the Riemann surface may $n$-fold cover a neighborhood of its image
on the conformal polygon, so that the cone angle on the
conformal polygon will be $\frac1n$ of the cone angle on the
Riemann surface. Formal \wei data are conformal polygons labelled
with vertex data chosen from this table.

Denote by $H$ finite nonbranch points of the minimal surface where
the symmetry curves intersect, by $C$ a catenoid end, by $P$ a planar end
with degree of the Gauss map equal to $3$, and with  $R$ a point which lies on just
one symmetry arc. An up-arrow ($\uparrow$) means that the Gauss map takes
the value
$\infty$, a down-arrow ($\downarrow$)  that the Gauss map takes the value
$0$. Then the divisors at the respective points are given by the table
below:
$$\matrix
 & G & dh & G dh & G^{-1} dh & \angle G dh & \angle G^{-1} dh & ex(G dh) &
ex(G^{-1} dh)\\
H\uparrow & -1 & 1 & 0 & 2 & \pi/2 & 3\pi/2 & -1 & 1\\
H\downarrow & 1 & 1 & 2 & 0 & 3\pi/2 & \pi/2 & 1 & -1\\
C\uparrow & -1 & -1 & -2 & 0 & -\pi/2 & \pi/2 & -3 & -1\\
C\downarrow & 1 & -1 & 0 & -2 & \pi/2 & -\pi/2 & -1 & -3\\
P\uparrow & -3 & 1 & -2 & 4 & -\pi/2 & 5\pi/2 & -3 & 3\\
P\downarrow & 3 & 1 & 4 & -2 & 5\pi/2 & -\pi/2 &3 & -3\\
R\uparrow & -1 & 1 & 0 & 2 & \pi & 3\pi & 0 & 2 &\\
R\downarrow & 1 & 1 & 2 & 0 & 3\pi & \pi & 2 & 0 &
\endmatrix
$$
\vglue8pt

The last two
columns just contain the vertex data $a_i$ and $b_i$ of the
formal \wei data
for $Gdh$ and $G^{-1}dh$. The last two rows record angles that are not quite
in the same pattern as the previous rows; the points $R$ have a 
surrounding neighborhood with only  a two-fold symmetry while the 
points $H$, $C$, and $P$ all have neighborhoods with a four-fold
symmetry. These points will only be important in \S8.2 when
we indicate how one might 
deform the surfaces with planar ends into those with catenoid ends.
(Full details will be given in \cite{WW2}.)

\demo{{\rm 3.5.} Examples of simple orthodisks}
In this section we give some motivating examples of minimal surfaces
and their corresponding orthodisks.  We begin first with some
discussion about moduli spaces of orthodisks, which we shall later
formalize at the outset of \S4. We then organize our examples
into families, depending upon whether the associated moduli spaces
are empty, singletons, or nontrivial.

To begin then, observe that from some formal \wei data, we may
draw orthodisk systems consisting of conjugate orthodisks
$\ogup$
and $\ogdn$ for 
the forms $Gdh$ and ${G^{-1}dh}$, respectively. We often have some 
freedom in deciding the lengths of the edges of the orthodisks,
even up to equivalence of conformal polygons.  Thus, the 
formal \wei data determine a moduli space of conjugate
orthodisk systems $\{\ogup, \ogdn \}$.
While the function theory of this moduli space is the principal subject 
of this paper, here we shall content ourselves with a few examples
for which the moduli space is either trivial or small.
(See \cite{WW2} for more examples.)
\enddemo

 {3.5.1}. {\it A singleton moduli space\/}.

\demo{Example {\rm  3.5.1}} The Costa surface ${\rm DH}_{0,0} \ (= CT_0)$.
This surface is the starting point for the investigations
in this paper. We seek a torus with two catenoid ends
($C_1$ and $C_2$), one planar
end $P_1$ and one finite point $H_1$ with vertical normal.
We also assume our standard eight symmetries: reflections about 
two orthogonal vertical planes and a rotation about a 
horizontal line. (Previous existence and uniqueness proofs
may be found in \cite{Cos} and \cite{Ho-Me}.)
From the table in \S3.4, we can write formal \wei data
for this surface, assuming that the points occur on the
real line in the order $C_1 - P_1 - C_2\break - H_1$. We will 
justify this assumption in \S3.6; however, below
we reproduce a computer image of a fundamental domain
of Costa's surface for the group generated by reflections in
vertical planes. This makes apparent the conformal polygon,
and the order of the special points on its boundary. 

A table corresponding to the one in \S3.4 is given below; it 
computes the cone angles for the forms $Gdh$ and $G^{-1}dh$ 
on the putative conformal polygon for Costa's surface.  These flat
structures are drawn below, assuming an additional symmetry about the
line $\{y=-x\}$. 
$$\matrix
 & G & dh & G dh & G^{-1} dh & \angle G dh & \angle G^{-1} dh & ex(G dh) &
ex(G^{-1} dh)\\
C_1\downarrow & 1 & -1 & 0 & -2 & \pi/2 & -\pi/2 & -1 & -3\\
P_1\uparrow & -3 & 1 & -2 & 4 & -\pi/2 & 5\pi/2 & -3 & 3\\
C_2\downarrow & 1 & -1 & 0 & -2 & \pi/2 & -\pi/2 & -1 & -3\\
H_1\downarrow & 1 & 1 & 2 & 0 & 3\pi/2 & \pi/2 & 1 & -1
\endmatrix
$$

Yet recall that
any conformal
quadrilateral with a symmetry across a 
diagonal is conformally equivalent to a square.

\figin{Zcostapatch}{500}

\centerline{A quarter of Costa's surface --- a conformal image of the orthodisks}
\vglue12pt

\figin{Zu}{500}

\centerline{Costa orthodisks}
\vglue12pt

Thus we then note that (as in Example
3.5.1), the moduli space of possible examples consists only of the 
singleton of a pair of square tori, so that the only element in
the moduli space is a reflexive pair. This establishes the
existence of this surface, by a proof that is somewhat distinct
from the other proofs of existence of this surface 
(\cite{Cos}, \cite{Ho-Ka}).

\vglue12pt{3.5.2.} {\it Empty moduli spaces\/}.

\demo{Nonexample {\rm 3.5.2(i)}} Catenoid with one handle.
In this example, we try to construct a minimal surface with two catenoid ends
and one handle. Now, a theorem of Schoen \cite{Sch} implies that 
any minimal surface of genus $g \ge 1$
with but two embedded catenoid ends cannot exist;
here we will require in addition
the eightfold symmetry present in all of our
examples. Yet, it is clear from an analysis of the
orthodisks pictured below that such a surface cannot
exist with $4$-fold symmetry, because there is no symmetric and
conjugate pair: the periods from $\overline{C_1C_2}$ to
$\overline{H_1H_2}$ are conjugate about the line 
$\{y=-x\}$, while the periods from $\overline{C_2H_2}$ to
$\overline{H_1C_1}$ are conjugate across the line
$\{y=-x\}$.

\figin{Zcathand}{500}

\centerline{Catenoid with one handle}
\vglue12pt

We will use this technique of showing nonexistence extensively
when we prove the nonexistence parts of the main theorem.
 \enddemo

{\it Nonexample} 3.5.2(ii). The Horgan surface ${\rm DH}_{0,1}$.
The second example is called the Horgan surface (see \cite{Ho-Ka}): To
visualize it,
start with one plane and two handles growing upward and downward and
perpendicular
to each other.
Both handles connect to catenoid ends.
This looks almost as shown in the following figure:

\figin{Zhorganbw}{600}
\vglue-12pt
\centerline{The Horgan surface?}
\vglue12pt

This pattern leads to a sequence
$$C_1 - P - C_2 -  H_2 - M - H_1 - C_1$$
(see \S3.6)
where $M$ denotes a previously unencountered 
regular point where the symmetry lines
cross. The  orthodisks  are as follows:

\centerline{\BoxedEPSF{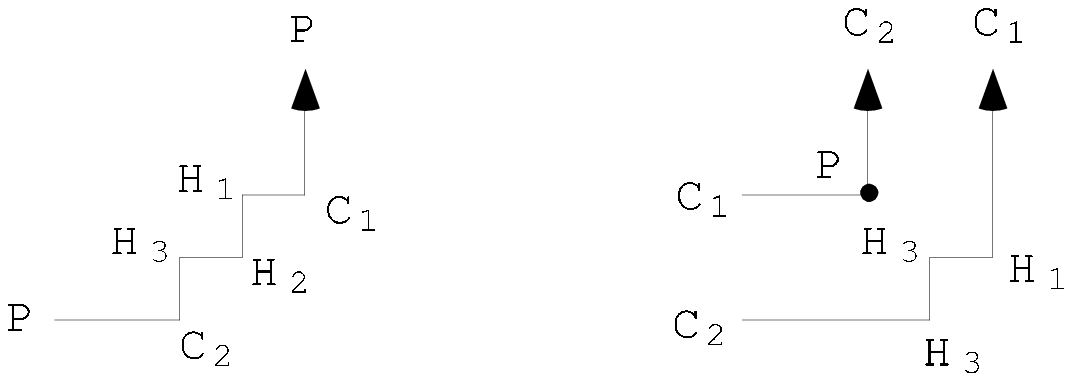 scaled 600}}

\vglue8pt
\centerline{Horgan orthodisk}
\vglue8pt

Recently, one of us \cite{W} showed that this surface cannot exist
by a method that is different from that used in Example 3.5.2(ii).
Therefore, this moduli space is empty.
We find it remarkable that our general existence proof fails at precisely one 
point for the case of this surface; see \S4, after Definition~4.1.1.

\demo{{\rm 3.6.} Orthodisks for the  Costa towers}
This section is intended to give a heuristical picture for the Costa towers and to derive the
orthodisk pairs from this picture.

As a first series of new examples, we introduce the orthodisks for
generalized Costa surfaces
which we call {\it Costa towers\/}:
\enddemo

{\it Definition} 3.6.1. A {\it Costa tower} of genus $g$ is a
complete minimal surface $CT_g$ 
of genus $g$ with two parallel catenoid ends and $g$ planar
ends which are all embedded. These surfaces have total curvature $4\pi
(g+2)$.

\demo{{R}emark {\rm 3.6.2}} The existence of generalized Costa surfaces is known for
$g\le 2$ and numerical evidence has been provided for $g=3$; see
\cite{Cos}, \cite{Bo}, \cite{Ho-Me}, \cite{Woh1}, \cite{Woh2}. All these
examples have (at least) a two-fold reflectional symmetry at two
perpendicular planes intersecting in the
$z$-axis. In their construction, this assumption makes the period problem
low dimensional. The genus $1$ Costa surface has the square torus 
as underlying Riemann surface and is known to admit
a deformation to all rectangular tori such that the planar end becomes
catenoidal.

We will briefly discuss these deformations from our point of view in
\S8.2.
\enddemo

To get an impression of how these surfaces might appear, imagine cutting a
catenoid by $g$ horizontal planes. For each plane, we have to resolve a
singularity which looks topologically like a horizontal plane cut by a vertical
cylinder. 
We replace a neighborhood of the singular locus in the complex
of cylinder and plane by one of the two shapes below; i.e.,
we replace the upper part of the cylinder by a $Y$-piece
whose lower two boundary components are glued to two holes in the plane and 
glue the boundary component of the lower part of the cylinder to another hole
in the plane in between the two first holes, and finally glue the 
boundary of the flat planar piece to the planar boundary of the
excised neighborhood.

\centerline{\BoxedEPSF{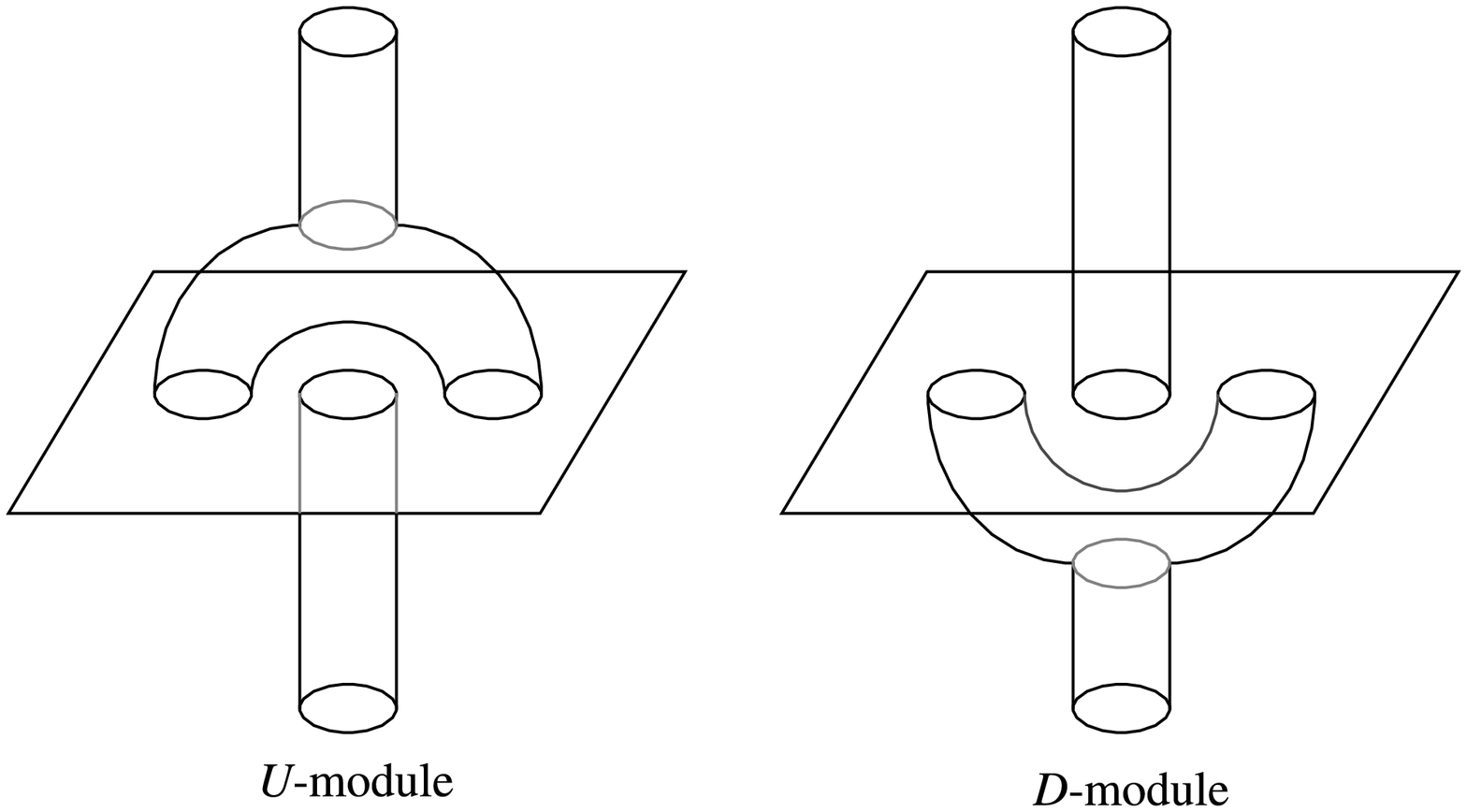 scaled 425}}
\vglue-36pt
\centerline{Building blocks for informal surface models}
\vglue4pt

We will call the left gadget a $U$-module (for {\it up}) and the right one a
$D$-module (for {\it down}). It is convenient to regard the planar
portions
of the figures as extending indefinitely.

Stacking $g$
copies of these modules over each other and finishing with catenoidal ends
yields a rough informal model of  candidates for a Costa tower of genus $g$.
These models are not very realistic, but they will suffice to write
down candidate formal \wei data for the surfaces.
\vglue6pt
{\it Remark {\rm 3.6.3}}. At first glance there are two more
possible constructions of this type, which are  obtained by rotating the
above pieces by $90^\circ$ around the
$z$-axis, but it turns out that these give no new \wei data, so we will
neglect them from
the beginning. 
\vglue6pt 

Despite this, there are still many possible ways to attach these modules
to each other, but we will not attempt to discuss all of them. Only for a
very distinctive construction will our proof establish the existence of minimal
surfaces corresponding to these models, and we have presently no contribution
to make  regarding the others. For these other models, no complete
existence or nonexistence proofs are known, but numerical experiments make it
doubtful that other module towers will produce more minimal surfaces.

Now we proceed to derive the \wei candidate data defined by a module tower.
Cutting a generalized Costa surface composed of these modules by its two
vertical symmetry planes
decomposes it into four congruent simply connected domains which we
want to describe using orthodisks. The first step is to recover  the
formal \wei data from the putative geometry:

\vglue6pt {\elevensc Lemma 3.6.4}. {\it Let $X$ be a generalized Costa surface consisting of
$g$ modules $X_i$ of type $U$. Denote by $P_i$ the point on $X$
corresponding to the planar end of module $i$ and by $H_i$ the saddle point
of that module $X_i$. 
Denote by $G$ the Gauss map and by $dh$ the height differential of
$X$. Then the divisors of $G$ and $dh$ are given by the table below}\/:

$$\align \noalign{\vskip-24pt}
(G) &=C_1^{1}P_1^{-1} H_1^{1}P_2^{1}H_2^{-1}\cdots P_g^{(-1)^g}H_g^{-(-1)^g}
C_2^{-(-1)^g} \\
\noalign{\vskip5pt}
(dh) &= C_1^{-1}P_1^{1} H_1^{1}P_2^{1}H_2^{1}\cdots P_g^{1}H_g^{1}
C_2^{-1} \\ \noalign{\vskip5pt}
(Gdh) &= H_1^2 P_2^2 \cdots \\
\noalign{\vskip5pt}
(G^{-1}dh) &= C_1^{-2} P_1^2 H_2^2 \cdots\thinspace .
\endalign$$
 
\demo{Proof} First we determine the direction of the normal vector
at all of the vertical points of $G$. We can
assume that the outer normal vector at the upper catenoid end $C_1$ is pointing
downward. As can be read off from the above figure, the normal will switch
between up and down
from planar end to planar end. Hence all odd labeled planar ends will have an
upward pointing normal, all even planar ends a downward pointing normal, and
the bottom catenoid end has normal pointing up if
and only if $g$ is even. Again by the figure above, the normal at the
saddle point of a module point is always in the opposite
direction as the normal at the planar end.
Using the information from the table at the end of \S3.4, we get the
claimed divisors.
\enddemo

{\it Definition} 3.6.5. Let $T$ be a sequence of length $g$ 
consisting of
symbols $U$ and $D$. A surface of type $T$ is a generalized Costa surface
consisting of modules $X_i$ where the handle of module $i$ grows upward
(resp.\
downward) if the symbol $i$ is $U$ (resp. $D$).

\nonumproclaim{Lemma 3.6.6} Let $X$ be a surface of type $T$. Then the formal
\wei data are
explicitly determined by the type $T$.
\endproclaim

{\it Proof}.
To see this, we first consider a surface where all modules are $U$-modules.
We follow the symmetry line given by the intersection with the $y=0$-plane,
beginning at the top catenoid end on the left. This line goes down to the
first planar end to the left. There the total angle is divided into four equal
pieces by the two symmetry lines meeting there. We continue now on the
symmetry line defined by the $x=0$-plane which goes down into the second
planar end. Continuing this process and switching to another symmetry line at
each vertex, we descend through all planar ends until we reach the bottom
catenoid end. From there, we ascend through all handles back to the top
catenoid end. We denote this closed path by
$$C_1 \to P_1 \to P_2 \to \ldots \to P_g\to C_2\to H_g\to \ldots\to
H_2\to H_1\to C_1.$$ Now exchange the $k^{{\rm{th}}}$ module by a
$D$-module.
As can be seen from the figure above, this affects only the entries $P_k$
and
$H_k$ in the path list; they are just \pagebreak exchanged. Hence for an arbitrary
surface
$X$ of type
$D$ we get a path sequence
$$C_1 \to A_1 \to A_2 \to \ldots \to A_g \to C_2 \to B_g \to \ldots \to B_2\to
B_1 \to C_1$$ where $A_k=P_k$ and $B_k=H_k$ if the module $X_k$ is of type $U$
and
$A_k=H_k$ and $B_k= P_k$ otherwise.
\vglue12pt

We now give some low genus examples to illustrate the domains. The images are
developed images of 
the flat structures $\ogup$ and $\ogdn$ for $Gdh$ and $G^{-1}dh$
(respectively) on
the Riemann surfaces: the domain is locally always to the
left of the curve. We indicate the corners with angle $5\pi/2$ by a fat dot.

The genus $1$ example using just one handle is the Costa surface; see
Example 3.5.1.

For genus $3$, there are essentially two distinct possible symmetric module
sequences
with types
$(UUU)$ and $(UDU)$;
i.e., there are only these two up to
replacing
a `U' by a `D' and vice versa.
\figin{Zuuu}{600}

\centerline{Orthodisks for the $UUU$ candidate surface  }

\figin{Zudu}{600}

\centerline{Orthodisks for the $UDU$ candidate surface }
\vglue12pt

We will only be able to prove the existence of the second one. Numerical
investigations
give little hope that the first one also exists.
Starting with the second surface, we will inductively construct all our
other examples.

To see how the surface of type $T=(UD)^nU$ evolves as $n$ increases, 
here are the orthodisks for $UDUDU$:

\figin{Zududu}{650}

\centerline{The $UDUDU$ candidate orthodisks for a genus $5$ surface}
\vglue12pt

Geometrically, our choice of candidate surfaces can be characterized by the property 
that the intersection of the surface with the symmetry planes contains one 
arc component which connects all the finite points with vertical normal (on the $z$-axes) while
all other arc components connect the ends. 

The pair of orthodisks for the general  surface of type
$CT_n=(UD)^nU$-surface of genus $g=2n+1$  are
indicated
below. These surfaces have $2$ catenoid ends denoted by $C_1, C_2$;
$g$ planar ends denoted by $\{P_k\}$; and $g$ finite
points with vertical normals at the
handles denoted by $\{H_k\}$. 
There are $4n+4$ vertices in each orthodisk, of which $n$ are branch
points.
\vglue-60pt
\figin{Zudn1}{650}

\centerline{ The $\ogup$ orthodisk for $CT_n$}

\centerline{\BoxedEPSF{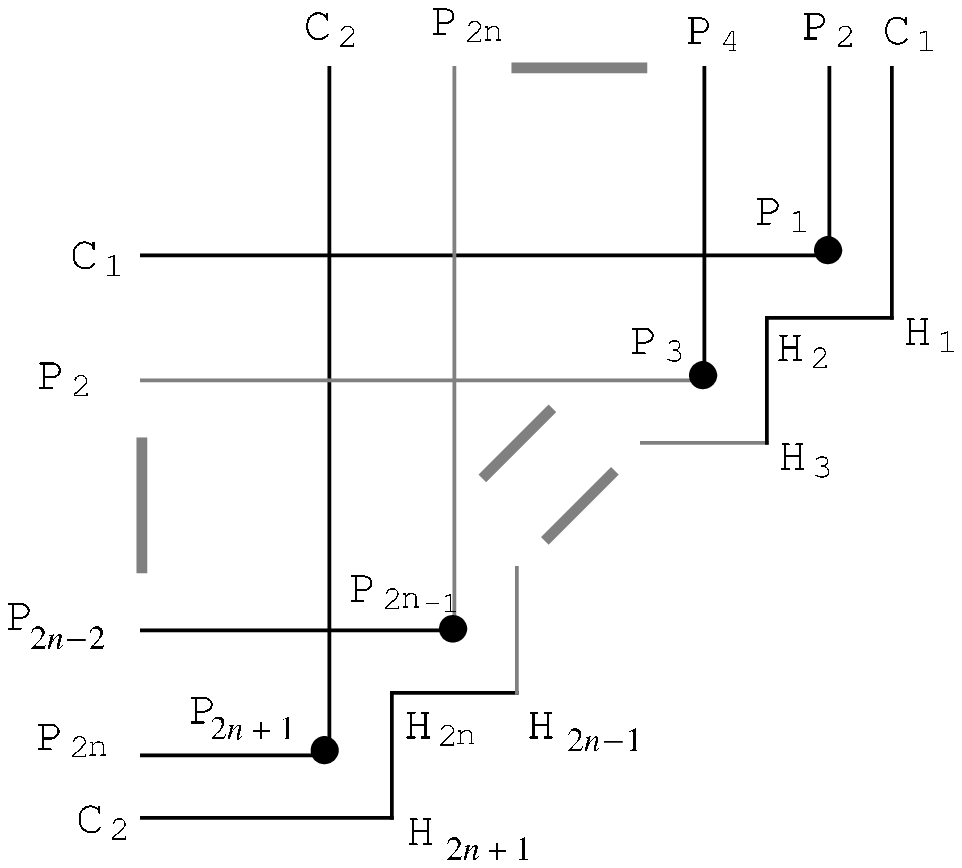 scaled 600}}
\vglue4pt
\centerline{The $\ogdn$ orthodisk for $CT_n$}

\demo{{\rm 3.7.}  More orthodisks by drilling holes}
In the previous section, we have described the candidate formal \wei data of a
family of minimal surfaces
generalizing Costa's surface for every odd genus greater than $1$. Now we
are going to add handles
(without enlarging the symmetry group) to all these surfaces in a symmetric
fashion. We call this process
{\it hole drilling}, and we will designate the resulting surfaces as 
${\rm DH}_{m,n}$ surfaces.

We start with an informal description.
For each odd integer $n \ge 3$ we describe a process of handle addition
which inductively
adds to the surface of type
$(UD)^nU$ an arbitrary number of handles without
adding ends.
Recall that our Costa-tower surfaces are composed of a stack of
modules called $U$- and $D$-module.

Now we allow two other kind of modules, called $F$- and $S$-module
(for front and side) which are vertical cylinders with handles
drilled through them from the front or the side, as indicated:
\vfill
\centerline{\BoxedEPSF{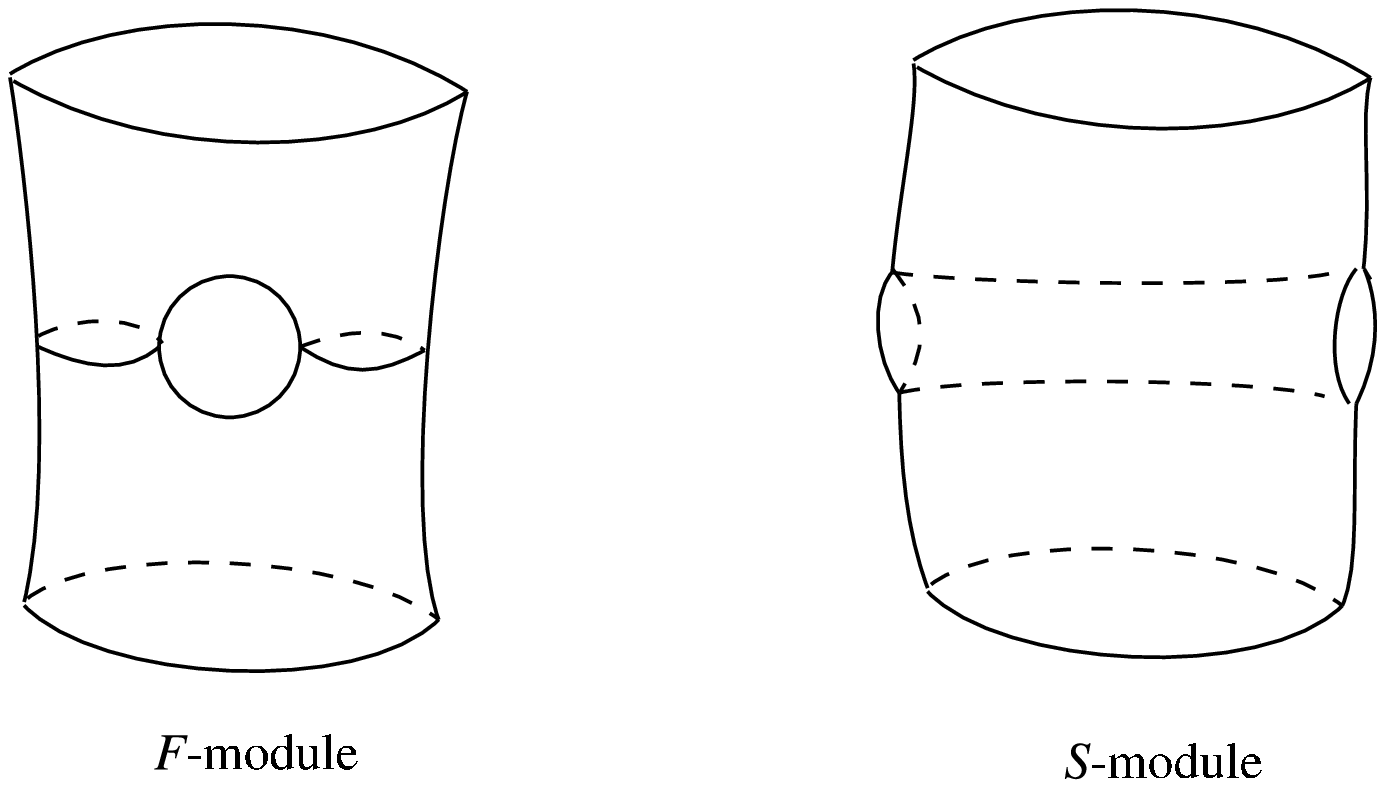 scaled 600}}
 \vglue-36pt
\centerline{Modules for drilling holes}
\eject

As in the case of the $U$- and $D$- modules, we are not able to prove existence
for surfaces made of arbitrary 
module sequences.

As for the Costa towers,  our choice of candidate surfaces is characterized by
the property  that the intersection of the surface with the symmetry planes contains one 
arc component which connects all the finite points with vertical normal (on the
$z$-axes) while all other arc components connect the ends. 

To be more specific, we will now add one handle to the $UDU$-surface.

There are apparently two possibilities according to the above rule, namely
$USDU$ and $UDFU$. But it is easy to see that the module sequences $SD$
and $DL$ generate equivalent formal
Weierstra\ss{}
data:
If one follows the intersection curves of these module sequences with the vertical
symmetry planes, one obtains the same pattern of special points at the $z$-axes,
hence the same formal \wei data.
So the topological significance of these modules does not impose much
restriction on the possible geometry of the surfaces.

\figin{Zsame}{600}

\centerline{Iterated hole drilling}

\demo{{R}emark {\rm 3.7.1}} This process does  work formally 
for $n=1$, but the candidate surface
would be Costa's surface with one handle added; this is known as the Horgan
surface.
See Nonexample 3.5.2(ii). 
\enddemo

We now describe the orthodisks arising from the handle addition process.
A module sequence of type $UF$ or $DU$ results in two consecutive finite
corners
of angle $\pi/2$ and $3\pi/2$ (or vice versa) of the orthodisks. Inserting
a handle according to the above 
rules results in adding two new finite corners
in between the old corners in a zigzag fashion. We illustrate this
for $USDU$ (compare with Example 3.5.1). 

\centerline{\BoxedEPSF{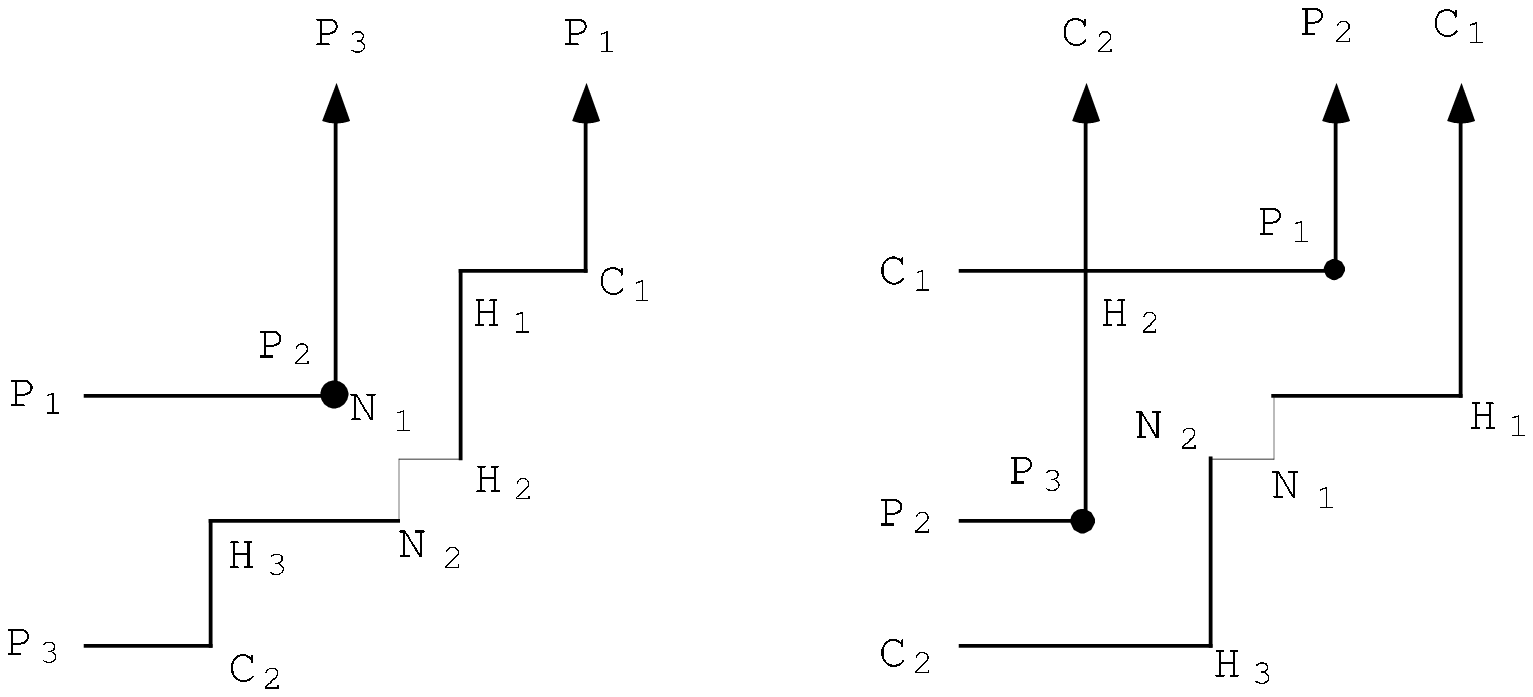 scaled 600}}

\centerline{Growing a handle by adding the two new corners $N_1$ and $N_2$}
\vglue8pt

Now that we have introduced these new module sequences, we can
extend our original one-parameter sequence of Costa towers to
a two-parameter sequence of orthodisks.
We introduce a new two-parameter family of orthodisks and give geometric
coordinates.  To set notation, we consider a pair of orthodisks with
$2n+1$ vertical
normals at the $z$-axis and
$2m+1$ planar ends. Together with the two catenoid ends we have $2m+2n+4$
Weierstra\ss{}
points so that the Riemann surface will have genus $m+n+1$. In the case
$m=n$ we are in the
situation of the Costa tower surfaces
of type $(UD)^nU$. We are aiming now for $n>m\ge0$, the cases in which
we will establish existence in the next three sections.

Here is a table depicting these surfaces and giving their
historical context:

$$\matrix
m\backslash n & 0 & 1 & 2 & 3 & 4 & 5 & 6 & 7\\
0 & \hbox{Costa} & \hbox{Horgan} & ? & ? &?&?&
? & ?
\\
1 &- & \hbox{$B$-$W$} & + & + & + & + & + & +
\\
2 &- & - & + & + & + & + & + &+
\\
3 &- & - & - & + & +  & + & + & +
\\
4 &- & - &- & - & + & + & + & +
\\
5 &- & - & - & - & - & + & + & +
\endmatrix$$

Meaning of the symbols:
\vglue4pt
\item{(1)} $?$: We have not been able to construct
these surfaces. Constructing one would imply the
existence of the others with higher genus. We conjecture
that these do not exist.\footnote{$^2$}{{\it Added in Proof}: This conjecture was confirmed in [MW].}
\vglue4pt\item{(2)} $-$: We will see that these do not exist. This is a very special
case of the
Hoffman-Meeks conjecture, and is the second portion
of the Main Theorem~B.
\eject
 \item{(3)} $+$: Dealt with in the following sections, this is the first portion of the main theorem.
\vglue4pt\item{(4)} $B-W$: This surface was found numerically by Boix and Wohlgemuth;
see \cite{Bo}.
 
\vglue-18pt
\section{The space of orthodisks}
\vglue-6pt
4.1. {\it Introduction}.
In this chapter, we will parametrize the space of pairs of conjugate orthodisks
by geometric data, set up a height function on this space which measures
the conformal
distance between the two orthodisks and prove its properness.
The main point in choosing the geometric coordinates lies in what we call
the {\it completeness condition:}
We require that whenever we have a sequence of points leaving every compact
set of the
coordinate space (that is, when a {\it geometric degeneration}
occurs), 
at least one of the
two orthodisks $\ogup$ and $\ogdn$ degenerates conformally.
After introducing complete geometric coordinates, the height function is set up
to measure differences in the extremal lengths
of cycles in a way 
which detects precisely the degenerations coming from a geometric
degeneration. 
Because it can easily happen that a geometric degeneration forces the
simultaneous conformal
degeneration of both orthodisks, we also need to  measure the rates of
growth or decay of extremal lengths.
This is done by a technical monodromy argument which is postponed to the final
section 4.7.  For the sake of clarity, we will now describe the height function and formulate
the monodromy theorem.

The height is constructed as a sum of terms, each of which
measures the possible conformal degeneration of the extremal length
of one cycle:

\demo{Definition {\rm 4.1.1}} Fix formal \wei data and consider a cycle $c$  in a conformal disk.
Given a pair of conjugate orthodisks $\ogup$ and $\ogdn$ with the chosen formal \wei data,
define the height of a cycle $c$ by
$$\height(c)=\left\vert e^{1/\ext_\ogup(c)}-e^{1/\ext_\ogdn(c)} \right\vert^	2+
\left\vert e^{\ext_\ogup(c)}-e^{\ext_\ogdn(c)} \right\vert^2.$$
\enddemo 

The definition of the height 
requires a choice of cycles $c_i$ and will be defined as

\centerline{${\displaystyle \height = \sum_{i} \height(c_i).}$}
\vglue4pt\noindent
This choice is restricted by two different requirements:
\vglue4pt
\item{(1)} Asymptotic computability.\ To prove properness, we need a precise
asymptotic formula in the case of geometric degenerations. Such formulas
are available only in special situations. For us, this means that we
have to choose simple cycles
or double cycles which are symmetric with respect to the
diagonal. Such cycles will be called {\it admissible} cycles.
\eject
\item{(2)} Avoidance of certain edges: In the next chapter, we need to compute
the derivative of the height with respect to changes in the coordinates which
can be described as `pushing edges'. There,  more cycles ending 
at the pushed edge yield more terms of the height function that 
we will need to control, so we aim not to have too many cycles
meeting a particular given edge.
\vglue12pt

These restrictions have the curious effect that it becomes easier to choose
the right cycles for surfaces of higher genus and with more topology; it is the lowest genus
case which causes most problems. And in fact, the reason that our 
present proof does not extend to prove
the existence of the nonexisting Horgan surface (see \cite{W}) 
is our inability to choose sufficiently many good cycles.

From the above form of the height function  it is clear that
we can prove properness when we can ensure that the extremal
length of a chosen cycle  tends to $0$ or $\infty$ for one
domain but remains bounded and bounded away from $0$ for the
other one. However,  we will encounter situations
where the extremal length of a cycle degenerates to zero (or infinity)
simultaneously in {\it both} domains.
Here we will need to invoke the following theorem, proven in 
\S4.7 as Lemmas 4.7.1 and 4.7.2 by a
monodromy argument:

\nonumproclaim{Monodromy Theorem 4.1.2} Fix formal \wei data and an admissible cycle $c$ {\rm (}\/see above\/{\rm ).}
Consider a sequence of  pairs of conjugate orthodisks 
$\ogup$ and $\ogdn$
such that either $c$ encircles a single edge in $\ogup$ shrinking
geometrically to zero and both $\ext_{X_1(p_n)}(\gamma)\to 0 $ and $\ext_{X_2(p_n)}(\gamma)\break \to 0$
or $c$ foots on an edge in $\ogup$ shrinking geometrically to zero and 
both $\ext_{X_1(p_n)}(\gamma)\to \infty $ and $\ext_{X_2(p_n)}(\gamma)\to \infty$.
Then $\height(c)\to \infty$ as $n\to\infty$.
\endproclaim

For the sake of a clear exposition, we have treated in detail the
special but most fundamental case $CT_3={\rm DH}_{1,1}$  in \S4.2 and 
\S4.3, postponing the general case
to the two subsequent sections, \S4.4 and \S4.5.
Historically, this surface was first discovered numerically
by Boix and Wohlgemuth (\cite{Bo}).

\demo{{\rm 4.2.} Geometric coordinates for the ${\rm DH}_{1,1}$ surface}
In this section, we describe geometric coordinates for our basic
surface ${\rm DH}_{1,1}$.

Denote 
$$\Delta=\{ (y,b,g)\in \Bbb R^3: 0<y,b,g<1, y+b+g=1\}$$
and by $\bar\Delta$ its closure. For each point in $\Delta$, we will associate
a pair of orthodisks.
We introduce three cycles, called yellow, green, and blue:

\centerline{\BoxedEPSF{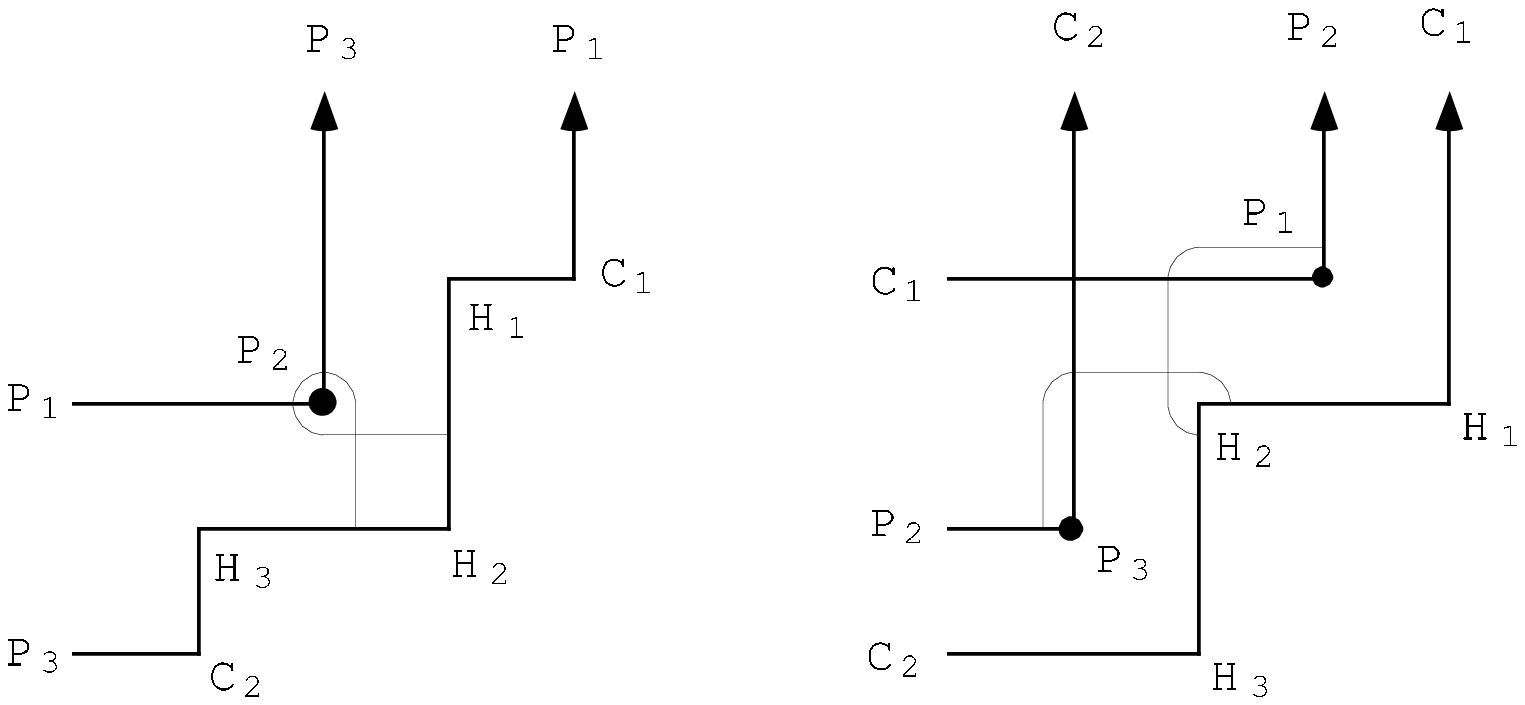 scaled 650}}

\centerline{yellow cycle}

\figin{Zu3blue}{650}
\vglue-24pt
\centerline{blue cycle}

\figin{Zu3green}{650}
\vglue-23pt
\centerline{green cycle}
\eject

Formally,
\vglue4pt
\item{(1)}
yellow consists of (undirected) arcs connecting $H_3 H_2$ with $P_1 P_2$
and symmetrically $H_2 H_1$ with $P_2 P_3$.
\vglue4pt\item{(2)}
blue consists of (undirected) arcs connecting $P_3 C_2$ with $H_3 H_2$ 
and symmetrically $C_1 P_1$ with $H_2 H_1$.
\vglue4pt\item{(3)}
green consists of (undirected) arcs connecting $P_1 P_2$ with $H_1 C_1$
and symmetrically $P_2 P_3$ with $C_2 H_3$.
\vglue4pt

Why should we choose these particular curves? At this stage, the
choice has little to do with the eventual geometry of the minimal
surface, and much to do with the conformal and combinatorial
geometry of the orthodisk: our goal is to meet criteria (1) and
(2) of the last subsection.

Now each orthodisk is  determined up to scaling by the lengths of the
periods of these cycles which we denote by $y,b,g$ subject to the condition
$y+b+g=1$. 
Hence for any triple $(y,b,g)\in \Delta$ we can form a pair of orthodisks
with these normalized period lengths. It is clear that these 
orthodisks will be conjugate.

We call $\Delta$ a geometric coordinate system for the formal \wei data
of type ${\rm DH}_{1,1}$. By this we mean the following:
\enddemo

{\it Definition} 4.2.1. Let formal symmetric \wei data be given. A
geometric coordinate system is an open
subset $\Delta$ of a Euclidean space
such that
for each point in $\Delta$, there is a pair of 
normalized symmetric conjugate
orthodisks with the given \wei data
such that the periods of the orthodisk are linear functions
in the coordinates of the point in
$\Delta$.
\vglue12pt

Here a pair of symmetric orthodisks is {\it normalized}
if the outer sheet boundary of the orthodisk $\ogup$ consists of finite edges
with total length equal to~$1$.

We will from now on only consider normalized orthodisks.

Here our geometric coordinate system $\Delta$ records the periods
for a set of cycles that span the homology of the (covering) surface,
hence provides sufficient information for determining all of the
periods via linear functions of the given ones.

In the rest of this section we will show that our coordinate geometric system
is {\it complete} in the following sense:

\demo{Definition {\rm 4.2.2}} A geometric coordinate system $\Delta$
is called {\it complete} if for any sequence
of points in $\Delta$ leaving any compact set, the conformal structure of
at least one of the orthodisks
degenerates.
\enddemo

This condition will ensure that our height function satisfies
the minimal necessary requirements for being proper.

\nonumproclaim{Lemma 4.2.3} The above geometric coordinate system $\Delta$ for
${\rm DH}_{1,1}$ is complete.
\endproclaim

{\it Proof}.
Consider a sequence $p_n\in \Delta$ with $p_n\to \partial \Delta$.
After choosing a subsequence, we can assume that $p_n\to p_0\in\partial\Delta$
Now we want to detect conformal degenerations of the orthodisk domains. For this
it will be convenient to use the extremal lengths of the cycle blue and another cycle
mauve introduced below:

\figin{Zu3mauve}{600}
\vglue-12pt
\centerline{mauve cycle}
\vglue12pt

Mauve consists of curve connecting $P_3 C_2$ with $ C_1P_1$.

Denote the extremal length of a cycle (say blue) in (say) the orthodisk $\ogup$
as usual by $\ext_\ogup({\rm blue})$.

\figin{Zdelta}{700}
\vglue-12pt
\centerline{Geometric coordinate simplex}
\eject

We need to distinguish six cases (which are indicated in the above figure), depending on the
location of
$p_0$ in
$\partial
\Delta$:
\vglue6pt
\item{(1)} $p_0 = (0,0,1)$: $\ext_\ogup({\rm blue})=0$,
\vglue6pt
\item{(2)} $p_0=(0,1,0)$: $\ext_\ogdn({\rm blue})=\infty,$
\vglue6pt
\item{(3)} $p_0=(1,0,0)$:  $\ext_\ogup({\rm blue})=0$,
\vglue6pt
\item{(4)} $p_0\in (>0,>0,0)$: $\ext_\ogdn({\rm mauve})=\infty$,
\vglue6pt
\item{(5)} $p_0\in (>0,0,>0)$: $\ext_\ogup({\rm blue})=0$,
\vglue6pt
\item{(6)} $p_0\in (0,>0,>0)$: $\ext_\ogup({\rm mauve})=\infty$.
\vglue6pt\noindent
Since the degeneration of an extremal length clearly signals the
degeneration of the conformal structure, we are done.

\demo{{\rm 4.3.} The height function for ${\rm DH}_{1,1}$}
Our aim is to define a proper function on $\Delta$ which
is zero if and only if the two domains 
$\ogup$ and $\ogdn$ given by the geometric coordinates of the
point are conformal. The conformal difference between the two domains is
naturally measured by expressions in the extremal length like the
Teichm\"uller distance, which, however, is not proper.

To obtain properness, we need to measure rather subtle growth
differences of extremal lengths.

Another important point for the choice of the height function is
that we need to be able to decrease the  height at noncritical points.
This requires that we control the first derivative of the height,
and this means heuristically  that we should use as few 
curve families as possible in the definition of the height.
\enddemo

{\it Definition} 4.3.1.
Consider on $\Delta$ the height function
$$\align
\SH_{1,1}=&\left\vert e^{1/\ext_\ogup({\rm blue})}-e^{1/\ext_\ogdn({\rm blue})} \right\vert^	2\\
&  +\left\vert e^{1/\ext_{\ogup}({\rm mauve})}-e^{1/\ext_\ogdn({\rm mauve})} \right\vert^	2 \\
&  +\left\vert e^{\ext_\ogup({\rm blue})}-e^{\ext_\ogdn({\rm blue})} \right\vert^	2\\
&  +\left\vert e^{\ext_\ogup({\rm mauve})}-e^{\ext_\ogdn({\rm mauve})}
\right\vert^	2.
\endalign$$
\vglue12pt

Clearly the conformal structure of a domain $\ogup$ or $\ogdn$ is determined by
the extremal lengths of blue and mauve, so that the height function is zero
if and only if  the two domains are conformal.

It is also clear that the height function admits at least the possibility
of being proper, because 
at a boundary point, at least one of the extremal lengths used is either $0$
or $\infty$ by Lemma 4.2.3. The main difficulty arises when the  
extremal lengths
for the same cycle degenerate for both orthodisks. This requires that we
measure extremal length growth rates and is responsible for the complicated
shape of the height function.

We can now prove, up to the crucial monodromy, Theorem~4.1.2.

\nonumproclaim{Theorem 4.3.2} The above height function $\SH_{1,1}$ is proper.
\endproclaim

\demo{Proof} Consider a sequence $p_n\in \Delta$ with $p_n\to \partial \Delta$.

After choosing a subsequence, we can assume that $p_n\to p_0\in\partial\Delta$
and that all extremal lengths of blue and mauve converge to some numbers in\break
$\Bbb R \cup  \{0,\infty\}$.

By Lemma 4.2.3,
at least one of the extremal lengths in the definition of the 
height function is $0$ or
$\infty$.
Now observe that in the cases (3), (4) and (5) of Lemma 4.2.3, 
$\ext_\ogdn({\rm mauve}) \to
\infty$ while
$\ext_\ogup({\rm mauve})$ remains bounded. Similarly, in case (6) 
$\ext_\ogup({\rm mauve}) \to \infty$
while
$\ext_\ogdn({\rm mauve})$ remains bounded.
This leaves us with the cases (1) and (2). Here both 
 mauve  extremal lengths go to $\infty$.

The point now is that we are able to control the growth rates of the 
extremal length
of  blue. It then follows from Theorem~4.1.2 that,
independently of the paths of approach to $\p \Delta$,
in case (1)
$$\vert e^{1/\ext_{\ogup}({\rm blue})}-
e^{1/\ext_{\ogdn}({\rm blue})} \vert^	2 \to\infty$$
while in case (2)
\vglue12pt
\hfill ${\displaystyle\vert e^{\ext_{\ogup}({\rm blue})}-e^{\ext_{\ogdn}({\rm blue})} \vert^
2\to\infty.}$
\enddemo

\vglue9pt
4.4. {\it Geometric coordinates for the ${\rm DH}_{m,n}$ surfaces.}
The goal of this section is to introduce a set of coordinates
for the moduli space of pairs of symmetric orthodisks with formal \wei data of type
${\rm DH}_{m,n}$
and to prove that they are complete.

Fix formal \wei data $a_i$ of a symmetric orthodisk  of type ${\rm DH}_{m,n}$ for
the rest of this section.

The corresponding pairs of  orthodisks have the following geometry:
The orthodisk $\ogup$ has an outer sheet bounded by a polygonal
arc with finite vertices
$$C_2 -H_{2n+1}-H_{2n} - \ldots - H_1- C_1$$
and $m$ interior sheets, each having one finite vertex, the branch
point $P_{2k}$.\eject

The orthodisk $\ogdn$ has an outer sheet bounded by a polygonal arc with
finite vertices
$$H_{2n+1}-H_{2n} - \ldots - H_1$$
and $m+1$ interior sheets containing  one finite branched vertex
$P_{2k+1}$ each.

A picture of such a pair of orthodisks is shown below for
${\rm DH}_{5,8}$.
In the following figures, the orthodisk $\ogup$
will always be the left one. We discard the
labeling of the vertices. 
\vglue-18pt
\figin{Zdh58}{700}
\vglue-18pt
\centerline{The ${\rm DH}_{5,8}$ orthodisks }
\vglue12pt

The first step is to set up a space of complete geometric coordinates
as we did
for the ${\rm DH}_{1,1}$ surface. We face here a new difficulty: 
it is not the case that for
every orthodisk $\ogup$ there is
necessarily a conjugate orthodisk $\ogdn$, and vice versa.
Thus,
the geometric coordinate space has to be an appropriate subset of the natural
parameter spaces of both orthodisks, and in addition we have to ensure that
the coordinates are complete.
We begin by an informal description with examples and then give a formal
definition. After we give the formal definitions, we prove the setup suffices
for our purposes in Lemmas~4.4.4 and 4.4.5.

The idea is to require that each branch point lie in a 
 box which is bounded by (extensions of) suitable orthodisk edges.
Not only will this allow for the existence of well-behaved moduli
 spaces,
but it will also provide for important degenerations of the domain
at the boundary.

In the following figures, these boxes are shaded in gray, and when necessary,
separated by additional (suggestive)
edges. When boxes overlap, they are shaded more darkly.
We also discard from now on the fat dots marking the interior branch point ---
all interior vertices are branch points.\eject

\centerline{\BoxedEPSF{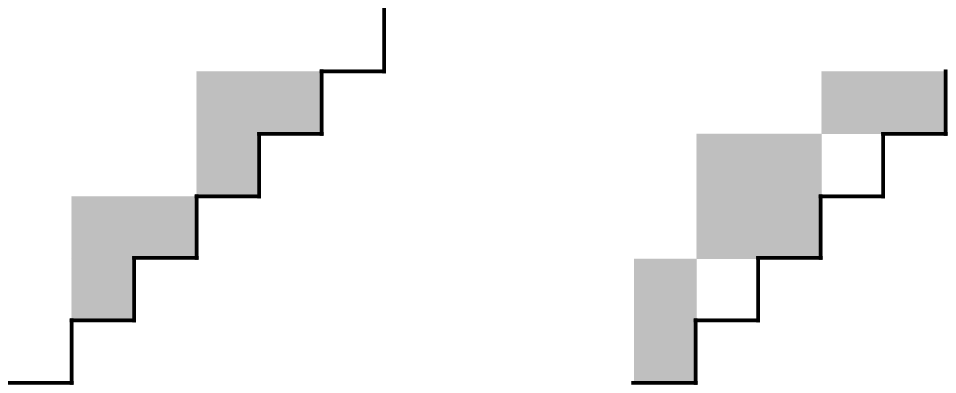 scaled 600}\qquad\BoxedEPSF{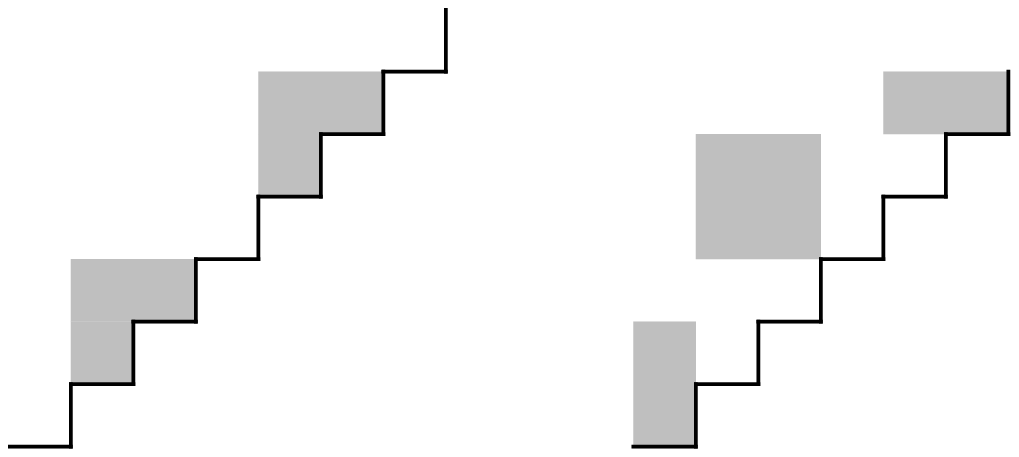 scaled 600}}
\vglue12pt
\centerline{Bounding boxes for ${\rm DH}_{2,4}$ and ${\rm DH}_{2,5}$ }
 \vglue18pt

\centerline{\BoxedEPSF{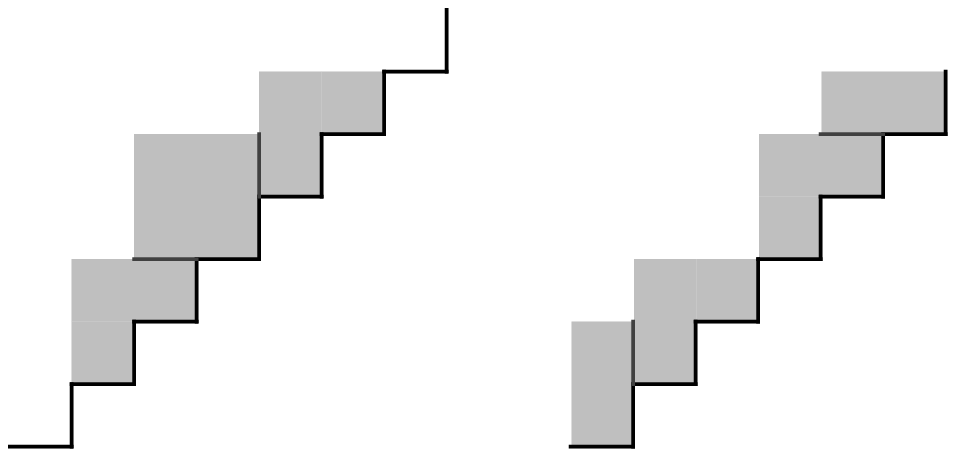 scaled 550}\qquad\BoxedEPSF{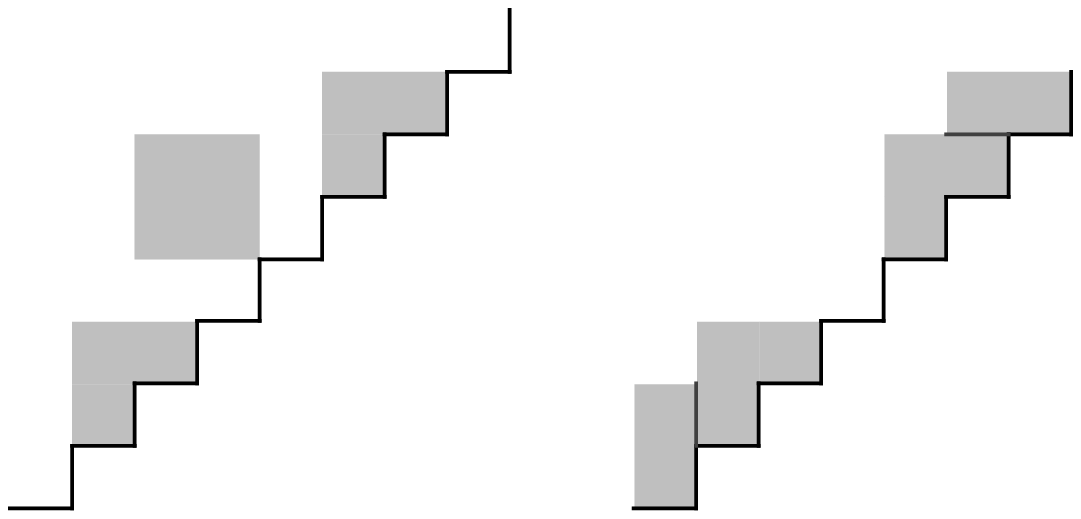 scaled 550}}
\vglue12pt
\centerline{Bounding boxes for ${\rm DH}_{3,5}$ and ${\rm DH}_{3,6}$ }
\vglue18pt

\centerline{\BoxedEPSF{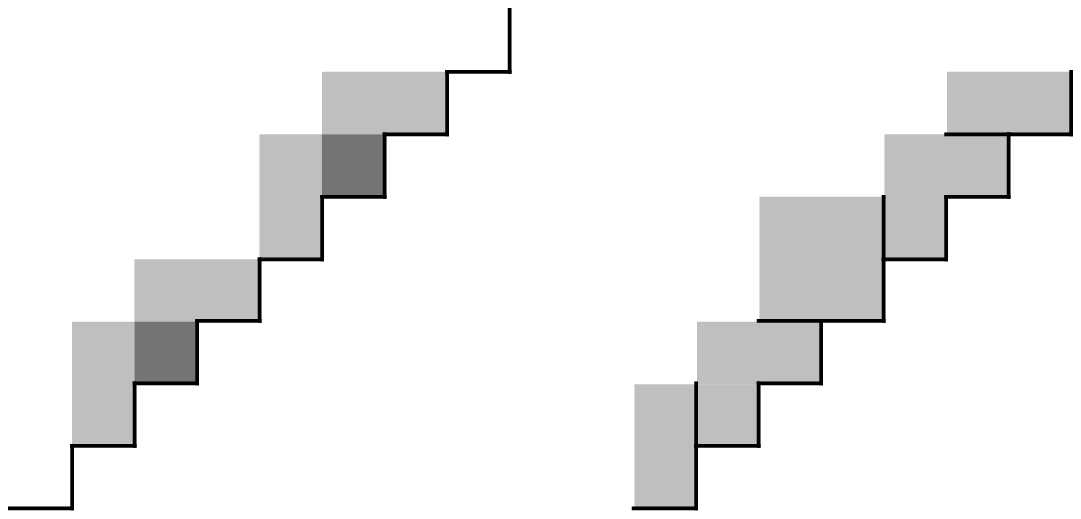 scaled 550}\qquad\BoxedEPSF{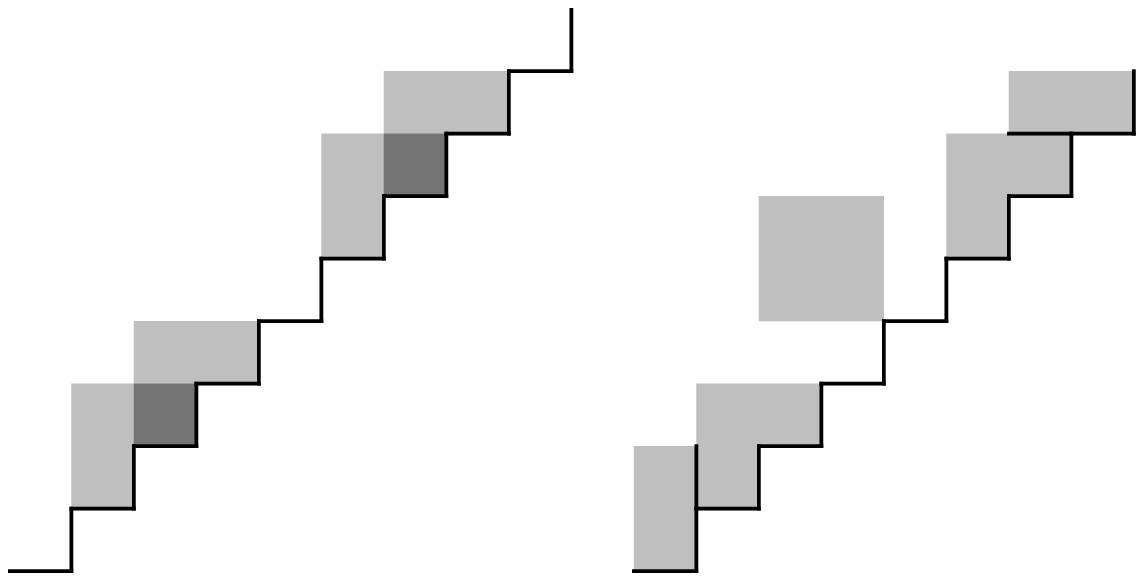 scaled 550}}
\vglue12pt
\centerline{Bounding boxes for ${\rm DH}_{4,6}$ and ${\rm DH}_{4,7}$ }
\vglue18pt
 \centerline{\BoxedEPSF{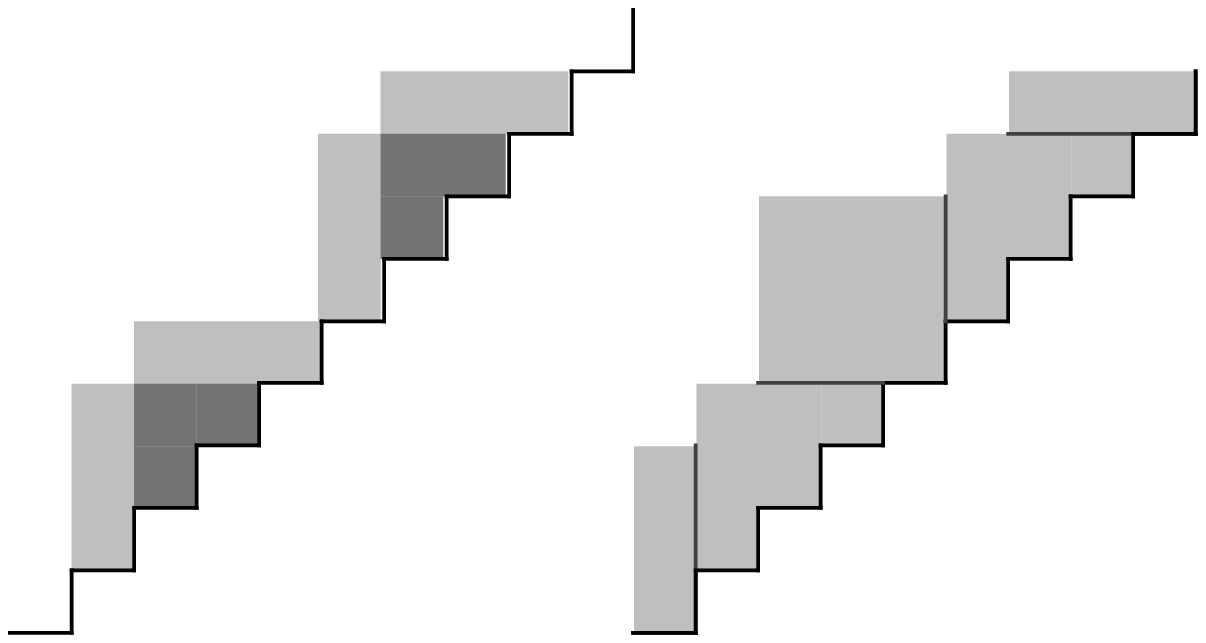 scaled 550}}
\vglue12pt
\centerline{Bounding boxes for ${\rm DH}_{4,8}$ }
\vglue12pt

As a general rule, there are $m$ boxes for $\ogup$ (one for each branch point) which are
intersections of finite rectangles with the orthodisks. For $\ogdn$, there are $m+1$
boxes, the two outermost of 
them being intersections of half-infinite rectangular strips with the
orthodisks.

As an illustration, we explain in detail how the boxes can be deduced
from the geometry for the ${\rm DH}_{2,4}$ surface:

In $\ogdn$, the edge $P_4P_5$ clearly has
to lie above the edge $C_2H_9$.
The conjugacy requirement 
for $\ogup$ forces that the
edge $P_4P_5$ lie to the right of the edge
$C_2H_9$. This forces, in $\ogup$, 
the edge $P_3 P_4$ to lie above the edge $H_8H_9$. Observe
that this is an inductive process to get border lines of the bounding boxes. In this most
simple case, the induction immediately terminates: by symmetry, the branch point $P_3$
has to lie on the diagonal in $\ogdn$, hence the edge $P_3 P_4$ lies left of $H_4 H_5$,
which causes, by conjugacy, the corresponding
edges to lie below each other in $\ogup$. This way we get
a backward induction in the general case, providing us with left and upper border
lines for bounding 
boxes for all branch points on and below the diagonal, up to the
outermost branch point in $\ogdn$ for which we do not get a left border line.
By symmetry, we also obtain such borders for the branch points above the diagonal.
Finally, the branch point on the diagonal is also bounded because the edge
$P_2P_3$ lies left of $H_1H_2$ in $\ogup$ which gives a bound for $P_3$
in $\ogdn$.

The idea now is to specify the set of geometric coordinates as follows:
First, record
the edge lengths of the boundary of the outer sheet of the orthodisk $\ogup$.
Then adjoin to those
the positions of the branch points within the bounding boxes, measured
as horizontal and vertical distances to the boundary.

In this ${\rm DH}_{2,4}$ example, the space of geometric coordinates
can be described as the space of $\ogup$ (left)
orthodisks with the requirement that the branch points lie inside the bounding boxes.
 Each
such orthodisk is given by the data for its outer sheet boundary zigzag and
the position of one point in the L-shaped box
which itself is defined by the boundary zigzag.
Using these data, it is possible to construct
a unique conjugate orthodisk $\ogdn$ such that
its branch points lie in the shaded bounding boxes.
Furthermore, these coordinates are complete;
any degeneration of the coordinates which
does not come from an edge degeneration in $\ogup$
stems from one branch point coming close to the
inner boundary of its shaded box. 
But this forces
a branch point in $\ogdn$ to come close to the
outer boundary, forcing a conformal degeneration.
This concludes the discussion of the
case ${\rm DH}_{2,4}$; the completeness of the coordinates for the general 
case of ${\rm DH}_{m,n}$ is proven in Lemma~4.4.5.

As a first step towards formally describing these boxes
in the general case, we introduce
families of cycles which will be useful both for specifying the
coordinates and for defining the height function. In principle, there are
four cases to distinguish depending on the parities of $m$ and $n$, but
our notation will hide most of the differences.

Denote
$$ k= [m/2] \qquad\hbox{and}\qquad d =[\frac{n-m}2]+1.$$

We first introduce cycles related to the outer sheet:
$$\alignedat3
\alpha_1 &: P_1C_1 &\to &\qquad H_1 H_2  \\
\alpha_2 &: C_1 H_1 &\to &\qquad H_2 H_3 \\
\alpha_3 &: H_1 H_2 &\to &\qquad H_3 H_4 \\
\cdots &: \cdots &\to &\qquad \cdots \\
\alpha_{2n}&: H_{2n-2}H_{2n-1} &\to &\qquad H_{2n}H_{2n+1} \\
\alpha_{2n+1}&: H_{2n-1}H_{2n} &\to &\qquad H_{2n+1}C_2 \\
\alpha_{2n+2}&: H_{2n}H_{2n+1} &\to&\qquad P_{2m+1}C_2.
\endalignedat$$

These cycles connect consecutive parallel edges of the outer sheet
boundary. Their periods are realized as the finite edges. The symmetry
condition on the orthodisk implies for these cycles
$$\vert \per \alpha_k \vert = \vert \per \alpha_{2n+3-k} \vert .$$

Next we introduce four families of cycles called
$\lambda_k, \rho_k, \upsilon_k$ and $\delta_k$ (for left, right, up, and
down) connecting  the horizontal and vertical
edges of the inner sheets
 to certain carefully chosen horizontal and vertical edges of the outer
sheet.
The periods of these cycles will be used to define the bounding boxes:

$$\alignedat3
\rho_1&: P_2P_3 &\to&\qquad H_1H_2 \\
\rho_2&: P_4 P_5 & \to &\qquad H_3H_4 \\
\cdots &: \cdots &\to&\qquad\cdots \\
\rho_k&: P_{2k}P_{2k+1} &\to&\qquad H_{2k-1}H_{2k} \\
\upsilon_1 &: P_1P_2 &\to&\qquad H_1 C_1 \\
\upsilon_2 &: P_3P_4 &\to&\qquad H_3H_2 \\
\cdots &: \cdots &\to&\qquad\cdots \\
\upsilon_k &: P_{2k-1}P_{2k} &\to&\qquad  H_{2k-1} H_{2k-2} \\
\delta_1 &: P_1P_2  &\to&\qquad H_{d+3}H_{d+2} \\
\cdots &: \cdots &\to&\qquad\cdots \\
\delta_k &: P_{2k-1}P_{2k}  &\to&\qquad  H_{2k+d+1}H_{2k+d}\\
\lambda_1 &: P_2P_3  &\to&\qquad H_{d+3}H_{d+4} \\
\cdots &: \cdots &\to&\qquad\cdots \\
\lambda_k &: P_{2k}P_{2k+1}  &\to&\qquad H_{2k+d+1}H_{2k+d+2} .
\endalignedat$$

In the case that $m$ is odd, we introduce in addition
cycles for the box at $P_{k+1}$. For odd $n$ we set
$$\alignedat3
\rho_{k+1} &: P_{m+1}P_{m+2} &\to&\qquad H_{n+1}H_{n+2} \\
\delta_{k+1} &: P_mP_{m+1} &\to&\qquad H_nH_{n+1} \\
\upsilon_{k+1} &: P_mP_{m+1} &\to&\qquad H_{n-d-2}H_{n-d-1} \\
\lambda_{k+1} &: P_{m+1}P_{m+2}&\to&\qquad H_{n+d+3}H_{n+d+4}
\endalignedat$$
and for even $n$:
$$\alignedat3
\rho_{k+1} &: P_{m+1}P_{m+2} &\to&\qquad H_{n}H_{n+1} \\
\delta_{k+1} &: P_mP_{m+1} &\to&\qquad H_{n+1}H_{n+2} \\
\upsilon_{k+1} &: P_mP_{m+1} &\to&\qquad H_{n-d-1}H_{n-d} \\
\lambda_{k+1} &: P_{m+1}P_{m+2}&\to&\qquad H_{n+d+2}H_{n+d+3} .
\endalignedat$$

This defines cycles emanating from the inner $\ogup$ sheet edges
for the sheets containing the branch points $P_2, \ldots, P_k$, that is,
for all sheets above the diagonal. For the sheets below the diagonal, we
employ symmetry: Denote by $\alpha'$  the image of a cycle $\alpha$
under reflection at the diagonal $y=-x$.
We then define
$$\align
\lambda_j &= \upsilon'_{m+1-j} \\
\rho_j &= \delta'_{m+1-j} \\
\upsilon_j &= \lambda'_{m+1-j} \\
\delta_j &= \rho'_{m+1-j} .
\endalign$$

This way all cycles are defined for $j=1,\ldots, m$.

The next two figures  show all the cycles for ${\rm DH}_{5,8}$. For visibility,
we have omitted both the names of the vertices and the fat dots indicating
the branch points; we have also shortened the edges of the inner sheets.

The first figure shows the $\upsilon$-cycles and the $\delta$-cycles
(thicker
and shaded),
the second the $\rho$-cycles and the $\lambda$-cycles (thicker and shaded).
Informally, for each inner sheet of the orthodisk $\ogup$ there is a cycle of each type
($\upsilon, \delta,\rho,\lambda$). The
$\upsilon$-cycles have periods which point upward in
$\ogup$ for the first half of the cycles in the upper-right part of the orthodisk,
and similarly the periods of  $\rho,\lambda,\delta$ point right, left and down.

\centerline{\BoxedEPSF{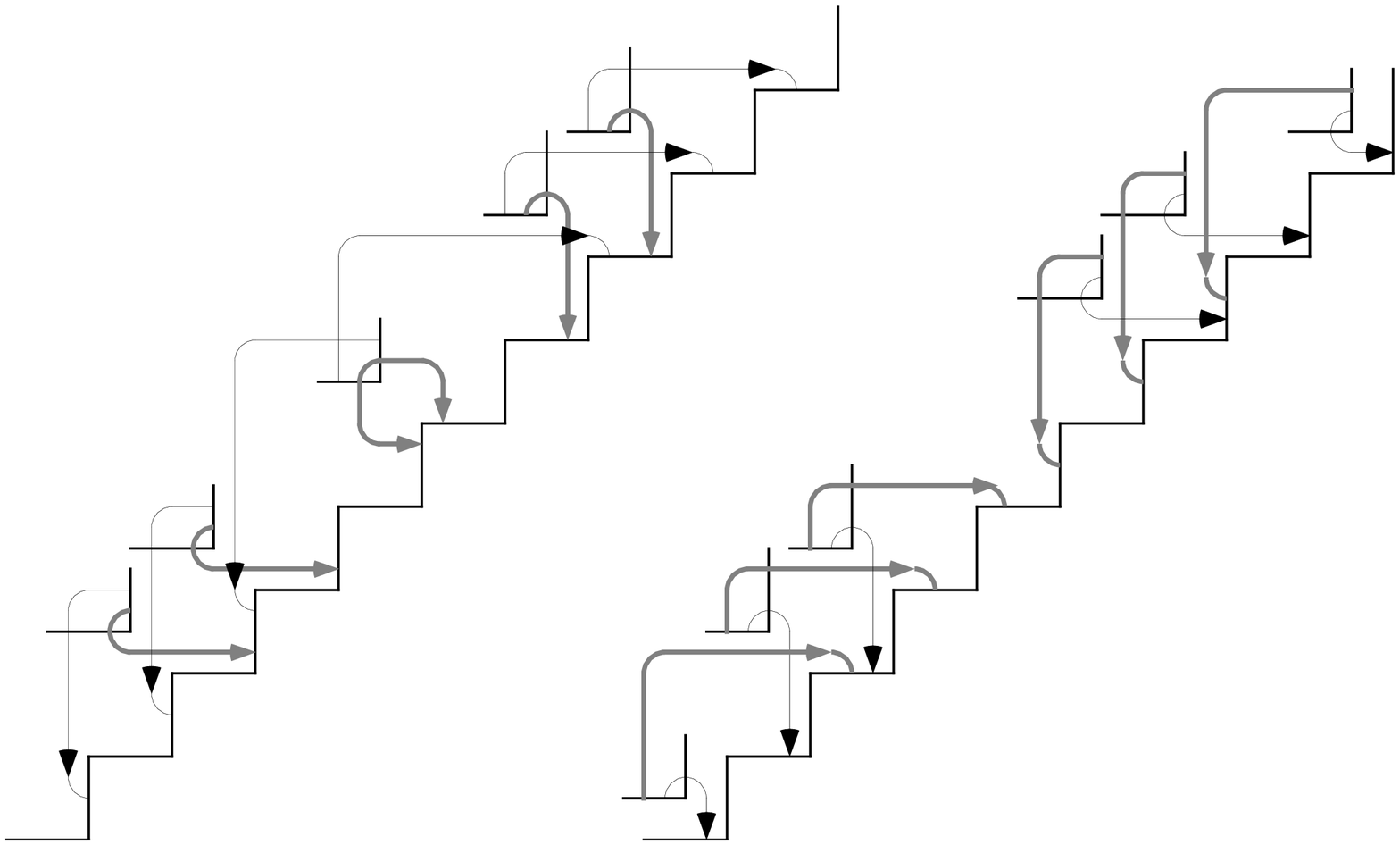 scaled 500}}
\vglue12pt
\centerline{$\upsilon$-cycles and the $\delta$-cycles for ${\rm DH}_{5,8}$}
 \vglue24pt

\centerline{\BoxedEPSF{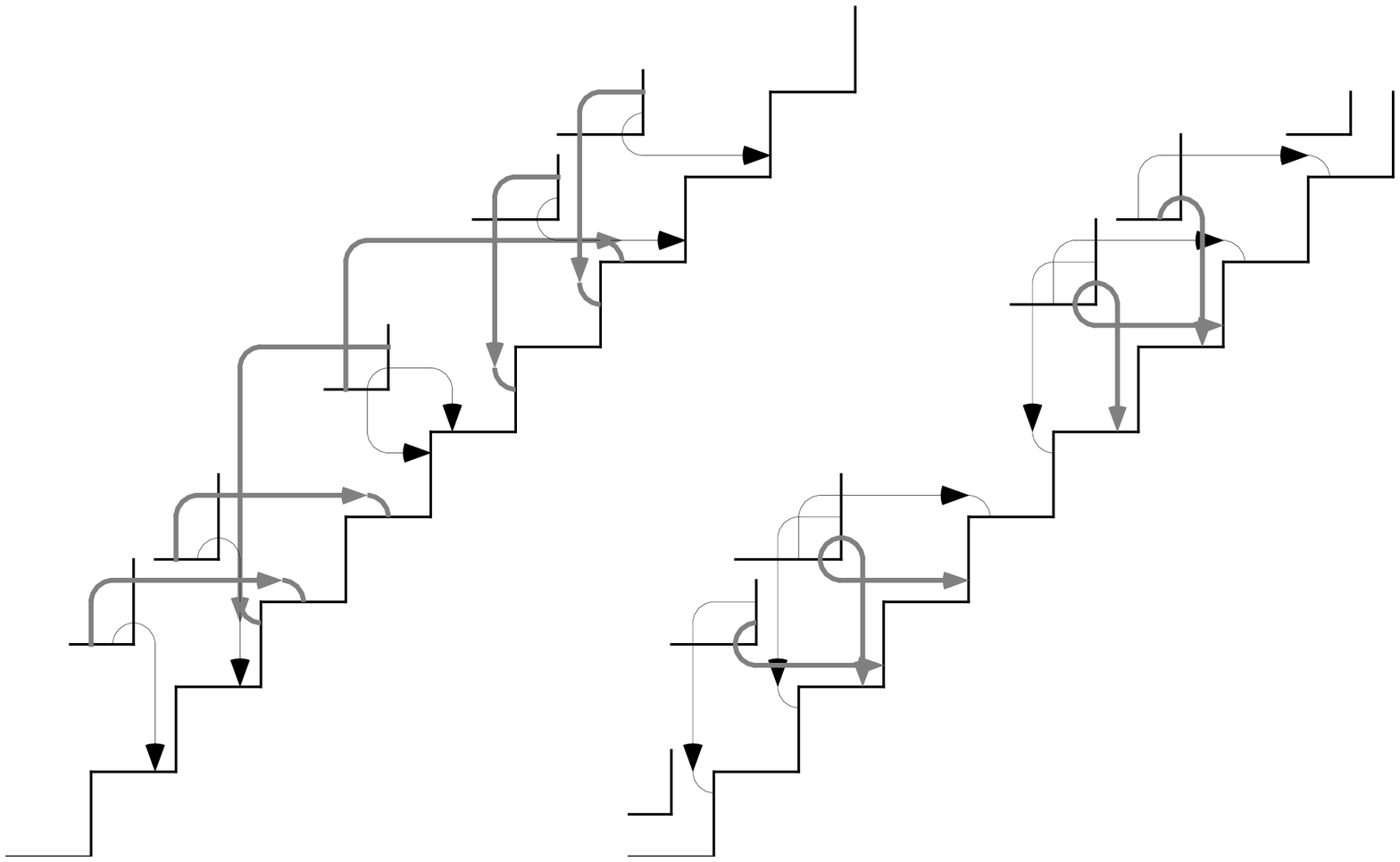 scaled 500}}
\vglue12pt
\centerline{$\rho$-cycles and the $\lambda$-cycles for ${\rm DH}_{5,8}$}
\vglue12pt

Now we are able to define the bounding boxes for the branch
points: 

\demo{Definition {\rm 4.4.1}} Denote by
$B_{2j}$ the rectangle defined by
$$\per \upsilon_j >0, \per \delta_j >0, \per \lambda_j>0, \per \rho_j>0$$
and by $B_{2j-1}$ the rectangle defined by
$$\per\rho_{j-1}>0,\quad \per\delta_{j-1}>0,\quad\per
\lambda_j>0,\quad\per\upsilon_{j}>0 .$$
\enddemo

Note that the boxes $B_1$ and $B_{2m+1}$ require the additional
cycle $\lambda_0:=\alpha_0$ whereas $\rho_0$ is undefined so that the
rectangle is open at one side.

To define the geometric coordinates, we impose two restrictions
on the branch points:
\vglue4pt
\item{(1)} $P_j \in B_j,$
\vglue4pt\item{(2)} $P_j$ must lie within the outer sheet.
\vglue4pt

While the second requirement is obviously necessary, the first
will imply that the geometric coordinates are complete.

\demo{Definition {\rm 4.4.2}} An orthodisk 
$\ogup$ or $\ogdn$ is called
{\it admissible} if it satisfies  the above condition for their respective
branch points. A pair of orthodisks is called admissible if the
orthodisks are conjugate and both are admissible.
\enddemo

{\it Definition} 4.4.3. The {\it geometric coordinates} of $\ogup$
are given by the periods of the cycles $\alpha_j, \rho_j, \lambda_j,
\upsilon_j,
\delta_j$ subject to the above condition. Similarly, the geometric
coordinates of $\ogdn$ are given by the periods of these cycles subject
to the above condition. This clearly defines two open subsets of a
Euclidean space of dimension $g+1=m+n$. Each subset parametrizes the
configuration space of admissible orthodisks $\ogup$ or $\ogdn$.
Denote by $\Delta$ the intersection of these subsets; this set
parametrizes admissible pairs of orthodisks and is called {\it geometric
coordinate space} for ${\rm DH}_{m,n}$.
 
\nonumproclaim{Lemma 4.4.4} $\Delta$ is an open cell of dimension $g-1$.
\endproclaim
 
\demo{Proof} To prove openness we just  observe that all imposed
conditions are open in a certain linear space. To compute its dimension,
we can easily get rid of the redundant equations of the definition of
geometric coordinates for (say)
$\ogup$: The cycles $\delta_j-\upsilon_j$ and $\lambda_i-\rho_i$
are a certain sum of cycles
$\alpha_j$. So to specify a symmetric orthodisk, we need just $m+n$
cycles. To prove that it is a nonempty cell, we construct a retraction
to an interior point as follows: First, 
by using convex combinations of the outer sheet edge lengths,  
we can easily deform a given pair
of admissible orthodisks to an admissible
pair of orthodisks with the outer sheet periods being of length $1$.
Now we move the edge $P_{2m}P_{2m+1}$ in $\ogdn$ down; this is
unobstructed in
$\ogdn$ but has to be accommodated by a movement of $P_{2m}P_{2m+1}$ to
the left in $\ogup$. We move so far to the left here, that we are left of
the edge $H_{2n-1}H_{2n-2}$. This allows us to move 
the edge $H_{2m-1}H_{2m}$ down
in $\ogup$. We move so far down that we are, in 
$\ogdn$, left of the edge $H_{2n-2}H_{2n-3}$, which enables us again to move
$H_{2m-2}H_{2m-1}$ down. We continue this process, doing the same
with the symmetric edges simultaneously. If we agree to move always to a
position $0.5$ apart from the left or lower boundary sheet, we will reach,
in this way, a canonical point in the space $\Delta$. This defines the
retraction.
\enddemo

\nonumproclaim{Lemma 4.4.5} $\Delta$ provides complete geometric coordinates.
\endproclaim
 
\demo{Proof}
We first have to show that each 
point in $\Delta$ gives a pair of
well-defined symmetric (normalized)
orthodisks. Using the coordinates of such a point, we see  immediately 
that one can
reconstruct the outer sheets for the $\ogup$ and the $\ogdn$ 
domains. Also, the positions
of the branch points are specified so that they lie within the outer sheets and their
respective bounding boxes. 
To obtain a correct pair of orthodisks, we have to make sure
that these branch points do not coincide as one might expect for instance in the ${\rm DH}_{5,8}$
case where the bounding boxes overlap. However, this can be excluded as follows:
If two branch points come close to each other, the periods of the two cycles connecting the
parallel edges of their respective sheets vanish. This causes in the other domain a vanishing period
between two inner 
sheets closer to the diagonal, and inductively we produce such a period degeneration on the diagonal,
switching between the two domains. But on the diagonal, this cannot happen, because the branch
point there is confined to the diagonal.

Finally, we have to show that when we have a sequence of admissible orthodisks
leaving the space of geometric coordinates, at least one of the conformal
structures of the
$\ogup$ or
$\ogdn$ orthodisks degenerates. The geometric degeneration means by definition
that for at least one of the orthodisks, either an edge of the outer sheet
boundary degenerates or a branch point hits the boundary of its
bounding box. (The exceptional case of the half-infinite bounding 
rectangle in the domain $\ogdn$ domain
does not actually allow the branch 
point to drift to infinity, because the corresponding period in the 
$\ogup$ domain is represented by a finite 
edge of the outer sheet, which is bounded by
normalization.) In the first case, we get a node and are obviously done.
In the second case, the branch point can either converge to some
point on the outer boundary (in which case we pinch a cycle and are also
done), or it can converge to a point on one of the virtual border lines of the
bounding box.
Then we  clearly have no conformal degeneration for {\it this} orthodisk sequence,
so we must prove that we get such a degeneration for the other family.
(The next figure shows two conjugate orthodisks with bounding boxes. Equally
labeled arrows indicate equal distances from orthodisk edges to the virtual border lines of
bounding boxes. As long as these distances are equal, the branched point can move freely within
their bounding boxes.)

Suppose, for concreteness, that we have a sequence of orthodisks $\ogup(n)$
where a branch point converges to an upper virtual border line of
its bounding box. This branch point $P_k$
belongs to an edge $E=P_kP_{k-1}$ which comes
close to this virtual line from below. By conjugacy, in the
corresponding 
$\ogdn$ orthodisks, the corresponding edge $E^*=P_kP_{k-1}$ comes 
close to a bounding border line to the right.
The other endpoint $P_{k-1}$ 
of this edge is a branch point for $\ogdn$ which hence
comes arbitrarily close to the outer sheet boundary. 
So we pinch a curve
in $\ogdn$ and get a conformal degeneration. The other cases are treated
similarly.
\enddemo

 \vglue12pt
\centerline{\BoxedEPSF{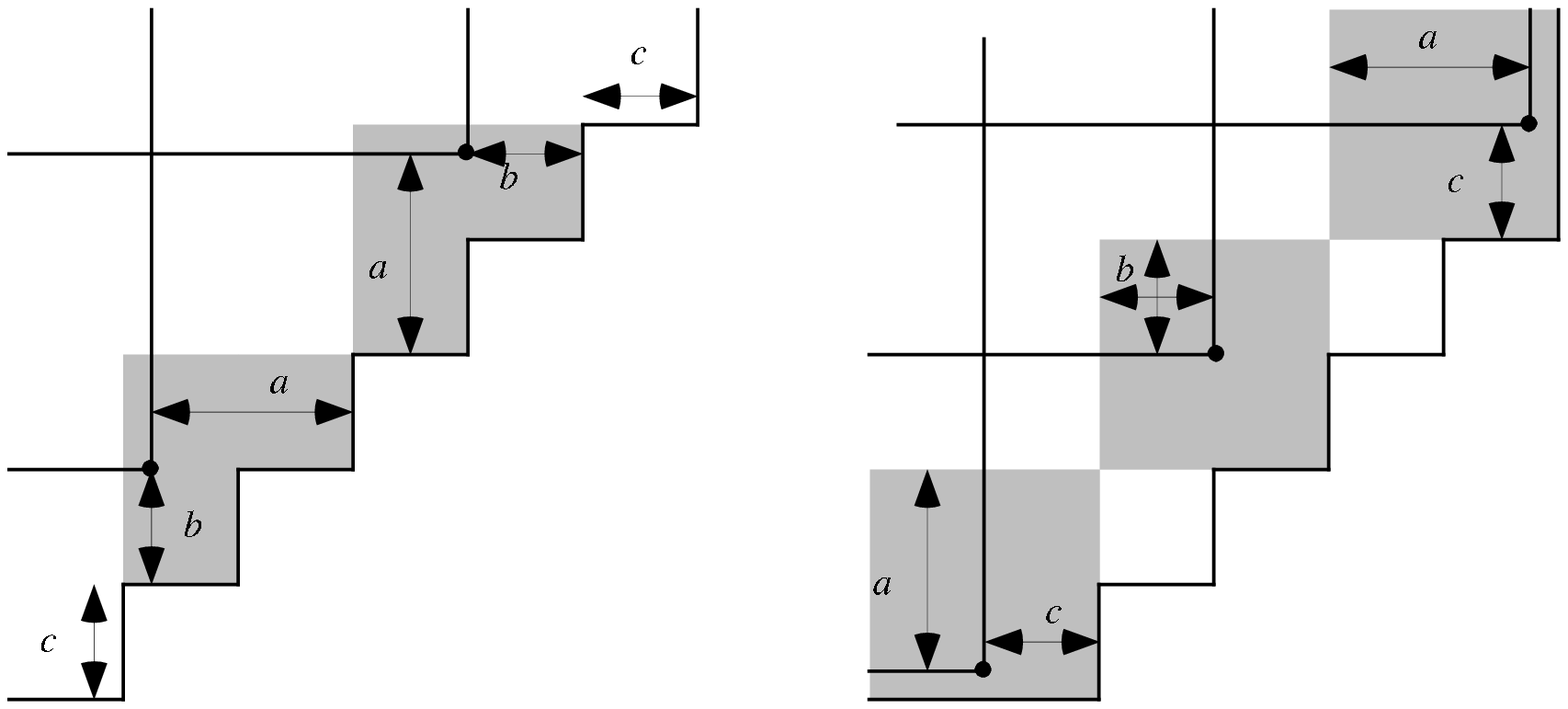 scaled 700}}
\vglue12pt 
\centerline{Conjugate orthodisks with coupled branch points}
\vglue12pt
\noindent

4.5.  {\it Height functions for} ${\rm DH}_{m,n}$.
In this section, we will provide a proper height function
defined on the geometric coordinate space $\Delta$ of the previous
section for the ${\rm DH}_{m,n}$ surfaces with $m+n>2$.

We begin by stating the basic requirements on a height function which
we require for the steps of the proof of Theorem~B(i) in \S\S4--6.

To estimate the rate of 
growth/decay of extremal lengths in terms of geometric degenerations,
we need an asymptotic expression for the extremal lengths 
of the cycles we use
in the height function. Such an expression (see \cite{Oht})
is only known to us
in the case that the cycle is a 
lift of an arc in the upper half-plane connecting two
disjoint edges. This is the case of a conformal quadrilateral. Informally,
an {\it admissible}
cycle in an orthodisk is one which can be reduced to such a cycle
as just described. More
precisely, we recall from \S4.1 the 

\demo{ Definition {\rm 4.5.1}} A cycle in an orthodisk is called {\it admissible}
if it is either simple and connects a pair of symmetric edges or has two
symmetric
components which do not cross the diagonal.
\enddemo

Using admissible cycles,  we will compose the height as a sum of the
following terms:

\demo{Definition {\rm 4.5.2}} Let $c$ be an admissible cycle. Define
$$\height(c)=\vert e^{1/\ext_\ogup(c)}-e^{1/\ext_\ogdn(c)} \vert^	2+
\vert e^{\ext_\ogup(c)}-e^{\ext_\ogdn(c)} \vert^2 .$$
\enddemo

Recall that the conformal polygon of a ${\rm DH}_{m,n}$-surface is 
(two-fold) symmetric
and has $2m+2n+4$ vertices; there are thus $m+n$ conformal moduli of
such
a shape.

In addition to the requirement that the cycles detect conformal degenerations,
we will also need some more properties which we will state 
here but only explain and use
in \S6. (There, we take a low genus reflexive orthodisk,
and append some thrice punctured spheres to it along nodes;
we then open the nodes to obtain high genus surfaces in 
a stratum $\SY \subset \Delta$ along which we flow in \S5 to a reflexive
orthodisk solution to our problem.) This 
`regeneration' in \S6  requires
a set of $m+n$
cycles for the height which 
provide conformal coordinates for the points in $\Delta$; moreover,
the set of curve systems,  
specified for the higher genus
surface,  which do not degenerate as we
pinch the surface to a noded surface (i.e.,
are either not pinched into the 
node or do not cross the node)
should also provide local coordinates for the
boundary stratum of lower genus surfaces. 

In particular,
there should be exactly two `bad' curves (which
degenerate as ${\rm DH}_{n+1,n+1}$ degenerates to ${\rm DH}_{n,n}$, i.e., with extremal
lengths which go to $0$ or $\infty$ as $P_1$ and $H_1$ converge)
in the $n=m$ case, and one bad curve in the $n \ne m$ case (i.e., which 
degenerates as $m+1 \to m$ in ${\rm DH}_{m+1, n}$)
with an extremal length that
goes to $\infty$ as the 
finite vertex nearest the central vertex approaches it.
Also, as we shall see in \S5, the technical requirements of 
the Noncritical Flow Step will restrict the edges which are
allowed for the feet of the bad curves; we postpone the
discussion of this until after Proposition~5.3.2 in \S5.3.)

For the definition of the height function, we distinguish two cases:
This distinction comes from the inductive design of the proofs of Theorem
A and B: In the case $CT_n={\rm DH}_{n,n}$, we use the edge $C_1H_1$ for
pushing, and this requires the cycles to be chosen so that no cycle
foots on this edge and precisely
one cycle encircles it. In the case that we want to prove existence of a
${\rm DH}_{m,n}$ with $m<n$, we need the central edges of the outer sheet
boundary zigzag for pushing.

Besides that, the choice of the cycles for the height function follows
the same general rules: By using cycles which connect infinite edges with
infinite edges, we detect degeneracies of inner sheets, and by using cycles
which encircle or foot on outer sheet edges, we enforce properness using
the monodromy lemma.

Finally, as a general rule, all our cycles will be symmetric with respect
to the diagonal.

\phantom{too weirtd}
\vglue-8pt
\item{(1)} $\mu$ connects $C_1P_1$ with $C_2P_{2m+1}$.

\phantom{too weirtd}
\vglue-8pt
\item{(2)} $\nu$ connects $C_1P_1$ with $H_1 H_2$ and
 $C_2P_{2m+1}$ with $H_{2n+1}H_{2n}$.

\phantom{too weirtd}
\vglue-15pt
\item{(3)} $\rho$ connects $H_n H_{n+1}$ with $P_{m-1}P_m$ and
 $H_{n+1}H_{n+2}$ with $P_{m+1}P_{m+2}$.
 
\phantom{too weirtd}
\vglue-8pt
\item{(4)} $\sigma$ connects $H_1C_1$ to $P_1 P_2$ and $H_{2n+1}C_2$ to $P_{2m} P_{2m+1}$.
\vglue4pt
\item{(5)} $\delta$ connects $P_1P_2$ to $H_nH_{n+1}$ and $P_{2m} P_{2m+1}$ to $H_{n+2}H_{n+3}$.
\vglue4pt
\item{(6)} $\gamma$ connects $H_{n-1}H_n$ to $H_{n+2}H_{n+3}$.
\vglue4pt
\item{(7)}  $\beta_k$ connects $P_1P_2$ with $P_{2k-1} P_{2k}$
and  $P_{2m} P_{2m+1}$ with $P_{2m+2-k}P_{2m+1-k}$.
\vglue4pt
\item{(8)}
$\alpha_k$ encircles $H_k H_{k+1}$ and  $H_{2n+2-k} H_{2n+1-k}$.
\vglue4pt
\item{(9)}
$\tau$ connects $P_1P_2$ with $P_{2m}P_{2m+1}$.
\vglue4pt

The cycle $\rho$ has been chosen to replace $\alpha_{n-1}$ which 
would interfere with
the pushing edge $H_1H_2$ for small $n$.

Now we define the cycles for the height function:

In the case $m=n=1$, select the cycles $\mu$ and $\nu$.
In the case $m=n>1$,
select the cycles $\mu,\nu, \gamma,\delta,\beta_2,\ldots,\beta_{m-1},
\alpha_2,\ldots,\alpha_{n-1}$. 
In the case $m=1<n$, select $\mu, \rho, \sigma,
\alpha_1,\ldots,\alpha_{n-2}$ and in the case $1<m<n$ select
$\mu,\rho,\sigma,\tau,
\alpha_1,\ldots,\alpha_{n-2},\beta_2,\ldots,\beta_{m-1}$.

Note that if $n=2$ there are no $\alpha$-cycles.
In all cases, we have selected $m+n$
cycles.

Observe that for $m=n=1$, the cycle $\mu$ becomes 
mauve and $\nu$ becomes blue so that
the height defined here coincides 
with the one defined in \S4.3 for ${\rm DH}_{1,1}$.

\demo{Definition {\rm 4.5.3}} The height for the ${\rm DH}_{m,n}$ surface is defined as the sum of the
heights of all selected cycles for the values of $(m,n)$.
\enddemo

{\it Example} 4.5.4. As an example, we show the cycles for the height
of the ${\rm DH}_{2,4}$ and the ${\rm DH}_{3,3}$ surface which are relevant for the properness proof:
\vglue-12pt
\figin{Zheight24}{625}
 
\centerline{The cycles for the height function for ${\rm DH}_{2,4}$ }
 \eject

\centerline{\BoxedEPSF{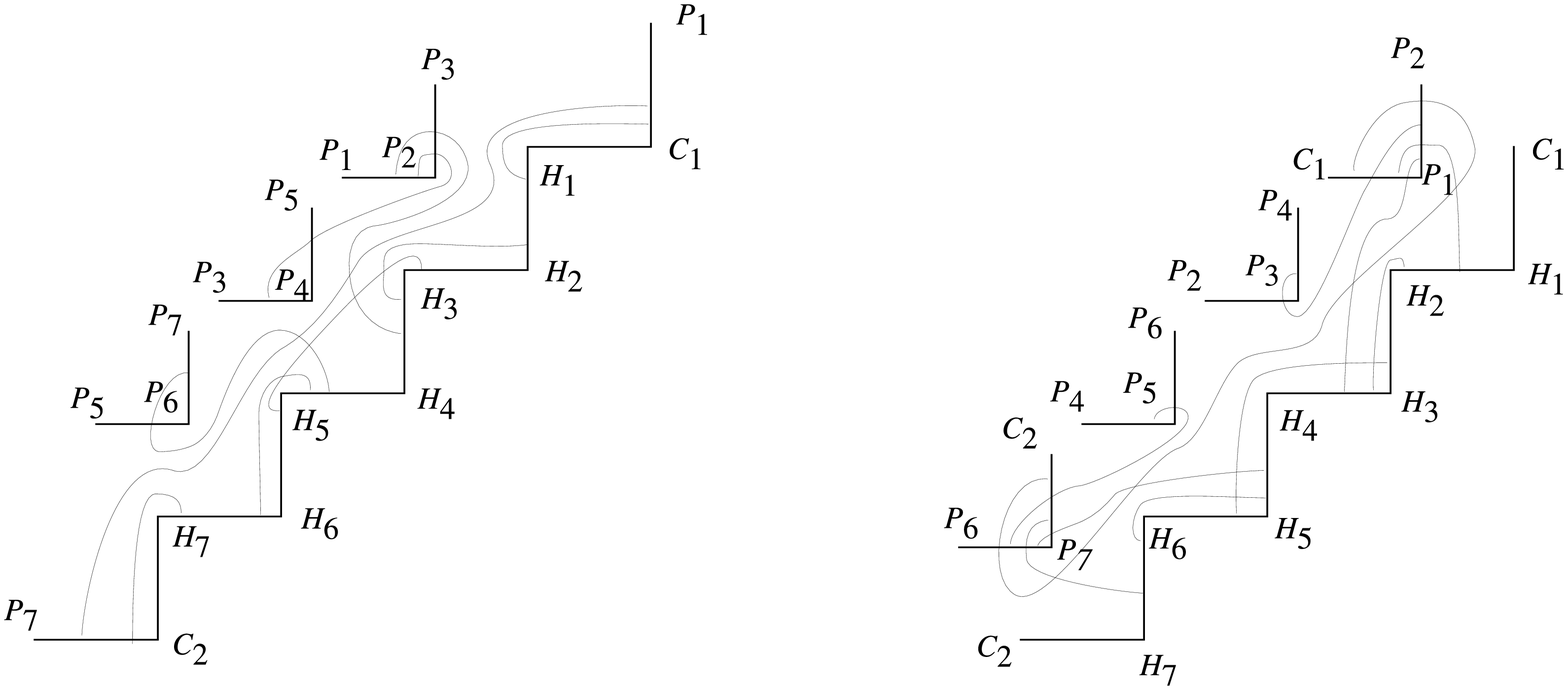 scaled 400}}
\vglue18pt
\centerline{The cycles for the height function for ${\rm DH}_{3,3}$ }
\vglue12pt

We record for later reference the following trivial consequence
of our using curves in the height function whose extremal lengths 
are conformal coordinates for the domains defined by the geometric 
coordinates.  

\nonumproclaim{Lemma 4.5.5} Let $X=\{\ogup, \ogdn\}$ be an 
orthodisk system in $\Delta$. Then $X$ is reflexive if and only 
if $\height(X)=0$.
\endproclaim

 4.6. {\it Properness of the height functions for} ${\rm DH}_{m,n}$.
In this section, we prove, modulo the postponed monodromy argument.

\nonumproclaim{Theorem 4.6.1}
The  height functions $\height$
from {\rm \S4.5} are proper.
\endproclaim

\demo{Proof}
To show that the  height functions from \S4.5 are proper, we need to prove
that for any
sequence of points in $\Delta$ converging to some boundary point, the
height goes to
infinity.  First, by Lemma~4.4.5, at least one
of the two orthodisks degenerates conformally.
We will now analyze the possible geometric degenerations. 

For this, we normalize the orthodisks in the sequence so that the finite edges
of the outer boundary of $\ogup$ has total length $1$. Such a scaling does
not affect the conformal structure.

Suppose first that the degeneration occurs with an inner sheet branch
point. As proven in Lemma~4.4.5, this implies that some branch point(s)
converges to a point on the outer sheet boundary.
Consider the  inner sheet with the smallest $k$ so that $P_k$ is a branch point of that sheet, for
which this degeneration happens.  \eject

If  $k=1$ and the edge $H_1C_1$ in $\ogup$ does not shrink to $0$, 
then $\ext \sigma \to 0$ in $\ogdn$
for $m<n$ and $\ext \nu \to 0$ in $\ogdn$, but the
corresponding extremal length remains positive in
$\ogup$. If the edge
$H_1C_1$ in
$\ogup$ does  shrink to
$0$, the monodromy argument can be applied to 
$\alpha_1$ ($m<n$) or to $\nu$ ($m=n$).

If $k>1$, the sheet containing $P_k$ is
separated from the outer sheet boundary by a $\beta$-cycle or (for a central inner sheet by the
cycle $\mu$) which will be pinched in one of the orthodisks but not the other, and which also will
not be pinched by other inner sheets degenerating. Hence the height of
this cycle will go to $\infty$.

Now let us assume that all inner sheets are bounded away from the outer sheet boundary and
that the geometric degeneration occurs only in the outer boundary sheet of $\ogup$. 

We begin by analyzing one execeptional case:
Suppose that the edge $C_1H_1$ in $\ogup$ shrinks to $0$.
The monodromy lemma, applied to the cycle $\alpha_1$ 
(in the $m<n$ case) or the cycle $\nu$
(in the $n=m$ case) proves the properness of the height function of this cycle.

Now suppose that the edge $C_1H_1$ does not shrink to $0$, 
and consider the smallest $k$ for which
the edge
$H_k H_{k+1}$ shrinks to $0$. The monodromy lemma can be 
applied to a neighboring cycle
$\alpha_k$ or $\alpha_{k\pm 1}$ This treats all $k$ 
but $k=n$ in the case $m<n$ and $k=1$ in the
case $m=n=2$.

Consider the case $m<n$ and suppose that the central 
edge $H_n H_{n+1}$ shrinks to $0$. This is
the only geometric degeneration occurring in the orthodisk. 
Hence we can apply the monodromy lemma to
$\rho$.  

Consider the case $m=n=2$ and suppose that the finite
edge $H_1 H_{2}$ shrinks to $0$. This is again
the only geometric degeneration occurring in the orthodisk.
Hence we can apply the monodromy lemma to
$\nu$.

This treats all possible cases, and the theorem is proven.
\enddemo

\pagegoal=50pc

 4.7. {\it A monodromy argument}.
In this section, we prove that the periods of orthodisks have incompatible
logarithmic
singularities in suitable coordinates and apply this to prove the monodromy theorem 4.1.2.
The main idea has already been used in \cite{WW}: to study the dependence of extremal lengths
of the geometric coordinates, it is necessary to understand the asymptotic dependence
of extremal lengths of the 
degenerating conformal polygons (which is classical and well-known, see
\cite{Oht}), and the asymptotic dependence of the geometric
coordinates of the degenerating conformal
polygons. This dependence is given by Schwarz-Christoffel maps which 
are well-studied
in many special cases. Especially, it is known that these maps possess  asymptotic
expansions in logarithmic terms. Instead of computing this expansion explicitly for the
two maps needed (which is possible but tedious), we use a monodromy argument to
show that the crucial logarithmic terms have a different sign for the two expansions.

Let $\Delta$ be a domain of dimension at least two defining geometric
coordinates for a pair
of orthodisks $X_1$ and $X_2$ corresponding to given formal \wei data as usual. Suppose
$\gamma$ is a cycle in the underlying conformal polygon which joins  edges
$P_1P_2$ and
$Q_1Q_2$ which are parallel (and hence nonadjacent) in the orthodisks. In our applications, $\gamma$
will be one of the $\alpha$-cycles used in the height function. Denote by
$R_1$ the vertex before  $Q_1$ and  by $R_2$
the vertex after   $Q_2$ and observe that by assumption, $R_2\ne P_1$ and
  $P_2\ne R_1$.
Introduce a second
cycle
$\beta$ which connects $R_1Q_1$ with $Q_2R_2$.
\figin{Zmono1}{600}
\vglue6pt
\centerline{Monodromy argument}
\vglue12pt

\pagegoal=48pc

We formulate the claim of Theorem 4.1.2 more precisely in the following two
lemmas:
 
\nonumproclaim{Lemma 4.7.1} Suppose that for a sequence $p_n\in \Delta$ with
$p_n\to p_0\in \partial\Delta${\rm ,} 
$\ext_{X_1(p_n)}(\gamma) \to 0$ and $\ext_{X_2(p_n)}(\gamma)\to 0$. Suppose furthermore that $\gamma$
is a cycle encircling a single edge which degenerates geometrically to $0$ as $n\to\infty$ in $X_1(p_n)$.
Then
$$\vert e^{1/\ext_{X_1(p_n)}(\gamma)}-e^{1/\ext_{X_2(p_n)}(\gamma)} \vert^
	2\to\infty .$$
\endproclaim

\nonumproclaim{Lemma 4.7.2} Suppose that for a sequence $p_n\in \Delta$ with
$p_n\to p_0\in \partial\Delta$ we have that
$\ext_{X_1(p_n)}(\gamma) \to \infty$ and $\ext_{X_2(p_n)}(\gamma)\to \infty$.
Suppose furthermore that $\gamma$
is a cycle with an endpoint 
on an edge which degenerates geometrically to $0$ as $n\to\infty$  in $X_1(p_n)$.
Then
$$\vert e^{\ext_{X_1(p_n)}(\gamma)}-e^{\ext_{X_2(p_n)}(\gamma)} \vert^
	2\to\infty .$$
\endproclaim

\demo{Proof} We first prove Lemma 4.7.1.

Consider the conformal polygons corresponding to the pair of orthodisks.
Normalize the punctures by M\"obius transformations so that
$$P_1=-\infty, \quad P_2=0, Q_1=\varepsilon, \quad  Q_2=1$$
for $X_1$ and
$$P_1=-\infty,  \quad P_2=0, Q_1=\varepsilon',  \quad Q_2=1$$
for $X_2$.
By the assumption of Lemma 4.7.1, we know that $\varepsilon, \varepsilon'\to 0$
as $n\to\infty$. We now apply the monodromy Corollary 4.7.5 below to the curve
$\varepsilon_0
e^{it}$ and conclude that either
$$\frac{\vert\per\beta(X_1)\vert}{\vert\per\gamma(X_1)\vert} +
\frac1\pi \log\varepsilon  \tag"(4.7.1a)"$$
is single-valued in $\varepsilon$ while
$$\frac{\vert\per\beta(X_2)\vert}{\vert\per\gamma(X_2)\vert} -
\frac1\pi \log\varepsilon' \tag"(4.7.1b)"$$
is single-valued in $\varepsilon'$, or the same statement holds for analogous quantities
with opposite signs. Without loss of generality we can treat the first case.

Now suppose that $\varepsilon'$ is real analytic (and hence single-valued)
in $\varepsilon$ near $\varepsilon=0$. Then use of the fact that
$X_1$ and $X_2$ are conjugate implies that
the absolute lengths of $\beta$ in $X_1$ and $X_2$ are equal, as are 
those of $\gamma$; hence
$$\frac{\vert\per\beta(X_1)\vert}{\vert\per\gamma(X_1)\vert}=\frac{\vert\per\beta(X_2)\vert}
{\vert\per\gamma(X_2)\vert} .$$
Thus
we see that, after subtracting (4.7.1b) from (4.7.1a),
$$\log(\varepsilon\varepsilon'(\varepsilon))$$
is single-valued in $\varepsilon$ near $\varepsilon=0$ which contradicts that
$\varepsilon, \varepsilon' \to 0$.

Now  Ohtsuka's extremal length formula states that 
for the current normalization of $X_1(p_n)$,
we have 
$$\ext(\gamma) = O\left(\vert\log \varepsilon\vert^{-1}\right)$$
(see Lemma 4.5.3 in \cite{WW}  and \cite{Oht}). We conclude
that
$$\vert e^{1/\ext_{X_1(p_n)}(\gamma)}-
e^{1/\ext_{X_2(p_n)}(\gamma)}\vert =O\left(\frac1{\varepsilon}-\frac1{\varepsilon'}\right)
$$
which goes to infinity,
since we have shown that $\e$ and $\e'$ tend to zero at 
different rates.  This proves Lemma 4.7.1.

The proof of Lemma 4.7.2 is very similar: For convenience, we normalize the
points of
the punctured disks such that

$$P_1=-\infty,\quad P_2=0, Q_1=1, \quad Q_2=1+ \varepsilon$$
for $X_1$ and
$$P_1=-\infty,\quad P_2=0,Q_1=1,\quad Q_2=1+\varepsilon'$$
for $X_2$.

By the assumption of Lemma 4.7.2, we know that $\varepsilon, \varepsilon'\to 0$
as $n\to\infty$. We now apply the monodromy Corollary 4.7.5 below to the curve
$1+\varepsilon_0
e^{it}$ and conclude that
$$\frac{\per\gamma(X_1)}{\per\beta(X_1)} + \frac1\pi \log\varepsilon$$
is single-valued in $\varepsilon$ while
$$\frac{\per\gamma(X_2)}{\per\beta(X_2)} - \frac1\pi \log\varepsilon'$$
is single-valued in $\varepsilon'$.
The rest of the proof is identical to the proof of Lemma 4.7.1.
\enddemo

To prove the needed Corollary 4.7.5, we need
asymptotic expansions of the extremal length in terms of the geometric
coordinates of
the orthodisks. Though not much is known explicitly about extremal lengths
in general,
for the chosen cycles we can reduce this problem to an asymptotic control of
Schwarz-Christoffel integrals. Their monodromy properties allow us to
distinguish their
asymptotic behavior by the sign of logarithmic terms.

We introduce some notation:
Suppose we have an orthodisk such that the angles at the vertices
alternate between $\pi/2$ and $-\pi/2$ modulo $2\pi$. Consider the
Schwarz-Christoffel map
$$F:z \mapsto \int_i^z (t-t_1)^{a_1/2}\cdot\ldots\cdot(t-t_n)^{a_n/2} $$
(see \S3.2)
from a conformal polygon with vertices at $t_i$ to this orthodisk. Choose
four  distinct vertices 
$t_{i},t_{i+1}, t_j, t_{j+1}$ so that $j\equiv i \pmod2$, ensuring that the edges
$t_{i}t_{i+1}$ and $t_{j}t_{j+1}$ are parallel in  the orthodisk geometry.
(See the figure below.)
Introduce a cycle $\gamma$ in the upper
half-plane connecting edge
$(t_{i},t_{i+1})$ with edge $(t_j, t_{j+1})$ and denote by $\bar\gamma$ the
closed cycle obtained from
$\gamma$ and its mirror image across the real axis. Similarly, denote by
$\beta$ the cycle connecting
$(t_{j-1}t_j)$ with $(t_{j+1}t_{j+2})$ and by $\bar\beta$ the cycle
together with its mirror image.
\figin{Zmono2}{600}
\vglue6pt
\centerline{Cycles for analytic continuation }
\vglue12pt

 Now consider the Schwarz-Christoffel period integrals
$$\align
F(\gamma) &= \frac12 \int_\gamma
(t-t_1)^{a_1/2}\cdot\ldots\cdot(t-t_n)^{a_n/2}, \\
F(\beta) &= \frac12 \int_\beta
(t-t_1)^{a_1/2}\cdot\ldots\cdot(t-t_n)^{a_n/2} ,
\endalign$$
as  multi-valued functions depending on the
{\it now complex} parameters $t_i$.

\nonumproclaim{Lemma 4.7.3} Under analytic continuation of $t_{j+1}$ around
$t_j$ the periods change their values as shown\/{\rm :}\/
$$\align
F(\gamma) &\rightarrow F(\gamma) + 2 F(\beta) , \\
F(\beta) &\rightarrow F(\beta )
\endalign$$
\endproclaim

\demo{Proof} The proof is the same as  \cite{WW, Lemma 4.4.1}:
the path of analytic continuation of $t_{j+1}$ around
$t_j$ gives rise to an isotopy of $\bbc$ which moves $t_{j+1}$ along this path.
This isotopy drags 
$\beta$ and
$\gamma$  to new cycles $\beta'$ and $\gamma'$.
 
Because the curve $\b$ is defined to surround $t_j$ and $t_{j+1}$, the
analytic continuation merely returns $\b$ to $\b'$. Thus,
because $\beta'$ equals $\beta$, their periods are also equal. 
On the other hand, the curve $\g$ is not equal to 
the new `dragged' curve $\g'$. To see this, note that
the period of $\gamma'$ is obtained by developing the flat  structure of the
doubled orthodisk along $\gamma'$. To compute
this flat structure, observe 
the crucial fact that the angles at the orthodisk 
vertices are either $\pi/2$ or
$-\pi/2$, modulo $2\pi$; thus the angles of the
doubled orthodisk equal $\pi$,  modulo $2\pi$. Thus the
arc $\gamma'$ develops into the union of the
arc $\gamma$ with the arcs $\beta$ and $\overline{\beta}$; in
particular, we see that the period of $\gamma'$ equals the 
period of $\gamma$ plus twice
the period of $\beta$.
\enddemo

Now denote  $\delta:= t_{j+1}-t_j  $ and fix all $t_i$ other than
$t_{j+1}$, which  we regard as the independent variable.

\nonumproclaim{Corollary 4.7.4} The function
$F(\gamma)-\frac{\log\delta}{\pi i} F(\beta)$ is single\/{\rm -}\/valued and
holomorphic near
$\delta=0$.
\endproclaim

\demo{Proof}  By definition, the function is locally holomorphic near $\delta=0$. By
Lemma~4.7.3  it is single-valued.  
\enddemo

Now, for the properness argument, we are interested
in the geometric coordinates --- these
are the absolute values of the periods.
More precisely we are interested in
$$\|F(\gamma) \| := \vert \ree F(\gamma) \vert + \vert \im
F(\gamma) \vert .$$ We translate
the above statement about periods into a statement about their respective absolute values.

We consider  two conjugate orthodisks parametrized by Schwarz-Christoffel
maps $F_1$ and $F_2$
defined on the same conformal polygon. Recall that we have constructed $\beta$ and $\gamma$ to be 
purely horizontal or vertical.

\nonumproclaim{Corollary 4.7.5} \hglue-5pt Either $\| F_1(\gamma)\|  -
\frac{\log\delta}\pi \|
F_1(\beta) \|$ or $\| F_1(\gamma)\|  + \frac{\log\delta}\pi
\|
F_1(\beta)\|$ is real analytic in $\delta$ for $\delta=0$. In the first
case{\rm ,}
$\| F_2(\gamma)\|  + \frac{\log\delta}\pi \|
F_2(\beta) \|$ is real analytic in $\delta${\rm ,} and in the second
$\| F_2(\gamma)\|  - \frac{\log\delta}\pi \|
F_1(\beta) \|$ is.
\endproclaim

{\it Proof}.
Recall from the tedious \S4.1 that the above periods
are linear combinations of the
geometric coordinates where the coefficients are just signs, 
in the sense of being 
elements of $\{1,-1,i, -i\}$.
Now by construction, $F_j(\beta)$ is purely real or purely imaginary;
moreover, the direction of $F_j(\beta)$ is $\pm i$ times the
direction of $F_j(\gamma)$. This,
together with Corollary 4.7.4, implies the first claim. Next,
note that (see for example the table in \S3.4) if we turn left at 
a vertex in $X_1$, we will turn right at the corresponding vertex
in the conjugate orthodisk $X_2$, and vice versa. Thus, if 
the directions of the corresponding edges for $\gamma$
in corresponding orthodisks
differ by $+i$, then the directions of the 
corresponding edges for $\beta$ will differ by a $-i$, and vice versa. 
This implies the second claim. 
\hfill\qed

\vglue-6pt
\section{The gradient flow}
\vglue-6pt

 5.1. {\it Overall strategy}.
 
\vglue6pt {5.1.1}. In this section we continue the proof of the
existence portions of  the main theorems, Theorem A and
Theorem B(i).  In the previous sections, we assigned to a
configuration $\SC$ (the ones $\SC=(UD)^kU$ and
$\SC={\rm DH}_{m,n}$ were of principal interest) a moduli space
$\Del=\Del_{\SC}$ of pairs of conformal structures $\{\ogup,
\ogdn\}$ equipped with geometric coordinates $\vec{\bold
t}=(t_i,\ldots,t_l)$.

In the last section, we defined a height function $\height$ on
the moduli space $\Delta$ and proved that it was a proper
function: as a result, there is a critical point  for the
height function in $\Del$, and our overall goal in the  next
pair of sections is a proof that one of these critical points
represents a reflexive orthodisk system in $\Del$, and hence,
by Theorem~3.3.5, a minimal surface of the configuration
$\SC$. Our goal in the present section is a description of
the tangent space to the moduli space $\Del$: we wish to
display how   infinitesimal changes in the geometric
coordinates $\vec{\bold t}$ affect the height function.  In
particular, it would certainly be sufficient for our purposes
to prove the statement:

\nonumproclaim{Model 5.1.1} If $\vec{\bold t}_0$ is not a reflexive
orthodisk system{\rm ,}  there is an element $V$ of the tangent
space $T_{\vec{\bold t}_0}\Del$ for which $D_V\height \ne 0$.
\endproclaim

This would then have the effect of proving that our critical
point for the height function is reflexive, concluding the
existence  parts of the proofs of the main theorem.

We do not know how to prove or disprove this model statement
in its full generality.  On the other hand, it is not
necessary for the  proofs of the main theorems that we do so.
Instead we will replace  this theorem by a pair of lemmas.
The lemmas each have two cases, with the division into cases
depending on whether, in the configuration $\SC = {\rm DH}_{m,n}$,
we have $n=m$ or $m<n$.

\nonumproclaim{Lemma 5.1.2} Consider a configuration
$\SC = {\rm DH}_{m,n}$. 

{\rm (i)} Suppose $n=m$ and $\SY\subset\Del$ 
is a real two\/{\rm -}\/dimensional subspace of $\Del$ which is defined 
by the equations ${\height(t_i)=0}$ for all coordinates $t_i$
with only two exceptions {\rm (}\/say $t_1$ and $t_2$\/{\rm ).}
Suppose also that the coordinates $t_1$ and $t_2$
refer to cycles which encircle neighboring edges{\rm ,} at most one
of which is a finite edge.  Then there is a sublocus $\SY^*$
along which ${\height(t_{i_1})=0}${\rm ,} and along that 
one\/{\rm -}\/dimensional sublocus{\rm ,} if
$\vec{\bold t}_0\in\SY^*$ has positive height {\rm (}\/i.e.\
$\height(\vec{\bold t}_0)>0${\rm ),}   there is an element $V
\in T_{\vec{\bold t}_0}\SY^*$ of the tangent space
$T_{\vec{\bold t}_0}\SY^*$ for which $D_V\height \ne 0$.
\vglue4pt
{\rm (ii)} If $m<n$ and if $\SY\subset\Del$ is a real
one\/{\rm -}\/dimensional subspace of $\Del$ which is defined  by the
equations ${\height(t_i)=0}$ for all coordinates $t_i$ with
only one exception {\rm (}\/say $t_1${\rm ),} and that exception refers to
a cycle which is not a finite edge bully. Then if
$\vec{\bold t}_0\in\SY$ has positive height{\rm ,}  then there is
an element $V \in T_{\vec{\bold t}_0}\SY$  of the tangent
space $T_{\vec{\bold t}_0}\SY$ for which $D_V\height \ne 0$.
\endproclaim

\nonumproclaim{Lemma 5.1.3} Consider a configuration
$\SC = {\rm DH}_{m,n}$.

{\rm (i)} If $n=m${\rm ,} then for every configuration
$\SC = {\rm DH}_{m,n}${\rm ,} there is a nonsingular analytic subspace
$\SY \subset \Del=\Del_{\SC}${\rm ,} for which 
$\SY = \{\height(t_i)=0 | i \ne i_1, i_2\}${\rm ,} for some choice
of $i_1${\rm ,} $i_2${\rm ,} where $t_{i_1}$ and $t_{i_2}$ refer to
cycles encircling neighboring edges{\rm ,} of which at most one
is a finite edge. 

{\rm (ii)} If $m<n${\rm ,} then for every configuration $\SC = {\rm DH}_{m,n}${\rm ,}
there is a nonsingular analytic subspace $\SY \subset
\Del=\Del_{\SC}${\rm ,} for which 
$\SY = \{\height(t_i)=0 | i \ne i_1\}${\rm ,} for some choice of
$i_1$ and the cycle referred to by $t_{i_1}$ is not a finite edge bully.
\endproclaim

A cycle $\G$ is a {\it finite edge bully} for the previous lemma if each
finite edge $E$ of one of the orthodisks has a neighbor $E'$ so that $\G$ 
foots on $E$ but not $E'$ or vice versa.

We introduce this concept because when $m+n$ is small, the
technical requirements of the proof (see Lemma~5.3.1 and
Definition~5.3.0) contrast with a shortage of sides; when
$m+n\ge5$ in case~(ii), no cycle
is a finite edge bully, and these hypotheses are vacuous.  In our proof
of Lemma~5.1.3,
we will produce a cycle that is not a finite edge bully.

Given these lemmas, the proofs of Theorems A and B(i) are
straightforward.

\vglue6pt
{\it Proofs of Theorem~{\rm A} and Theorem~{\rm B(i)}}. Consider the locus
$\SY$ (or the sublocus $\SY^*$) guaranteed by Lemma~5.1.3. By
Theorem~4.5, the height function $\height$ is proper on
$\SY$ (resp. $\SY^*$), and has a
critical point on $\SY$. By Lemma~5.1.2, this critical point
represents a point where $\height=0$, i.e., a reflexive
orthodisk by Lemma~4.5.5.
\vglue4pt

The proof of Lemma 5.1.2 occupies the current section while the
proof of Lemma 5.1.3 is given in the following section.\eject

  {5.1.2}. Let us discuss informally the proof of
Lemma~5.1.2. Because angles of corresponding vertices in the
$\ogup\leftrightarrow\ogdn$ correspondence sum to
$0\bmod2\pi$, the orthodisks locally topologically fit
together along corresponding edges, so that conjugacy of
orthodisks requires corresponding edges to move in different
directions: if the  edge $E$ on $\ogup$ moves ``out'', the
corresponding edge $E^*$ on $\ogdn$  moves ``in'', and
vice versa (see the figures below). Thus  we expect that if
$\g$ has an endpoint on $E$, then one of the extremal lengths
of $\gamma$ decreases, while the other extremal length of
$\gamma$ on the other orthodisk  would increase: this will
force the  change of the height $\height(\gamma)$ of
$\gamma$ to have a definite sign, as desired.  This is the 
intuition behind Lemma~5.2.1; a rigorous argument requires us
to actually compute derivatives of relevant extremal lengths
using the formula (2.2) and we will need to avoid some
technical difficulties involving singularities of the
holomorphic quadratic differentials arising in (2.2). 
(We also will need to do some preliminary asymptotic analysis
and an intermediate value theorem argument to find the 
sublocus $\SY^*$ along which to apply this argument
in case (i).) We do
this by displaying, fairly explicitly, the deformations of
the orthodisks (in local coordinates on $\ogup/\ogdn$) as well
as the differentials of extremal lengths, also in 
coordinates. After some preliminary notational description in 
\S5.2, we do most of the computing in \S5.3.  Also in 
\S5.3 is the key technical lemma, which relates the formalism 
of formula (2.2), together with the local coordinate descriptions
of its terms, to the intuition we just described.

\vglue4pt 5.2.  {\it Deformations of ${\rm DH}_{1,1}$}.
We will prove Lemma~5.1.2 by first considering the special
case of $\SC = {\rm DH}_{1,1}$.  We will then note that this case,
  defined in terms of precisely two curve systems and
involving all the types of geometric situations encountered in
all of the configurations $\SC$, is only slightly different
from the most general case: the arguments we  give here
extend immediately to prove the lemma in general, even 
though the exposition, in this specific low-dimensional case,
can be much more concrete.

Here, we will describe the notation and
results for this special case where  $\SC={\rm DH}_{1,1}$.

We begin by recalling the notation of \S4.2: there the
geometric coordinates were labeled $\Del = (y,b,g)$ corresponding
to the curve systems labeled yellow, blue, and green (see
the figures in \S4.2). We introduced a fourth curve system, mauve,
which we used with blue to define the height function.

\nonumproclaim{Theorem 5.2.1} There exists a reflexive orthodisk
system for the configuration ${\rm DH}_{1,1}$.
\endproclaim

{\it Proof}. In this case, the statement of Lemma~5.1.3(i) is vacuous
as the moduli space $\Del=\Del_{1,1}$ is already two-dimensional.
So we are left to prove Lemma~5.1.2 in this case. The plan is
simple, as we consider a path $\Gamma \subset \Del$ which has two
limits: towards one end $E$, the size of the period 
$H_1H_2-C_1P_1$ (and $H_3H_2-P_3C_2$) decays to zero, and 
towards the other end $E'$, the period $P_1P_2-H_1C_1$ (and
$P_2P_3-C_2H_3$) decays to zero. We
\vglue4pt
{\elevensc Claim 5.2.2}. $$\sgn(\ext_\ogup({\rm blue})-\ext_\ogdn({\rm blue}))\neq
\sgn(\ext_\ogup({\rm blue})-\ext_\ogdn({\rm blue})).$$
 
Assuming the claim, we continue with the proof of 
Theorem~5.2.1. As the sign of $\ext_\ogup({\rm blue})-\ext_\ogdn({\rm blue}))$
changes along the arc $\Gamma=\{\Gamma(s)\}$, we know that there is a 
point $\Gamma(s_0)$ at which 
$\ext_\ogup({\rm blue})-\ext_\ogdn({\rm blue}))$ changes sign. Of course,
this function $\ext_\ogup({\rm blue})-\ext_\ogdn({\rm blue}))$ is the difference
of real analytic regular functions of the coordinates on
$\Delta$ (see expanded discussion of the analyticity of the
curve at the end of the proof of 
Lemma~6.2), and so there is an analytic path through
$\Gamma(s_0)$ along\break which 
$\ext_\ogup({\rm blue})=\ext_\ogdn({\rm blue}))$.
In keeping with the notations of\break Lemma~5.1.2 and 5.1.3,
we denote that analytic path by $\SY^*$.

Our next task is to find a point on $\SY^*$ for which 
$\SH({\rm mauve})=0$: at such a point, we would have
$\SH({\rm mauve})=0=\SH({\rm blue})$. As $\SY^*$ is one-dimensional,
and mauve is not a finite-edge bully, we in effect
prove Lemma~5.1.2(i). The idea is to make an
``infinitesimal push'' along the 
edges $H_1H_2$ (and $H_2H_3$); we push 
$H_1H_2$ (and $H_2H_3$) ``into''
$\ogup$ and then the requirement that 
$\ogup$ and $\ogdn$ should remain conjugate then 
forces us to make an equal push of
$H_1H_2$ (and $H_2H_3$) ``out of''
$\ogdn$. By the principle of domain monotonicity
of extremal length, we expect the change 
in $\ogup$ to increase the extremal length
$\ext_\ogup({\rm mauve})$ and we expect the 
change in $\ogdn$ to decrease the extremal length
$\ext_\ogdn({\rm mauve})$.  The rigorously defined 
formulae and the technical Lemma~5.3.1 will
support this intuition (indeed, will show that
it is true to first order), and we find that
if $\SH({\rm mauve}) \ne 0$ at a point 
$\vec{\bold t}_0\in\SY^*$, then there is an element $V
\in T_{\vec{\bold t}_0}\SY^*$ of the tangent space
$T_{\vec{\bold t}_0} \SY^*$ for which $D_V\height \ne 0$.
This completes the proof of Lemma~5.1.2 in this case,
as well as Theorem~5.1.2, up to the proof of the
claim.

\centerline{\BoxedEPSF{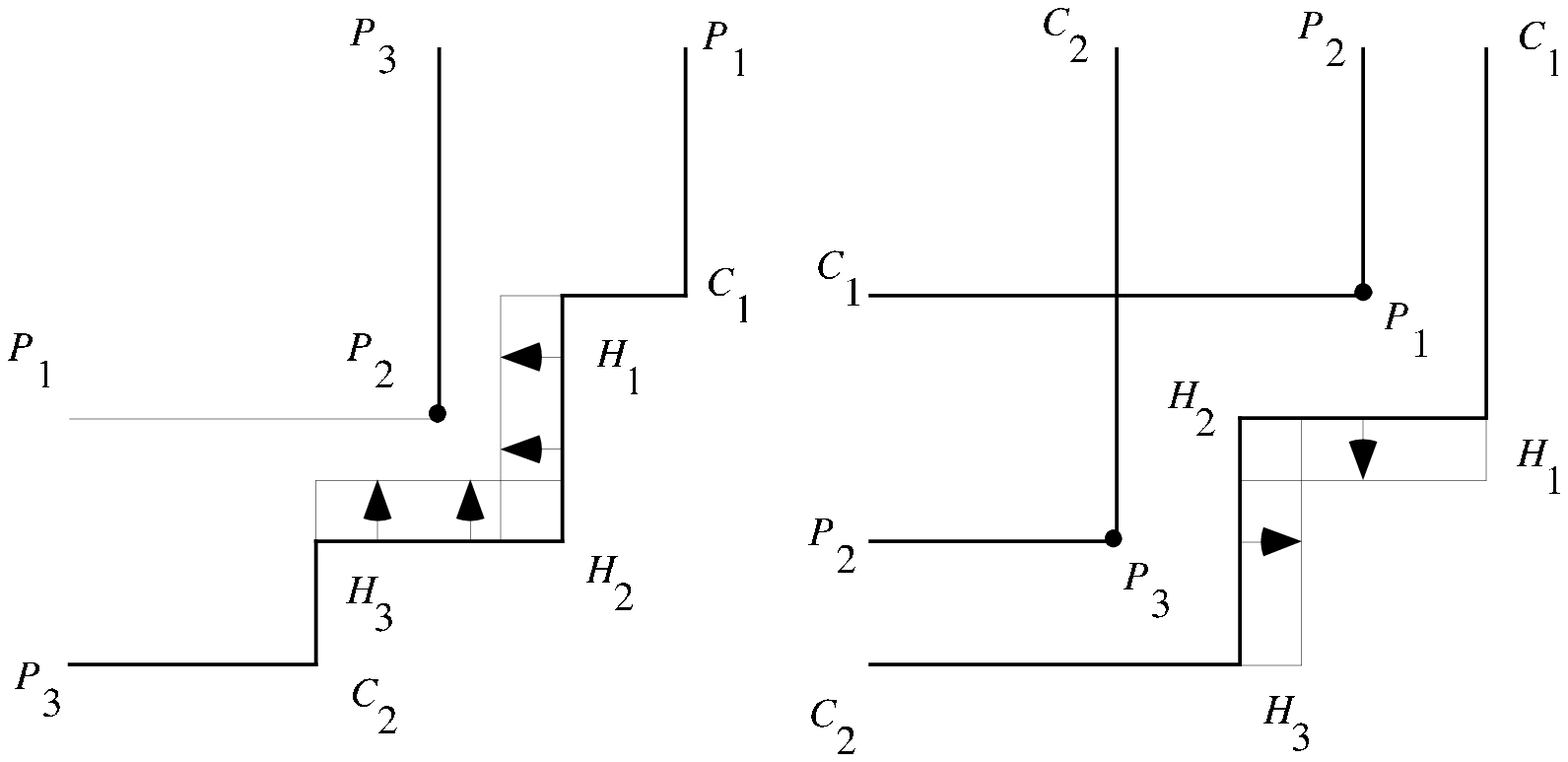 scaled 600}}
 
\centerline{Pushing edges along the locus $\SY^*$}
 \eject

{\it Proof of the Claim~{\rm 5.2.2}}. It is evident that as the
period $P_1P_2-H_1C_1$ tends to zero (with the
period $H_1H_2-C_1P_1$ remaining bounded from below),
the extremal length $\ext_\ogup({\rm blue})$ tends to a 
nonzero finite number, while the extremal length
$\ext_\ogdn({\rm blue})$ tends to infinity.

On the other hand, if we first send the period
$P_2P_3-H_1H_2$ to zero and then allow
the period $H_1H_2 - P_1C_1$ to tend to zero, 
we see that $\ext_\ogup({\rm blue})$ is equal to the
extremal length of ${\rm blue}$ on the quadrilateral component
$P_1-P_2=H_2-H_1-C_1-P_1$ of a noded surface,
while $\ext_\ogdn({\rm blue})$ is given in terms of 
the extremal length of ${\rm blue}$ on the
pentagonal surface 
$P-2-H_2-H_1-C_1-P_1-P_2$. Yet the asymptotics
of these extremal lengths of these 
explicit shapes are readily (if tediously)
computable, and we find that
$\ext_\ogup({\rm blue})$ decays to zero at a slower
rate than $\ext_\ogdn({\rm blue})$ as we send the
period $H_1H_2 - P_1C_1$ to zero.  Thus by allowing
$\Gamma$ to have ends near these two degenerate loci,
we obtain the claim.
\vglue12pt

 5.3. {\it Infinitesimal pushes}.
We need to formalize the previous discussion. As always we are
concerned with relating the Euclidean geometry of the 
orthodisks (which corresponds directly with the periods of the
\wei data) to the conformal data of the domains $\ogup$ and
$\ogdn$. From the discussion above, it is clear that the
allowable  infinitesimal motions in $\Del$, which is
parametrized in terms of the Euclidean geometry of 
$\ogup$ and $\ogdn$,  are given by infinitesimal changes in
lengths of finite sides or in distances between sheets, with
the changes being done simultaneously on $\ogup$ and $\ogdn$
to preserve conjugacy. The link to the conformal geometry is
the formula (2.2): a motion which infinitesimally transforms 
$\ogup$, say, will produce an infinitesimal change in the 
conformal structure. Tensorially, this tangent vector to the
moduli space of conformal structures is represented by a
Beltrami differential. Later, formula (2.2) will be used,
together with knowledge of the cotangent vectors $ d
\ext_\ogup(\cdot)$ and $ d \ext_\ogdn(\cdot)$, to determine
the derivatives of the relevant extremal lengths, hence the
derivative of the height.

\vglue12pt   5.3.1. {\it Infinitesimal pushes}.
Here we explicitly compute the effect of infinitesimal pushes
of  certain edges on the extremal lengths of relevant cycles. 
This is done by explicitly displaying the infinitesimal
deformation and then using this formula to compute the sign
of the derivative of the extremal lengths, using formula
(2.2). There will be four different cases to consider, all
of  which are apparent in the ${\rm DH}_{1,1}$ orthodisks. This is
why we can,  without loss of generality, confine the
computations to this case. To be concrete, we will choose for
each case either the orthodisk $\ogup$ or the orthodisk
$\ogdn$.

The three sides  $H_1C_1$, $C_1P_1$ and $H_2H_1$ are geometrically
different and require separate treatments. Thus we break
our discussion up into cases.\eject

\noindent \hangindent=42pt\hangafter=1{Case A.}  Finite {\it noncentral} edges 
(see also \cite{WW}). Example: $H_1C_1$ in $\ogup$.
 
\vglue4pt
\noindent \hangindent=42pt\hangafter=1{Case B.} Infinite edge whose finite vertex is unbranched.
Example: $H_1C_1$ in $\ogdn$.
\vglue4pt
\noindent \hangindent=42pt\hangafter=1{Case C.} Infinite edge whose finite vertex is a branch point
off the symmetry line.
Example: $C_1P_1$ on $\ogdn$.
\vglue4pt
\noindent \hangindent=42pt\hangafter=1{Case D.} An edge (finite or infinite)
and its symmetric side meet in a
corner. Example: $H_3H_2$ and $H_2H_1$.
\vglue4pt

For each case there are two subcases, which we can describe as
depending on whether the given sides are horizontal or
vertical. The distinction is, surprisingly, a bit important,
as together with the fact that we do our deformations in
pairs, it provides for an important cancellation of
(possibly) singular terms affecting the formulation of 
Lemma~5.3.1. We defer this point for later (see Remark 5.3.3),
while here we  begin to calculate the relevant Beltrami
differentials in the  cases.

Also, each infinitesimal motion might require two different types
of cases, depending on whether the edge we are deforming 
on $\ogup$ corresponds on $\ogdn$ to an edge of the same type
or a different type. We will thus compute the Beltrami differentials
only for a single domain, either $\ogup$ or $\ogdn$.
\vglue6pt
Case A. Here the computations are  analogous to those that we found in 
[WW]; they differ only in orientation of the boundary of the orthodisk.
\figin{Zcasea}{600}
\vglue6pt
\centerline{Beltrami differential computation --- Case A}
\vglue12pt

We first consider the case of a horizontal finite side; as in
the figure above, we see that the neighborhood of the
horizontal side of the orthodisk in the plane naturally
divides into six regions which we label $R_1,\ldots,R_6$. Our
deformation $f_\e = f_{\e, b, \d}$ differs from the identity
only in such a neighborhood,  and in each of the six regions,
the map is affine.
In fact we have a two-parameter family of these deformations, all
of which have the same infinitesimal effect, with the parameters
$b$ and $\d$ depending on the dimensions of the supporting neighborhood, 

$$
{\ninepoint
\hbox{$f_\e(x,y)$} = \cases
\(x,\e +\f{b-\e}b y\),   &\{-a\le x\le a, 0\le y\le b\}=R_1\\
\(x,\e +\f{b+\e}b y\),   &\{-a\le x\le a, -b\le y\le0\}=R_2\\
\(x,y+\f{\e+\f{b-\e}b y-y}\d(x+\d+a)\),   &\{-a-\d\le x\le-a,
0\le y\le b\} = R_3\\
\(x,y-\f{\e+\f{b-\e}b y-y}\d(x-\d-a)\),   &\{a\le x\le a+\d,
0\le y\le b\} = R_4\\
\(x,y+\f{\e+\f{b+\e}b y-y}\d(x+\d+a)\),   &\{-a-\d\le x\le-a,
-b\le y\le0\} = R_5\\
\(x,y-\f{\e+\f{b+\e}b y-y}\d(x-\d-a)\),   &\{a\le x\le a+\d,
-b\le y\le0\} = R_6\\
(x,y)   &\text{otherwise} 
\endcases}\tag"(5.1a)"
$$
where we have defined the regions $R_1,\dots,R_6$ within the
definition of $f_\e$. Also note that here the orthodisk contains the arc
$\{(-a,y)\mid0\le y\le b\}\cup\{(x,0)\mid-a\le x\le a\}\cup
\{(a,y)\mid-b\le y\le0\}$.
Let $E$ denote the edge being pushed, defined above as $[-a,a]\times\{0\}$.

Of course $f_\e$ differs from the identity only on a neighborhood
of the edge $E$, so that $f_\e$ takes the symmetric orthodisk to
an asymmetric orthodisk. 
We next modify $f_\e$ in a neighborhood of the reflected
(across the $y=-x$ line) segment $E^*$ in an analogous way with a
map $f^*_\e$ so that $f^*_\e\circ f_\e$ will preserve
the symmetry of the orthodisk. 

Our present conventions are that the edge $E$ is horizontal; this forces
$E^*$ to be vertical and we now write $f^*_\e$ for such a
vertical segment; this is a straightforward extension of the
description of $f_\e$ for a horizontal side, but we present the
definition of $f^*_\e$ anyway, as we are crucially interested in
the signs of the terms. So set
$${\ninepoint
\hbox{$f^*_\e$} = \cases
\(-\e +\f{b-\e}bx, y\),   &\{-b\le x\le0, -a\le y\le a\}=R^*_1\\
\(-\e +\f{b+\e}bx, y\),   &\{0\le x\le b, -a\le y\le a\}=R^*_2\\
\(x - \f{-\e+\f{b-\e}b x-x}\d(y-\d-a), y\),   &\{-b\le x\le0, a\le
y\le a+\d\}=R^*_3\\
\(x + \f{-\e+\f{b-\e}b x-x}\d(y+\d+a), y\),   &\{-b\le x\le0,
-a-\d\le y\le-a\}=R^*_4\\
\(x - \f{-\e+\f{b+\e}b x-x}\d(y-\d-a), y\),   &\{0\le x\le b, a\le
y\le a+\d\}=R^*_5\\
\(x + \f{-\e+\f{b+\e}b x -x}\d(y+\d+a), y\),   &\{0\le x\le b, -a-\d\le
y\le-a\}=R^*_6\\
(x,y)   &\text{otherwise} .
\endcases} \tag"(5.1b)"
$$
Note that under the reflection across the line $\{y=-x\}$, the
region $R_i$ gets taken to the region $R_i^*$.
 
Let $\nu_\e = \f{\(f_\e\)_{\bar z}}{\(f_\e\)_z}$
denote the Beltrami differential of $f_\e$, and set
$\dot\nu=\f d{d\e}\bigm|_{\e=0}\nu_\e$. Similarly, let $\nu^*_\e$
denote the Beltrami differential of $f^*_\e$, and set $\dot\nu^*=\f
d{d\e}\bigm|_{\e=0}\nu^*_\e$. Let $\dot\mu=\dot\nu+\dot\nu^*$. Now
$\dot\mu$ is a Beltrami differential supported in a bounded domain
in one of the domains $\ogup$ or $\ogdn$. 
We begin by observing that it is easy to compute that
$\dot\nu=[\f d{d\e}\bigm|_{\e=0}\(f_\e\)]_{\bar z}$ 
evaluates near $E$
to
$${\ninepoint
\dot\nu = \cases
\f1{2b},   &z\in R_1\\
-\f1{2b},   &z\in R_2\\
\f1{2b}[x+\d+a]/\d+i\(1-y/b\)\f1{2\d}=\f1{2b\d}(\bar z+\d+a+ib),
&z\in R_3\\
-\f1{2b}[x-\d-a]/\d-i\(1-y/b\)\f1{2\d}=\f1{2b\d}(-\bar z+\d+a-ib),
&z\in R_4\\
-\f1{2b}[x+\d+a]/\d+i\(1+y/b\)\f1{2\d}=\f1{2b\d}(-\bar z-\d-a+ib),
&z\in R_5\\
\f1{2b}[x-\d-a]/\d-i\(1+y/b\)\f1{2\d}=\f1{2b\d}(\bar z-\d-a-ib),
&z\in R_6\\
0  &z\notin\supp(f_\e-\id) .
\endcases}\tag"(5.2a)"
$$
We further compute
$$
\dot\nu^* = \cases
-\f1{2b},   &R^*_1\\
\f1{2b},   &R^*_2\\
\f1{2b\d}(i\bar z-\d-a+bi)    &R^*_3\\
\f1{2b\d}(-i\bar z-\d-a-bi)   &R^*_4\\
\f1{2b\d}(-i\bar z+\d+a+bi)   &R^*_5\\
\f1{2b\d}(i\bar z+\d+a-bi)    &R^*_6 .
\endcases\tag"(5.2b)"
$$

\vglue6pt
Cases B, C.
For both infinite sides the deformations are the same, here
represented in the horizontal case. (We 
defer the case
of infinite sides that meet along the symmetry line $\{y=-x\}$
until we treat Case D, where it will fit more naturally.)
\vglue-12pt
\figin{Zcasebc}{600}
 
\centerline{Beltrami differential computation --- Cases B and C}
$$
f_{\e;b,\d}(x,y) = \cases
\(x,\e +\f{b-\e}b y\),   &R_1=\{x \le 0, 0 \le y \le b\}\\
\(x,\e +\f{b+\e}b y\),   &R_2=\{x \le 0, -b \le y \le 0\}\\
\(x,y-\f{\e+\f{b-\e}b y-y}\d(x-\d)\),   
&R_3= \{0 \le x \le \d,0 \le y \le b\} \\
\(x,y-\f{\e+\f{b+\e}b y-y}\d(x-\d)\),   
&R_4= \{0 \le x \le \d,-b \le y \le 0\}\\
(x,y)   &\text{otherwise}.
\endcases
$$

Thus our infinitesimal Beltrami differential is given by

$$
\dot\nu_{b,\d} = \cases
\f1{2b},   &z\in R_1\\
-\f1{2b},   &z\in R_2\\
\f1{2b\d}(-\bar z+\d-ib),
&z\in R_3\\
\f1{2b\d}(\bar z-\d-ib),
&z\in R_4\\
0  &z\notin \supp(F_{\e;b,\d}-\id) .
\endcases
$$

The formulas for a push of a vertical edge are analogous.
\vglue12pt
Case D. We have separated this case out for purely expositional 
reasons. We can imagine that the infinitesimal push that 
moves the pair of consecutive sides along the symmetry line
$\{y=-x\}$ is the result of a composition of a pair of 
pushes from Case A or from Case C; i.e., our diffeomorphism $F_{\e;b,\d}$
can be written  $F_{\e;b,\d}= f_{\e} \circ f_{\e}^*$,
where the maps differ from the identity in the {\it union}
of the supports of $\dot\nu_{b,\d}$ and 
$\dot\nu_{b,\d}^*$.

It is an easy consequence of the chain rule applied to this 
formula for  $F_{\e;b,\d}$ that the infinitesimal Beltrami differential
for this deformation is the sum $\dot\nu_{b,\d} + \dot\nu_{b,\d}^*$
of the
infinitesimal Beltrami differentials 
$\dot\nu_{b,\d}$ and 
$\dot\nu_{b,\d}^*$ defined in formula (5.2)
for Case A (even in a neighborhood of the vertex along the
diagonal where the supports of the differentials $\dot\nu_{b,\d}$ and 
$\dot\nu_{b,\d}^*$ coincide).

\vglue12pt  5.3.2. {\it Derivatives of extremal lengths}. In this section, we combine the computations of $\dot\nu_{b,\d}$
with formula (2.2) (and its background in \S2)
and some easy observations on the nature of the quadratic
differentials $\Phi_{\mu}= \frac12 d \ext_{(\cd)}(\mu)\bigm|_{\cd}$ 
to compute the derivatives of extremal lengths under our
infinitesimal deformations of edge lengths.

We begin by recalling some background from \S2. If we are given
a curve $\gamma$, the extremal length of that curve on an orthodisk,
say $\ogup$, is a real-valued
$C^1$ function on the moduli space of that orthodisk.
Its differential is then a holomorphic quadratic differential
$\Phi_{\gamma}= \frac12d\ext_{(\cd)}(\gamma)\bigm|_{\ogup}$
on that orthodisk; the horizontal foliation of $\Phi_{\g}$ consists of
curves which connect the same edges in $\ogup$ as $\g$,
since $\Phi_{\g}$ is obtained as the pullback of the quadratic
differential $dz^2$ from a rectangle where the image of $\g$ connects
the opposite vertical sides.
We compute the derivative of the 
extremal length function using formula (2.2); i.e.
$$
\(d\ext(\gamma)\bigm|_{\ogup}\)[\nu] =
4\ree\int_{\ogup}\Phi_\gamma\nu .
$$

It is here where we find that we can actually compute the 
sign of the derivative of the extremal lengths, hence the
height function, but also encounter a subtle technical
problem. The point is that we will discover that, in suitable
circumstances, just the topology of the curve $\gamma$ on
$\ogup$ will determine the sign of the derivative on an edge
$E$, so we will be able to evaluate the sign of the integral
above, if we shrink the support of the Beltrami differential
$\dot\nu_{b,\d}$ to the edge by sending $b$, $\d$ to zero. 
(In particular, the sign of  $\Phi_\gamma$ depends precisely
on whether the foliation of $\Phi= \Phi_\gamma$ is parallel
or perpendicular to $E$, and on whether $E$ is horizontal or
vertical.) We then need to know two things: 1) that this
limit exists, and 2) that we may know its sign via
examination of the sign of $\dot\nu_{b,\d}$ and $\Phi_\gamma$
on the edge $E$. We phrase this as follows:

\demo{Definition {\rm 5.3.0}} $\Phi$ is admissable on an edge
$E$ if and only if  at a vertex $v$ of $E$ and in a conformal parameter
$z$ centered at $v$, we have $|z|\Phi(z)=o(1)$.
\enddemo

\nonumproclaim{Lemma~5.3.1} Let $\Phi$ be admissable for the edge
$E$. {\rm(1)}
$\lim_{b\to0,\d\to0}\ree\int\Phi\dot\nu$ exists{\rm ,} is finite
and nonzero. {\rm(2)} The foliation of $\Phi$ is either
parallel or orthogonal to the interior of the segment which is
$\lim_{b\to0,\d\to0}(\supp\dot\nu)${\rm ,} and {\rm(3)} The
expression $\Psi\dot\nu$ has a constant sign on that segment
$E${\rm ,} and the integral {\rm (2.2)} also has that {\rm (}\/same\/{\rm )} sign.
\endproclaim

In the statement of the lemma, the foliation refers to the horizontal
foliation of the holomorphic quadratic differential
$\Phi=\Phi_\gamma$, whose core curve is~$\gamma$.

This lemma provides the rigorous foundation for the intuition
described in the paragraph of strategy \S5.1.2.

\vglue12pt 5.3.3. {\it Proof of the technical lemma and a remark on
its limitations}.

\demo{Proof of Lemma~{\rm 5.3.1}}  
The foliation of $\Phi$, on say $\ogup$, lifts to a foliation
on the punctured sphere (which we will denote by $S_{Gdh}$), 
symmetric about the reflection about the equator. This
proves the second statement. The third statement follows from
the first (and from the above discussion of the topology of
the vertical foliation of $\Phi_{\g}$,  once we prove that
there is no infinitude coming from either the neighborhood of
infinity of the infinite edges or the regions $R_3$ and $R_4$
for the finite vertices. This finiteness will follow from the
proof of the first statement. Thus, we are left to prove the
first statement.

In the case of an edge with a finite vertex, we estimate the
contribution to $\ree\int\Phi\dot\nu$ of a neighborhood $N$
of a vertex from regions $R_1$ and $R_2$ (and $R^*_1$ and
$R^*_2$) as
$$
\split
|\ree\int\Phi\dot\nu|  &\le\int^b_0\int_{N\cap \partial\Omega_{Gdh}}\f1{2b}\cd
|\Phi|dxdy\\
&= \f1{2b}\int^b_0\int_{N\cap \partial\Omega_{Gdh}}o(x^2 + y^2)^{-1/2}dxdy
\endsplit
$$
which converges by the dominated convergence theorem.  
The other regions are analogous.

In the case of an edge with an infinite vertex, the argument
is somewhat easier, as the possible growth of $|\Phi|$ as we
tend towards the vertex has an {\it a priori\/} bound of rapid decay.
We compute this decay by observing that in the coordinate $w$
of the plane, the quadratic differential $\Phi$ has expansion
$\Phi=c_1/w +\{\text{h.o.t.\/}\}dw$, and the Schwarz-Christoffel map from the
half-plane to the orthodisk has integrand
$\om=\{c_0z^{-3/2}+\text{h.o.t.\/}\}dz$. It is then straightforward
that in the coordinate $z$ of the orthodisk, the quadratic
differential $\Phi$ decays faster than $|z|^{-2}$. \enddemo

Computations like those in the final paragraph of the above
proof prove the

\nonumproclaim{Proposition 5.3.2} The holomorphic quadratic
differential $\Phi$ is admissable on an edge $E$ if $\Phi$
relates to the vertices $v_1$ and $v_2$ of $E$ in one of
the following ways\/{\rm :}\/
\vglue4pt
\ritem{(i)} The orthodisk has an angle of $\pi/2$ of $v_i$.
\vglue4pt
\ritem{(ii)} The vertex $v_i$ is at infinity.
\vglue4pt
\ritem{(iii)} The orthodisk has an angle of $3\pi/2$ at
$v_i${\rm ,} and the foliations of $\Phi$ are either both
parallel to the edges incident to $v_i$ or both orthogonal
for those edges.
\vglue4pt
\ritem{(iv)} The vertex $v_i$ is a branch point where the
orthodisk makes an angle of $5\pi/2${\rm ,} and the holomorphic
quadratic differential has a singularity at $v_i$.
\vglue4pt \noindent 
In the last case {\rm (iv),} a singularity arises if and only if
the foliation is parallel to the edges incident to $v_i$
and a leaf emanates from $v_i${\rm ,} or the foliation is
orthogonal to those edges and at least two leaves emanate
from $v_i$.
\endproclaim

{\it Proof}.  As an example, we discuss (iii). Here, an
appropriate branch of the map $z\mapsto z^{2/3}=w$ takes a
neighborhood of $v_i$ in the orthodisk to a\break neighborhood of
the origin in the upper half-plane. The condition on\break the
foliation on the proposition allow us to conclude \cite{Str}
that the holomorphic quadratic differential $\Phi$ has the
expansion in the upper half-plane as\break $(c_0+c_1z+\text{h.o.t.\/})dz^2$,
from which the statement follows.
\enddemo

\demo{{R}emark {\rm 5.3.3}} In light of the previous proposition, we
see that it is quite common to have an edge which is not
admissable for a holomorphic quadratic differential: for
instance we need only consider the differential of extremal
length for a cycle which has feet on an isolated edge with a
vertex where the orthodisk has angle at least $3\pi/2$. Thus,
it is a significant problem to find extremal lengths so that
Lemma~5.3.1 may be invoked, and this restriction
informed our choice of
height function after Definition~4.5.2. In particular,
we were careful not to have the special cycles referred to in 
Lemma~5.1.2 have feet adjacent to edges where we would later push.

Moreover, in the cases where an edge $E$ is not admissable
for a given holomorphic quadratic differential, we cannot ---
in general --- expect to have control over
$\ree\int_E\Phi\dot\nu$. This follows from an analysis
similar to that in the proof of the previous proposition or
in the final paragraph of the proof of Lemma~5.3.1: by
choosing a branch of the map from a neighborhood of a vertex
to the upper half-plane, it is easy to see that there is no
guarantee that $\lim_{b,\d\to0}\ree\int\Phi\dot\nu_{b,\d}$
exists. However, a choice of consistent branches of those
maps from a neighborhood $N(v)$ of a vertex $v$ and from
the reflection $N(v^*)$ of $N(v)$ across the central
diagonal show that $\ree\[\int_{N(v)}\Phi\dot\nu_{b,\d}
+\int_{N(v^*)}\Phi\dot\nu^*_{b,\d}\]=O(1)$
as $b$, $\d\to0$; i.e., the contributions from the terms in
the expansions of $\Phi$ of order $z^{-\a}$ from $\a\ge1$
cancel in pairs. Thus, while we then find that
$\lim_{b,\d\to0}\ree\int\Phi\dot\nu_{b,\d}$ exists --- which
we expect from $\dot\nu_{b,\d}$ always representing a unique
tangent vector to $\Delta$ --- we also find that the
contributions of the terms $z^{-\a}$($\a\ge1$) to the sign of
$\Phi$ on nearby edges must be ignored in evaluating
$\lim_{b,\d\to0}\ree\int\Phi\dot\nu_{b,\d}$.
\enddemo

{\it Conclusion of the proof of Theorem}~5.2.1. 
Conjugacy of the domains $\ogup$ and and $\ogdn$ requires that
if we push an edge $E \subset \partial\overline\ogup$ into
the domain $\ogup$, we will change the Euclidean geometry of
that domain in ways that will force us to push the
corresponding edge $E^* \subset \partial\overline\ogdn$  out
of the domain $\ogdn$. To see this, observe that for finite
length edges, the fact that the angles add to multiples of 
$2\pi$ (see the table in \S3.4) means that there is a pair 
of homeomorphisms of $\ogup$ and $\ogdn$, respectively, which
preserves the verticality or horizontality of the bounding
edges, but allows the two-image domains $\ogup^*$ and
$\ogdn^*$, respectively, to fit together to form a locally
Euclidean space.  (We can also require the homeomorphisms to
preserve the lengths of finite boundary edges; the domains
will continue to fit together as above.) We first discuss the
situation of pushes of finite boundary edges. Conjugacy
requires that $\ogup$ and $\ogdn$ have the same Euclidean
lengths, and for finite edges, this length is but a single
measurement in the Euclidean complex $\ogup^* \cup \ogdn^*$.
But in this complex $\ogup^* \cup \ogdn^*$, a push out of a
finite edge of $\ogup^*$ is simultaneously a pushing in of the
corresponding edge of $\ogdn^*$, and vice versa. Thus, 
because our homeomorphisms could be taken to preserve a
neighborhood of the boundary of the  domains $\ogup$ and and
$\ogdn$, we see that  a push out of a finite edge of 
$\ogup$ is simultaneously a pushing in of the corresponding
edge of $\ogdn$, and vice versa. 

The corresponding statement
for infinite edges is argued  analogously, once we remark
that in the above argument about finite edges, (i) we never
made use of the lengths of the  (finite) edges being pushed,
so this quantity is irrelevant to the discussion, and (ii) we
also never made formal use of the cone angles being at finite
points in the plane -- the arguments applied equally well to
boundaries of $\ogup$ and $\ogdn$ fitting together at
$\infty$.

We apply this general reasoning to the case of $\SC = {\rm DH}_{1,1}$.
Our only pushing is along the locus $\SY^*=\{\SH({\rm blue})=0\}$,
and we seek to push $H_1H_2$ and $H_2H_3$.  Here the measured
foliations for the holomorphic quadratic differentials 
$\Phi_{\rm mauve}=\frac12d\text{Ext}_{\Omega.}({\rm mauve})$ are parallel
to the sides being pushed and to the neighboring sides.  
Proposition~5.3.2(i),(iii) then assert that the quadratic
differentials are admissible for those edges; Lemma~5.3.1
then implies that since the Beltrami differentials
$\dot\nu_{\Omega.}$ have opposite signs on
corresponding edges, so do
$\Phi_{{\rm mauve}}\dot\nu_{\Omega.}$; thus
$$
\sgn d\text{Ext}_{\ogup}({\rm mauve})[\dot\nu_{\ogup}] = - 
\sgn d\text{Ext}_{\ogdn}({\rm mauve})[\dot\nu_{\ogdn}].
$$

This concludes the proof of Lemma~5.1.2 in this case, and
hence also the proof of Theorem~5.2.1.

\demo{Conclusion of the proof of Lemma~{\rm 5.1.2}} 
We begin with statement (i).
The general case of this statement is precisely
analogous to the case of Theorem~5.2.1; indeed, we can regard the 
locus $\SY$ from Lemma~5.1.2 as being a version of the moduli space
$\Del$ studied in Theorem~5.2.1. 
The coordinates $t_{i_1}$ and $t_{i_2}$ from the lemma
are given by periods $H_1C_1-P_1P_2$ exactly as in the treatment
of $\Delta_{1,1}$. In particular, the asymptotic estimates
calculated in the proof of Claim~5.2.2 continue to hold exactly
as before, as the setting for those estimates has changed in
only insubstantial ways.  (To see this, note that the
only changes are the we do not necessarily have that $P_2$ tends to
$H_2$ on $\ogup$ (which does not affect the estimates since
we continue to have $|H_1C_1-P_1P_2|$ bounded away from zero), and
the angle at $P_2$ on $\ogup$ has changed 
(also immaterially as $P_2$ is at
some bounded distance from where the degeneracy is occurring).)
 
Moreover, because the rest of the terms
of $\height_{\SC}$ vanish along the locus $\SY$ to second order in the
deformation variable, we see that any deformation of the orthodisk
will not alter the contribution of these terms to $\height_{\SC}$.
Thus the only effect of an infinitesimal deformation of an orthodisk
system on $\SY$ to the height function $\height_{\SC}$ is to the
terms $\height(t_{i_1})$ and $\height(t_{i_2})$, which we control as in
Theorem~5.2.1. 

For statement(ii) of the lemma, we deform as we did along the locus 
$\SY^* \subset \SY$, here pushing on a finite edge which is parallel
to the foliation for $\Phi_{t_i}$ (in the obvious notation) if and
only its neighbors are. (It is worth remarking that in the
case of ${\rm DH}_{1,2}$, we will not be able to find an edge which is
finite on both $\ogup$ and $\ogdn$; in this case we will push on
$H_1C_1$, which, relative to the cycle $H_2H_3-P_1P_2$ is admissible
by Proposition~5.3.2(i),(ii).  Thus the proof above for the case 
$m=n$ also extends to this case.)

This concludes the proof of the lemma.
\enddemo

\vglue-12pt
\section{Regeneration}

In the previous section we showed how we might reduce the height
function $\height_{\SC}$ at a critical
point of a smooth locus $\SY$, where the locus  $\SY$ was 
defined as the null locus
of all but two of the heights, say $\height(t_{i_1})$ and
$\height(t_{i_2})$. In this section, we prove Lemma~5.1.3, which 
guarantees the existence of such a locus~$\SY$.

We will continue with the rather concrete exposition we used in
\S5; i.e., we will prove the existence of the genus five Costa tower,
${\rm DH}_{2,2}$ by using Theorem~5.2.1 (the existence of the genus
three Costa tower, ${\rm DH}_{1,1}$) to imply
the existence of a locus $\SY \subset \Delta = \Delta_{2,2}$ ---  
Lemma~5.1.2 and the properness Theorem~4.5.1 then prove the
existence. The arguments we give will immediately
generalize to prove the existence of a genus $2n+1$ Costa tower
${\rm DH}_{n,n}$, given the existence of a genus $g=2n-1$ Costa tower 
${\rm DH}_{n-1,n-1}$.

Indeed, our proof of Theorem~A is by induction:  
\vglue6pt
{\it Inductive Assumption} A. There exists a genus $2n-1$ Costa tower,\break
${\rm DH}_{n-1,n-1}$.
\vglue6pt
Thus, all of our surfaces are produced from only slightly less
complicated surfaces; this is the general principle of `handle
addition' referred to in the title.

For concreteness, our present goal is the proof of the statement: 

\nonumproclaim{Theorem~6.1} There is a reflexive orthodisk system
for the configuration ${\rm DH}_{2,2}$.
\endproclaim

{\it Start of proof}. Let us use the given height 
$\height_{2,2}$ for ${\rm DH}_{2,2}$ and consider how the
height $\height_{1,1}$ for ${\rm DH}_{1,1}$ relates to it, 
near a solution $X_1$ for the ${\rm DH}_{1,1}$
problem. (As we observed, the situation extends to general ${\rm DH}_{n,n}$
and ${\rm DH}_{n-1,n-1}$ with simple added notation.)

Our notation is given in \S4.4 and is recorded in 
the diagrams below: for instance,
the curve system $\d$ connects the edges $P_1P_2$ and $H_2H_3$.

We are very interested in how an orthodisk system might
degenerate. One such degeneration is shown in the next
figure, where the  branch point $P_1$ has collided with the
vertex $H_1$ on the orthodisk $\ogdn$. This has the effect,
on the orthodisk $\ogup$ of causing $H_1C_1$ to collapse and
forcing the line $P_1P_2$ (respectively, $P_4P_5$) to lie
directly across from (resp., over) the point $H_1 = C_1$
(resp., $C_2=H_5$). The point is that the degenerating family
of (pairs of) Riemann surfaces in $\Del_{\SC}$ limits on (a
pair of) surfaces with nodes. (We recall that a surface with
nodes is a complex space where every point has a neighborhood
complex isomorphic to either the disk $\{|z|<1\}$ or a pair
of disks $\{(z,w) \in \bbc^2| zw=0\}$.) In the case of the
surface corresponding to $\ogdn$, the components of the noded
surface (i.e. the regular components of the noded surface in
the complement of the nodes) are difficult to observe, as the
flat structures on the thrice-punctured sphere components are 
simply single points; on the other hand, in the case of the
surface corresponding to $\ogup$, the nodes are quite visible, as 
coming from a pair of sheets on the orthodisk, with one
infinite point and one  finite point (drawn in paler lines
on the figure below on the right). 
\vfill
\centerline{\BoxedEPSF{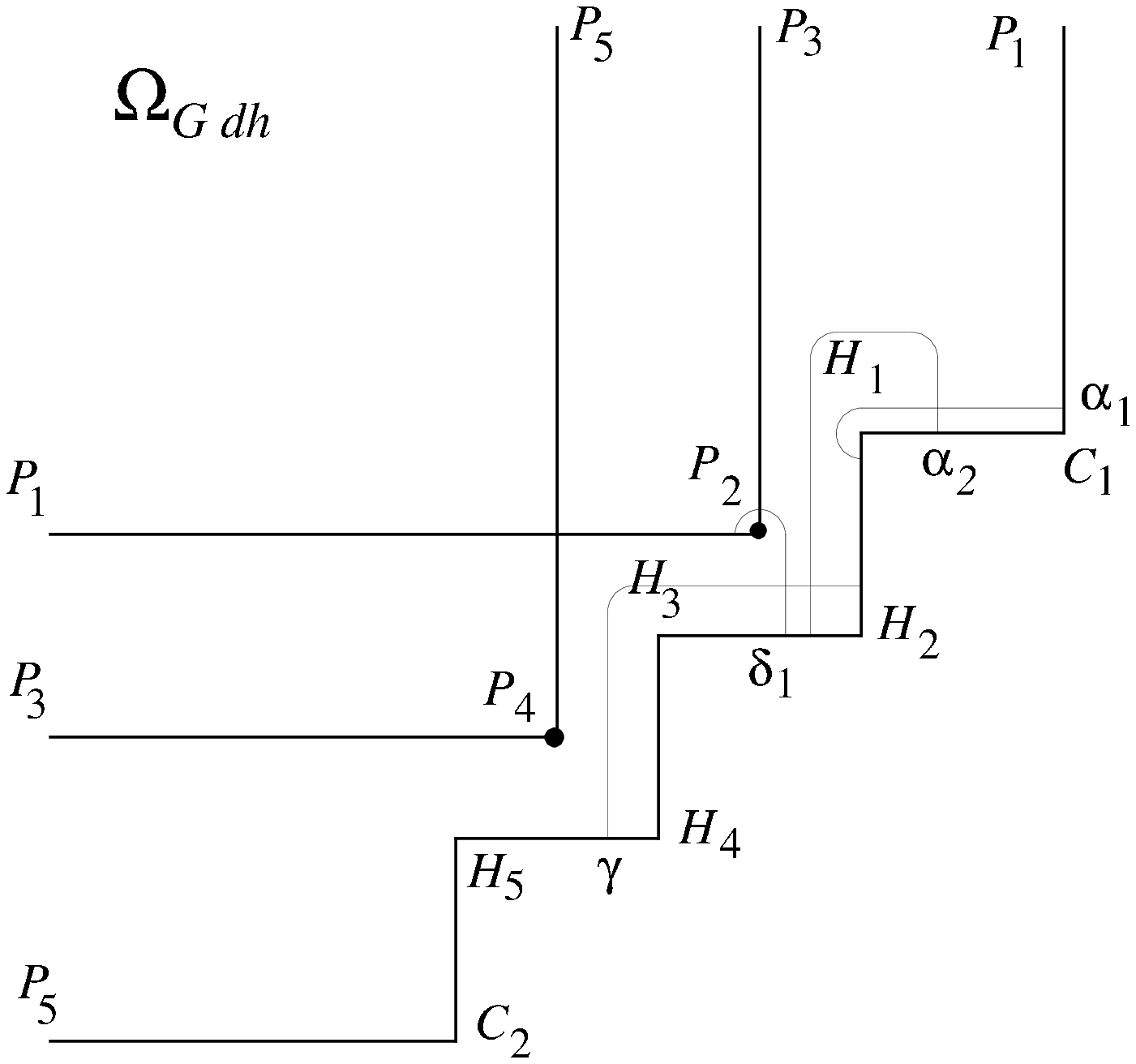 scaled 450}\qquad\BoxedEPSF{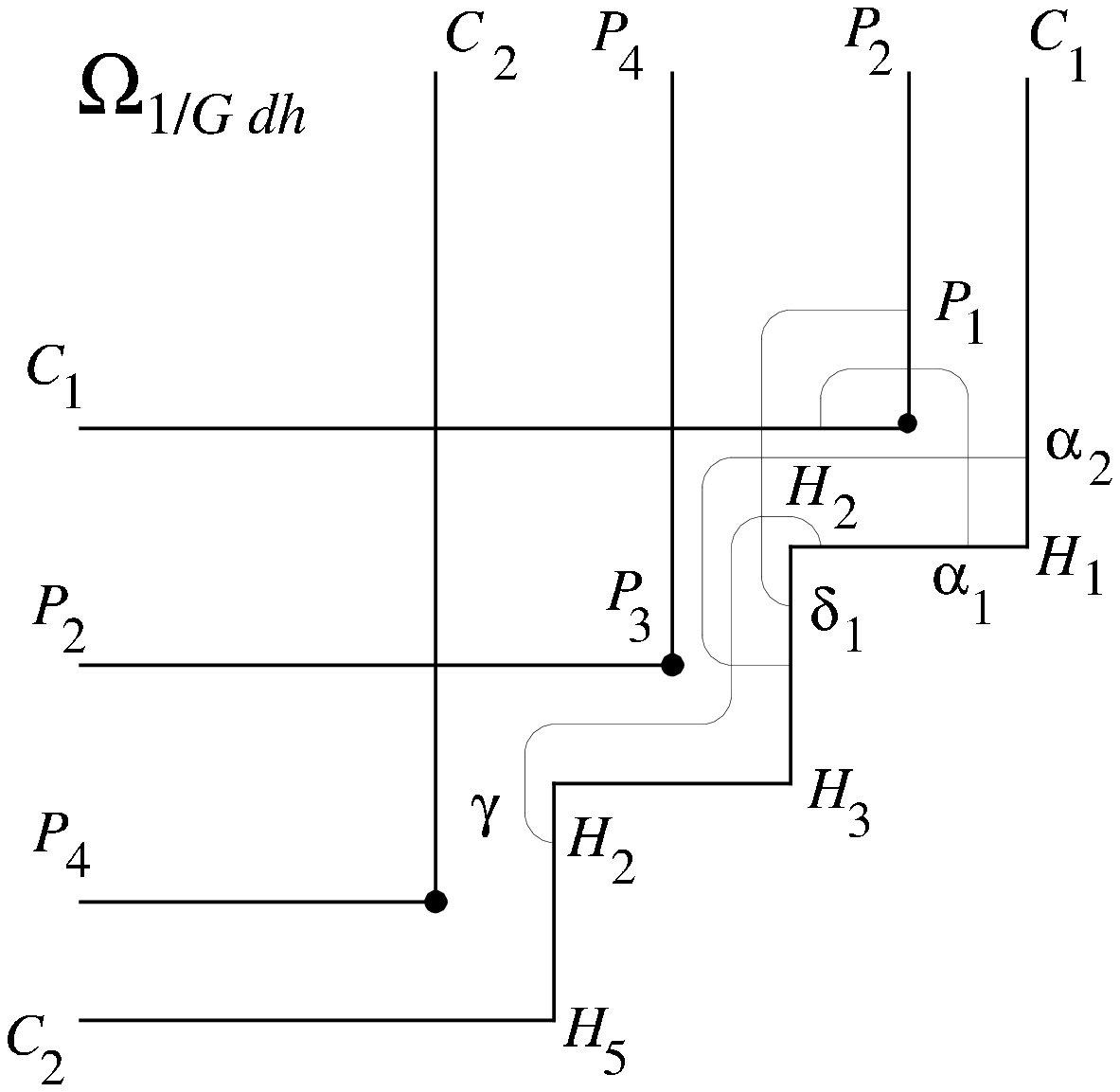 scaled 450}}

\vglue12pt
\centerline{Curve system used for regeneration}
\vfill

  {\it Remark}. In all of this, it is important to keep in
mind that our flat structures represent a Riemann surface
with a meromorphic  one-form. Thus while there are no
nontrivial holomorphic one-forms on a  sphere (so that the
flat structure is  given by the trivial form $0\cdot dz$, and
is hence a point), there are nontrivial {\it mero}morphic
forms on a thrice-punctured sphere with poles only at the
punctures. The flat structures for those that occur in our
present situation  are represented in  our method by a bigon,
with one finite and one infinite vertex, connected by one
horizontal and one vertical side.
\eject

An important issue in this section is that some of our curves
cross the pinching locus on the surface, i.e. the curve on
the  surface which is being collapsed to form the node. In
particular, in the diagram, the curves $\mu$ and $\nu$ are
such curves, and so their depiction in the degenerated figure is,
well, degenerate; the curves either surround a point, or
connect a point and an edge, or connect a pair of edges on
opposite sides of a node.
\vglue6pt
\centerline{\BoxedEPSF{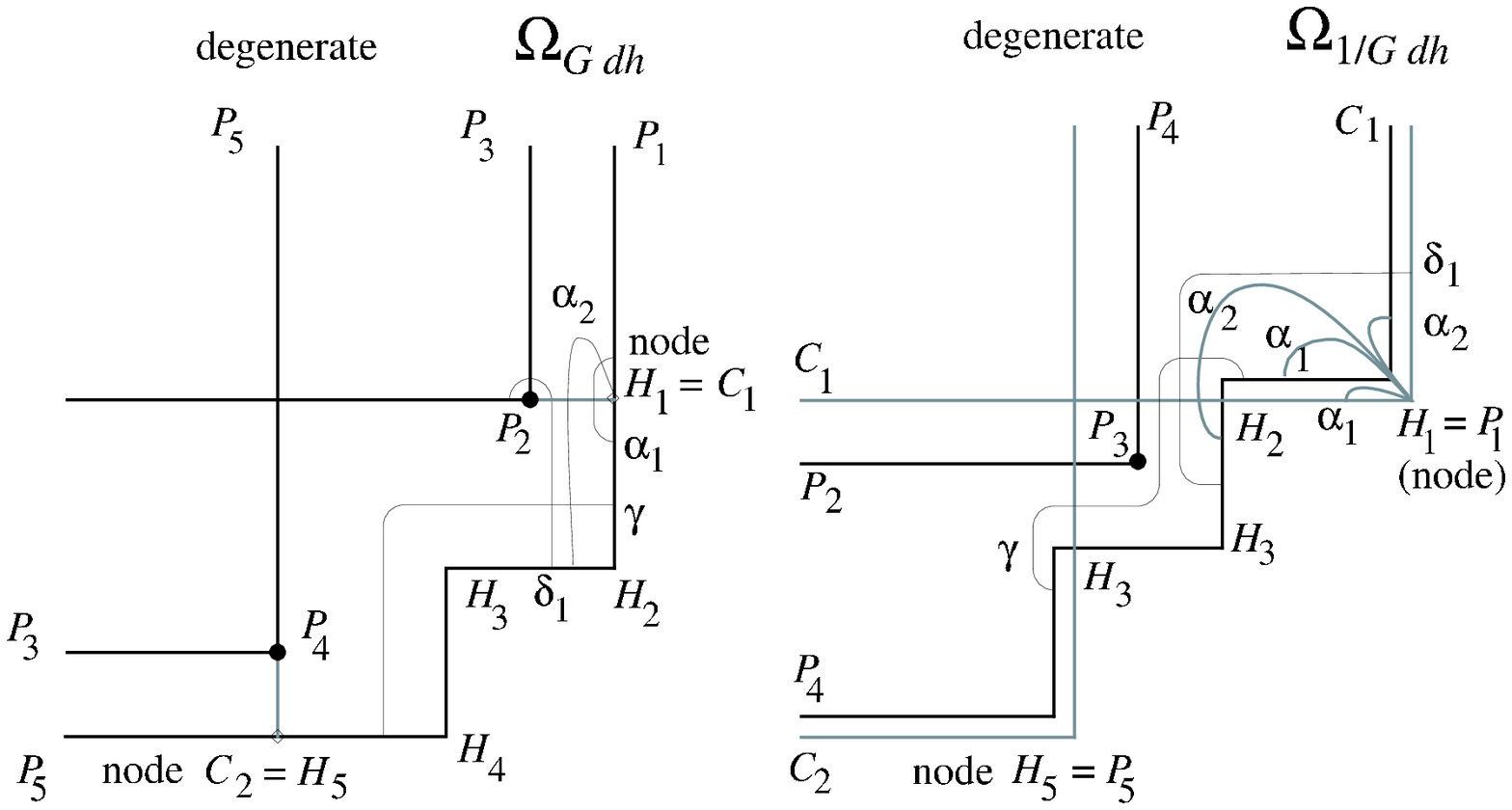 scaled 600}}
\vglue12pt
\centerline{Degenerated orthodisks for the Costa towers }
\vglue12pt

When we remove the thrice-punctured sphere components, we
observe  (compare with the figures in \S5.2 and \S3.6) that
we are left with an orthodisk system for the lower genus
configuration $\SC_{{\rm DH}_{1,1}}$, up to possibly relabeling
$\ogup$ as $\ogdn$ and vice versa.

Our basic approach is to work backward from this
understanding of degeneration --- we aim to ``regenerate'' the
locus $\SY$ in $\Del_{{\rm DH}_{2,2}}$ from the solution $X_1 \in
\Del_{{\rm DH}_{1,1}}$.

We focus on the curves $\d$ and $\g$, ignoring the
degenerate curves $\mu$ and $\nu$. Let $\Delta_2$ be the
geometric coordinate simplex for ${\rm DH}_{2,2}$.

(In the general case for $\Delta_{{\rm DH}_{n,n}}$, there are $2n-2$
nondegenerate curves (say, for notational convenience $\{
b_1,\ldots, b_{2n-2}\}$),  and two degenerate curves, $\mu$ and
$\nu$. The case of ${\rm DH}_{m, n \ne m}$ requires a separate but
parallel treatment.)

We restate Lemma~5.1.3 in terms of the present (simpler)
notation.

\nonumproclaim{Lemma 6.2 {\rm (Regeneration)}} There is a smooth
two\/{\rm -}\/dimensional analytic closed locus
$\SY\subset\ov{\Delta_n}$ so that
$\ext_\ogup(b_i)=\ext_\ogdn(b_i)$ on $\SY$\/{\rm ,}  and $\SY$ is
proper in $\Delta_n$.
\endproclaim

\demo{Proof}  We again continue with the notation for $n=2$,
as the general situation follows with just added notation.

As putatively defined in the statement of the lemma, the locus
$\SY$ would be clearly closed, and would have nonempty
intersection with $\ov{\Delta_2}$ as $\ov{\Delta_2}$ contains
the solution $X_1$ to the genus 3 problem.\vglue2pt

We parametrize $\ov{\Delta_2}$ near $X_1$ as
$\Delta_1\x[(0,\e)\x(0,\e)]$ and consider the map
$$
\Phi:(X,(t_1,t_2)): \Delta_1\x[(0,\e)\x(0,\e)]\lra\bbr^2
$$
given by
$$
(X,(t_1,t_2))\mapsto(\ext_\ogup(\d)-\ext_\ogdn(\d),
\ext_\ogup(\g)-\ext_\ogdn(\g)).
$$
Here, coordinates $t_1$ and $t_2$ refer to a specific choice
of normalized geometric coordinate; i.e.,
$t_1=\im(P_5C_2\lra H_5H_4)$ and
$t_2=\ree (C_2H_5 \lra P_4P_5)$,
where the periods $P_5C_2\lra H_5H_4$ and $C_2H_5 \lra P_4P_5$
are measured on the domain $\ogup$ (see the figure below). 
In terms of these coordinates, we note that whenever either
$t_1 =0$ or $t_2=0$, we are in a boundary stratum of
$\ov{\Delta_2}$. The locus
$\{ t_1>0, t_2>0\} \subset \ov{\Delta_2}$ is a neighborhood
in $\text{Int}(\Delta_2)$ with $X_1$ in its closure.
\vglue12pt
\centerline{\BoxedEPSF{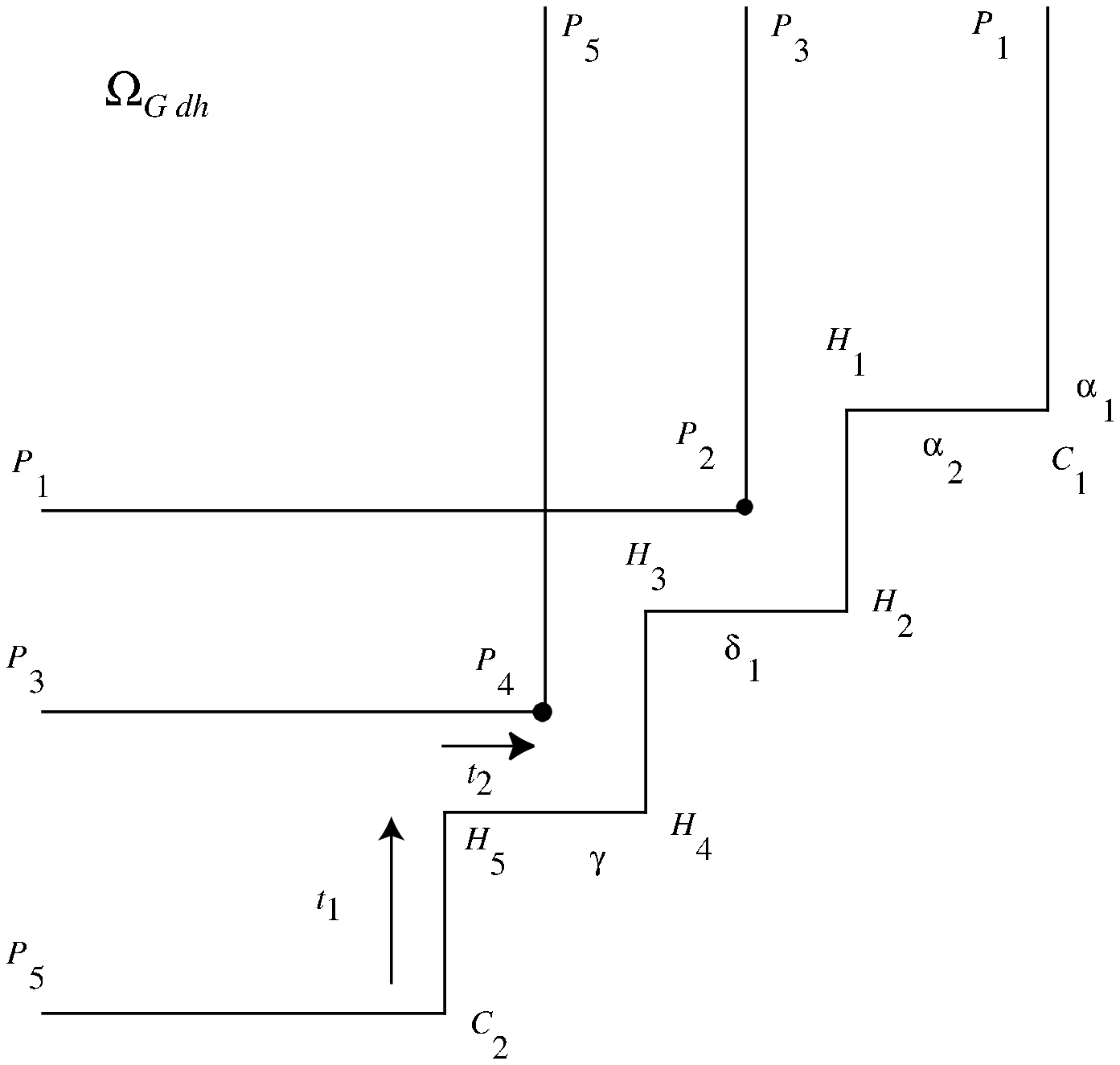 scaled 560}}
\vglue6pt
\centerline{Geometric coordinates $t_1$ and $t_2$ for $\Delta_2$}
\vglue6pt

Note that $\Phi(X_1,(0,0))=0$ as $X_1$ is reflexive.

Now, to find the locus $\SY$, we apply the implicit function
theorem, which  says that if
\vglue4pt
(i) the map $\Phi$ is differentiable, and
\vglue4pt
(ii) the differential $d\Phi\bigm|_{T_{X_1}\Delta_1}$ 
is an isomorphism onto $\bbr^2$,

\vglue4pt\noindent then there exists a
differentiable family $\SY\subset\ov{\Delta_2}$ for which
$\Phi\bigm|_\SY\equiv0$.
\eject

In order to show the differentiability condition (i), 
we require a clarification as to the meaning
of the curve system
$\g$, for points in $\Delta_1\subset\p\ov{\Delta_2}$. Here on
$X_1$ we consider the distinguished points $x$ and $x'$ 
on the $\ogdn$
orthodisk, as shown below. 
\figin{Zreg3}{525}
 
\centerline{Definitions of the points $x$ and $x'$}
 \vglue12pt

By definition, these points $x$ and $x'$ on $\ogdn$ of the
solution $X_1$ are  determined from continuing the levels
$C_1P_1$ and $P_3C_2$ until they meet the axis at $H_1C_1$ and
$C_2H_3$, respectively.  Then $\g$ is defined on
$\Delta_1\subset\p\ov{\Delta_2}$ as the curve family
connecting $H_1x$ and $H_3x'$. We extend the notion of the
points $x$ and $x'$  to the interior of $\Delta_2$ by
declaring there that $x=H_1$ and $x'=H_5$, respectively. The
definition of $\g$ then extends continuously to
$\ov{\Delta_2}$.

We now prove the differentiability (condition (i)) of $\Phi$.
As the locus\break of
$\ov{\Delta_2}\in\wh{\ov\SM_5} \x\wh{\ov\SM_5}$ is
differentiable (here $\wh{\ov\SM_5}$ refers to a smooth cover
of the relevant neighborhood of $X_1\subset\ov\SM_5$, where
$\ov\SM_5$ is the Deligne-Mostow compactification of the
moduli space of curves of genus five) the theorem of
Gardiner-Masur [GM] implies that $\Phi$ is differentiable, as
we have been very careful to choose curves $\{\d, \g\}$
which are nondegenerate in a neighborhood of $\ov{\Delta_2}$
near the genus three solution $X_1$, with both staying in a
single regular component of the noded surface.

We are left to treat (ii), the invertibility of the 
differential $d\Phi\bigm|_{T_{X_1}\Delta_1}$. To show that
$d\Phi\bigm|_{T_{X_1}\Delta_1}$ is an isomorphism, we simply
prove that it has no kernel. To see this, choose a tangent
direction in $T_{X_1}\Delta_1$; as we may regard
$\Delta_1\subset T_{0,8}\x T_{0,8}$ as a subspace of the
product of two \tec spaces of $8$-times punctured spheres, we
may regard a tangent direction as a pair
$(\nu_{Gdh},\nu_{G^{-1}dh})$ of Beltrami differentials, each
representing a tangent direction to the points $[\ogup]\in
T_{0,8}$ and $[\ogdn]\in T_{0,8}$, respectively. Yet at
$X_1$, the points $[\ogup]$ and $[\ogdn]$ represent the
identical point in $T_{0,8}$, so we begin by computing how
the Beltrami differentials $\nu_{Gdh}$ and $\nu_{G^{-1}dh}$
relate to one another. To this end, consider how an
infinitesimal push in the sense of \S5.3 on an edge $E$
defines Beltrami differentials $\nu_{Gdh}$ and
$\nu_{G^{-1}dh}$. Of course the conjugacy of $\ogup$ and
$\ogdn$ provides, via the formulas of \S5.3 the basic
defining relation that if $\nu_{Gdh}$ has local expansion
$\nu_{Gdh}(z)=\f1{2b}\f{d\bar z}{dz}$ near an interior point
of an edge $E$, then
$\nu_{G^{-1}dh}(\z)=\f1{2b}\f{d\bar\z}{d\z}$. However, since
$X_1$ is reflexive, we may also assume, in this particular
case, the existence of a conformal map $\z:\ogup\lra\ogdn$
which preserves the vertices. Such a map takes vertical sides
to horizontal sides by construction and this has the local
expansion $\z=\pm i|c|z+0(|z|^2)$ near an interior point of an
edge. We therefore compute the pullback of $\nu_{G^{-1}dh}$
to $\ogup$ as
$$
\split
\nu_{G^{-1}dh}(\z)\f{d\bar\z}{d\z} &= \nu_{G^{-1}dh}(\z)
\bar\z'/\z'\f{d\bar z}{dz}\\
&= (-\nu_{G^{-1}dh}(\z) + \text{h.o.t.\/})d\bar z/dz\\
&= -\f1{2b}\f{d\bar z}{dz}
\endsplit
$$
along the edge $E$.  As we found in Remark 5.3.3, a similar computation (using symmetry to cancel apparent singularities) holds
near the vertices.  Thus
$\z^*\nu_{G^{-1}dh}=-\nu_{Gdh}$, and since any deformation of $X_1$ is given by a linear
combination of such infinitesimal pushes, we conclude that
$[\nu_{Gdh}]=-[\nu_{G^{-1}dh}]$ as elements of the tangent
space $T_\ogup T_{0,8}=T_\ogdn T_{0,8}$. Thus, any conjugacy
preserving deformation of $X_1$ destroys the conformal
equivalence of $\ogup$ and $\ogdn$ to the order of the
deformation, a statement which implies that
$\ker(d\Phi)=\{0\}$, as required.

To finish the proof of the lemma, we need to show that
$\SY\bigm|_{\Delta_2}$ is an analytic submanifold of $T_5^{\phantom{|}}\x
T_5$, where $T_5$ is the \tec space of  genus five curves.
Now, there are two sets of analytic coordinates on
$\Delta$; one comes from using the extremal lengths
on $\ogup$ of the height cycles, and one comes from
using the extremal lengths on $\ogdn$ of the 
height cycles. As the map between the two coordinate
systems is analytic with an analytic inverse
(by Ohtsuka's formulae), we see that the locus
$\SY$, being defined as a level set of linear 
functions in these coordinates, is analytic as well.

This concludes the proof of Lemma~6.2 for the case $n=2$ and
hence also the proof of Theorem~6.1. We have already noted
that the argument is completely general, despite our having
presented it in the concrete case of ${\rm DH}_{2,2}$; thus, by
adding more notation, we have proved Lemma~6.2 in full
generality.
 \enddemo

{\it Conclusion of the proof of Lemma}~5.1.3. All that is
left in the proof of Lemma~5.1.3 is the remaining case (Case
(ii)) of $m<n$. In this case, the degeneration of the genus
$g+1$ case to the genus $g$ case is given by a much simpler
operation than that of the degeneration described earlier in
this section: here the degeneration consists of collapsing
the two central finite edges of the outer  sheet of both
$\ogup$ and $\ogdn$ (see the figures below). Thus, the
corresponding regeneration consists of removing a
neighborhood of the central finite edges of the outer sheet
of both $\ogup$ and $\ogdn$, and reconnecting the boundary
with a pair of short edges (see the figures below and compare
to the identical operation in [WW]). 
\figin{Zreg4}{750}

\centerline{Regenerating orthodisks for ${\rm DH}_{m,n}$, case $m<n$ }
\vglue12pt

The effect of this surgery is to replace the central vertex
with a triple of vertices which lie close to each other, as
well as changing the cone angle of the central vertex. The
proof of Lemma~5.1.3 is then virtually identical to that of
the case of $m=n$: the only change is that instead of there
being two curves $\mu$ and $\nu$ that degenerate under the
degeneration, there is, in the present case, just a single 
degenerate curve $\rho$. Then with that change, the implicit
function theorem argument given above goes through unchanged,
up to substituting $\rho$ for $\mu$ and $\nu$, and
obtaining a locus $\SY$ which is one-dimensional, rather than
the two-dimensional locus $\SY$ we obtained in the previous
case $m=n$. 
To satisfy the requirement that the cycle corresponding to 
the $t_i$ not be a finite edge bully, we choose our 
new cycle to foot on the central finite edge and on 
an infinite edge. Once we do this, the furthest 
finite edge $E=H_1C_1$ (and, unless $(m,n) =(1,2)$, its
immediate neighbors) will have neighbors which are
parallel to the foliation for $\Phi_{t_i}$ if and only if
that foliation is parallel to $E$. 
This concludes the proof of Lemma~5.1.3, and
hence of the Main Theorems A and B(i) (see the discussion
after Lemmas~5.1.2 and 5.1.3),  up to the proof that
${\rm DH}_{m,n}$ is embedded outside a compact set, which we
treat in \S9. \eject

\section{Nonexistence of ${\rm DH}_{m,n}$ for $m>n$}

In the previous sections, we have seen that it is possible to
add inductively a  handle or a plane together with a handle
to the minimal surfaces of type ${\rm DH}_{m,n}$, starting at
${\rm DH}_{1,1}$. This raises the question whether one can just add
planar ends. If we could  do this with a single
surface ${\rm DH}_{n,n}$, this would (assuming its embeddedness)
provide a counterexample to the Hoffman-Meeks conjecture. In
this section, we will show that this is not the case, at
least not in our setting of symmetry assumptions and chosen
module sequences.

More precisely, we show that there are no surfaces of type
${\rm DH}_{m,n}$ for $m>n$. This proves the second statement of
Theorem~B. Recall that the class of surfaces excluded by this
statement consists only of those construed from module
sequences of type ${\rm DH}_{m,n}$ which were characterized among
other possible module sequences as follows: The intersection
of the surface with the vertical symmetry planes has one
connected component joining all  points with vertical normal
on the $z$-axes; all other components join just  two
consecutive ends.

Before proving the general statement, we give the basic
example:

\vglue8pt {\it Example} 7.1.  {\it There is no} ${\rm DH}_{1,0}$.
If there were an example of such a surface, then we 
could compute, as in \S3, that the corresponding domains
$\ogup$ and $\ogdn$ would have the following shapes:

\centerline{\BoxedEPSF{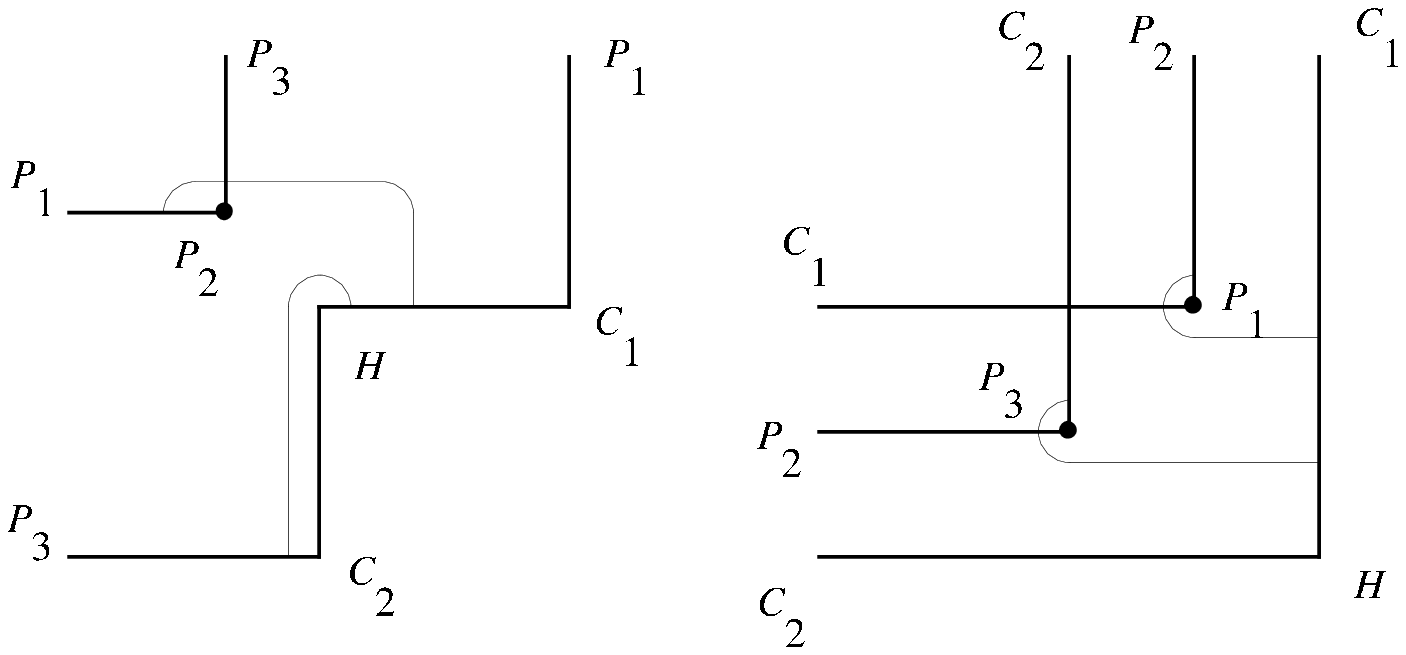 scaled 650}}
 
\centerline{Nonexistence of ${\rm DH}_{1,0}$}
\vglue6pt
\pagegoal=50pc 
For these domains to be conjugate, we require that the
corresponding  periods from $\ogup$ and $\ogdn$ sum to lie on
one line;  in all of our examples so far, this line has been
the 45-degree line $\{y=x\}$. In the present case, we first
determine this conjugation line by looking at the periods of
the cycle $P_3 C_2 \to H C_1$. In $\ogup$, the period points
upward while in $\ogdn$, it points to the right. So the
conjugation line is necessarily the $\{y=x\}$ diagonal. Now
look at the cycle $HC_1\to P_1P_2$. By symmetry, its period
must point upward in $\ogup$. By conjugacy, its period in
$\ogdn$ must point to the right. This is impossible, since
the branch point $P_1$ cannot lie outside the outer sheet.

\nonumproclaim{Theorem 7.2} There is no minimal surface of type
${\rm DH}_{m,n}$ with $m>n$.
\endproclaim
\pagegoal=48pc

{\it Proof}. Assume that we have a pair of reflexive
orthodisks with $m>n$. Depending on the parities of $m$ and
$n$, there are four  cases to distinguish, but the argument is
the same in all cases. For the sake of concreteness, let us
look at the case $m$ even and $n$ odd. As in the example
above, the conjugation line is necessarily the $y=x$ axis.
Now (in this parity case) the  orthodisks  have the form,
near the symmetry axis, as in the diagrams below: 
\figin{Zinddh}{600}
 
\centerline{General nonexistence by induction}
\vglue12pt

By symmetry, the cycle $A: H_{n+1}H_n \to P_mP_{m+1}$ has in
$\ogdn$ a period pointing upward, so in $\ogup$ a period
pointing to the right. To keep the branch point $P_m$ within
the outer sheet, this forces the period of the cycle
$B: H_n H_{n-1} \to P_{m-1}P_m$ to point upward in $\ogup$
and hence to the right in $\ogdn$.  This in turn forces the
period of the cycle $C: H_{n-1}H_{n-2}\to P_{m-2}P_{m-1}$ to
point upward in $\ogdn$. Proceeding inductively like this, we
will reach at some point the conclusion that the cycle $H_1C_1 \to
P_{m-n}P_{m-n+1}$ (which exists because  $m>n$) points
upward in $\ogup$ and hence to the right in $\ogdn$. However,
this clearly forces  the branch point $P_{m-n}$ to lie outside
the outer sheet. The other cases are completely analogous, 
starting always at the symmetry diagonal on the orthodisk with
a central inner sheet.

\vglue-8pt

\section{Extensions and generalizations}
\vglue-5pt

 8.1. {\it Higher dihedral symmetry}. It is a well-known phenomenon that a surface like our
${\rm DH}_{m,n}$ with $4$-fold dihedral symmetry has companions with
higher symmetry.  In the framework of our approach, this is
quite easy to achieve --- we only have to change the setup in
\S3 while the rest of the arguments in \S\S4--6 remain valid.
The orthodisks which reflect, in their orthogonal geometry,
the $4$-fold dihedral symmetry have to be replaced by {\it
skewdisks} whose interior or exterior angles are $2\pi/d$
where $d$ is the order of the dihedral symmetry group. 

Recall that our orthodisks represent a piece of the surface
obtained by cutting the surface by its two vertical symmetry
planes with the metric induced by a meromorphic $1$-form. In
the same way, a surface with $2d$-fold dihedral symmetry and
an additional rotational symmetry around (say) the $x$-axis
is cut into $2d$ pieces, which we represent by skewdisks.

Vice versa, given such a skewdisk, one can take its double to
obtain a skewsphere and then take a $d$-fold cyclic branched
cover branched over the vertices. This covering should be
chosen  so that every second edge of the original orthodisk
becomes a branch cut, and so that each of the $d$ copies of the
skewspheres is glued to two other spheres along these cuts.
This is needed to guarantee that the lifted cone metrics
indeed define meromorphic $1$-forms, i.e. have no linear
holonomy.

Here is an example of two such skewdisks for $d=4$ which
represent \wei candidate data for the
${\rm DH}_{1,1}$ surface with $8$-fold dihedral symmetry:
\vglue6pt
\centerline{\BoxedEPSF{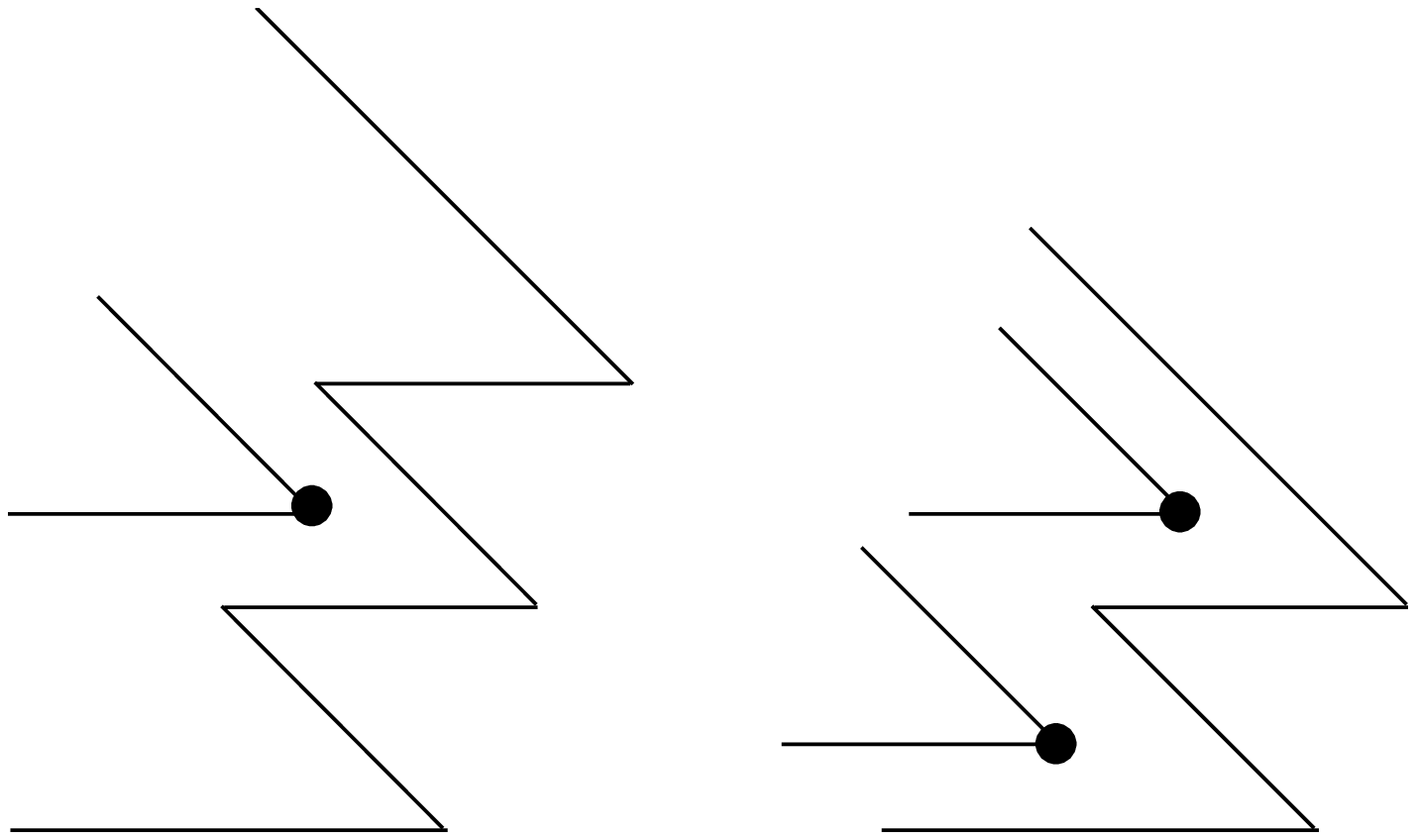 scaled 280}}
\vglue5pt
 
\centerline{Skewdisk for surfaces with $8$-fold dihedral symmetry}
\vglue6pt

All of the definitions and claims from \S3 apply to this
situation, with but a few modifications. The vertex data
become odd integers $a_i\equiv \pm 1 \pmod {2d}$ where $a_i$
represents a vertex with angle $\pi(a_i+2)/d$. Under the
construction of the Riemann surface, this becomes a cone angle
$2\pi(a_i+2)$ which corresponds to a zero of order $a_i+1$ of the
corresponding meromorphic $1$-form. The sufficient condition
in Lemma 3.3.4 becomes  $a_i+b_i\equiv 0 \pmod {d}$, with
similar straightforward changes in the other affected results
in \S3. Also note that by definition of the skewdisks and the
associated Riemann surfaces, the divisors of $G$ and $dh$
balance so that we obtain the \wei representation of a
complete regular minimal surface.
\pagegoal=50pc

In \S4, the construction of the height function is
unchanged, and for the proof of its properness
only \S4.7 needs to be changed slightly: in Corollary~4.7.4 we get a different phase factor due to the fact that
the angles  under consideration are no longer
right angles. But all we need in this section   is that
\vglue4pt
\item{(1)} the angles are $\equiv \pm \pi/d \pmod {2\pi}$,
\vglue4pt\item{(2)}
corresponding angles of $\ogup$ and $\ogdn$ always add to $0
\pmod {2\pi}$.
\eject

For Sections 5 and 6 we note that the skewdisks are just
affinely distorted orthodisks; this gives complete and
explicit control over the behavior of all extremal lengths,
so that the tedious computations of \S5 do not need to be
repeated.
\pagegoal=48pc

Thus we arrive at the next result:

\nonumproclaim{Theorem 8.1.1}  For every triple of integers
$(d,m,n)$ with $d\ge2${\rm ,}\break $1\le m\le n$ there exists a complete
minimal surface ${\rm DH}_{d,m,n}\subset\bbe^3$ of genus\break
$(m+n+1)(d-1)$ which is embedded outside a compact set with
the following properties\/{\rm :} it has $2n+1$ vertical normals at
flat points of order $d${\rm ,} $2m+1$ planar ends{\rm ,} and two catenoid
ends. The symmetry group is generated by a dihedral group of
order $2d$ and a rotational symmetry around a horizontal
axis. 
\endproclaim

 8.2. {\it Deformations with more catenoidal ends}.
The surfaces of type ${\rm DH}_{m,n}$ described in Theorem B had two
catenoid ends and $2m+1$ planar ends. Planar ends can be viewed as
catenoid ends with a vanishing coefficient of growth (compare
\cite{Kap}), so it is natural to conjecture that the surfaces
described in Theorem~B are each elements of a family of
surfaces obtained by deforming some of the planar ends into
catenoid ends of varying growth rates, at least for growth 
rates near zero. In this section, we prove a theorem 
confirming that conjecture within the context of the surfaces  studied
in this paper, i.e., those surfaces with exactly eight
symmetries.

We begin with some notation and we define the {\it growth rate}
$\alpha_i$ of a catenoidal or planar end $E_i$ to be the residue of
the form $dh$ at that end. The surfaces  created in this 
section will have a symmetry about a central line; thus the 
growth rate $\a_i$ of the $i^{\rm th}$ end $E_i$ will be the 
negative of the growth rate $\a_{2m+4-i}$ of the symmetric
end $E_{2m+4-i}$.

For a minimal surface with catenoidal or
planar ends, let $p$ denote the number of planar ends and $k$ denote
the number of catenoidal ends.

\nonumproclaim{Theorem 8.2.1} Choose growth rates
$\a_1,\dots,\a_m\in\bbr$ near zero. Then there is a complete
minimal surface ${\rm DH}_{m,n}(\a_1,\dots,\a_m)$ in $\bbe^3$ with
$2n+1$ vertical normals{\rm ,} $k$ catenoidal ends and $p \ge 1 $
planar ends{\rm ,} where $k+p=2m+3$.  The ends have growth rates 
$1${\rm ,} $\a_1,\dots,\a_m, 0, -\a_m, \dots, -\a_1$, $-1$.
\endproclaim

The proof follows the outline of the previous work, but
involves a  new type of special point of the type $R$
described in the  table in \S4. Space considerations require
us to defer the complete proof to \cite{WW2}. There we also
discuss how to combine the constructions of Theorems~8.1.1
and 8.2.1.

\demo{{\rm 8.3.} Embeddedness aspects of ${\rm DH}_{m,n}$}
In this section we collect some information about how
embedded our minimal surfaces ${\rm DH}_{m,n}$ are. First, for a
complete minimal surface of finite total curvature to be
embedded, the ends must be planar or catenoidal, all parallel
and at different height. For different catenoidal ends which
open  in the up direction, the growth rate of the higher
end must be larger than the growth rate of the lower end so
that they do not intersect. A similar statement holds for
catenoidal ends opening in the down direction.

These conditions are necessary but not sufficient for
embeddedness.
\enddemo

 {\it Definition} 9.1. We call a complete minimal surface of finite
total curvature satisfying the above condition {\it eventually embedded}.
\vglue12pt

In the next result, we state that all of the surfaces 
${\rm DH}_{m,n}$ (with two catenoidal ends) are eventually embedded,
completing the proof of Theorem~B. In fact, we prove a
slightly more general statement:

\nonumproclaim{Theorem 9.2} Suppose there is a complete minimal
surface of finite total curvature given by $\ogup$ and $\ogdn$
orthodisks such that the formal \wei data $a_i$ and
$b_i$ satisfy
\vglue4pt
\ritem{(1)} $a_i, b_i \in \{-3,-1,1,3\}${\rm ,}
\vglue4pt\ritem{(2)} $2+a_i+b_i \le \vert a_i-b_i \vert${\rm ,}
\vglue4pt\ritem{(3)} $\a_i + b_i \equiv 0 \pmod 2$.
\vglue4pt\noindent 
Suppose  that there are only two catenoidal ends{\rm ,} one
necessarily pointing up and the other down. Finally  assume
that all planar ends occur in only one component of the
orthodisk boundaries from which the two catenoidal vertices
have been removed. Then the minimal surface is eventually
embedded.
\endproclaim

{\it Proof}.
The first condition ensures that we only have catenoidal ends
or planar ends; the second is equivalent to completeness.
Because we only have two catenoidal ends, we do not have to
consider growth rates. The only thing which remains to be
checked is that the ends are at different heights. Using the
third condition and applying Lemma 3.3.4, we see that $dh$ has
only imaginary periods. Clearly the two catenoidal ends
represent the only infinite vertices of the $dh$ orthodisk.
By assumption, all other ends are on one component of the
orthodisk boundary without the catenoidal vertices removed, so
that they have to arrange themselves in order. But the $dh$
integrates to the height coordinate, and so the ends are all at
different height.
\vglue6pt

There is some further evidence that our surfaces are embedded: The Costa
surface is embedded; moreover, there is a periodic surface,
due to Callahan, Hoffman and Meeks  (see Example 3.5.6),
which is also embedded, and which is, at least in a weak
sense, a limit of the Costa towers ${\rm DH}_{m,m}= CT_{2m+1}$.

Theorem~9.2 shows that the surfaces ${\rm DH}_{m,n}$ are embedded
outside of a compact set.  We do not yet have a proof that
these surfaces, in general, are embedded, yet we state next
that any  lack of embeddedness is somewhat inessential. Again
space considerations require us to defer the argument to
\cite{WW2}, where we give two proofs, one exploiting the
algebraic aspects  (in particular we use the Arf invariant --
see \cite{KS}, \cite{Pin}, \cite{HH})  of our construction, and one
exploiting the geometric aspects of our  construction.

\nonumproclaim{Proposition 9.3} The surface ${\rm DH}_{m,n}$ is regularly
homotopic in $\bbe^3$ to an embedding via a compactly supported regular
homotopy that fixes neighborhoods of the critical points of $dh$.
\endproclaim

\AuthorRefNames [hallooo]
\references

[Blo] \name{D. Blo{\ss}},
 Elliptische Funktionen und vollst\"andige Minimalfl\"achen,  Ph.D.\  Thesis,
 Freie Universit\"at,  Berlin,1989.

[Bo] \name{E. Boix}
 Thesis,
  \'Ecole Polytechnique, Palaiseau, France, 1994.

[CHM] \name{M. Callahan, D. Hoffman}, and \name{W.\ H.\  Meeks, III},
 Embedded minimal surfaces with an infinite number of
ends, {\it Invent.\  Math.\/}
 {\bf 96} (1989), 459--505.

[Cos] \name{C. Costa},
 Example of a complete minimal immersion in ${\scriptstyle\bbr}^3$ of
genus one and three embedded ends,
{\it Bull.\ Soc.\ Brasil Mat.\/} {\bf 15} (1984), 47--54.

[EL] \name{J. Eells} and \name{L. Lemaire},
 Deformations of metrics and associated harmonic maps,
{\it Proc.\ Indian Acad.\ Sci.\ Math.\ Sci.\/}
 {\bf 90} (1981), 33--45.

[FLP] \name{A. Fathi, F. Laudenbach}, and \name{V.\  Poenaru},
{\it Traveaux de Thurston sur les Surfaces},  {\it Ast{\rm \'{\it e}}risque}
 {\bf 66--67}, Soc.\ Math.\ France, Paris, France.

[Gar] \name{F.\ P.\ Gardiner}, 
{\it Teichm{\rm \"{\it u}}ller Theory and Quadratic Differentials},
 Wiley Interscience, New York, 1987.

[GM] \name{F.\ P.\  Gardiner} and \name{H. Masur},
 Extremal length geometry of Teichm\"uller space,
{\it Complex Variables Theory Appl.\/} {\bf 16} (1991), 209--237.

[HH] \name{J.\ Hass}  and \name{J.\ Hughes}, 
 Immersions of surfaces in $3$-manifolds,
{\it Topology} {\bf 24} (1985), 97--112.

[HM] \name{J. Hubbard} and \name{H. Masur},
 Quadratic differentials and foliations,
{\it Acta Math.\/}  {\bf 142} (1979), 221--274.

[Ho-Ka]  \name{D. Hoffman} and \name{H. Karcher},
 Complete embedded minimal surfaces of finite total curvature,
{\it Geometry} V  (R.\ Osserman, ed.) {\it Encyclopaedia  Math.\ Sci.\/}
 {\bf 90}, Springer-Verlag, New York, 1997, 5--93, 267--272.

[Ho-Me] \name{D. Hoffman} and \name{W.\ H.\  Meeks, III},
 Embedded minimal surfaces of finite topology,
{\it Ann.\ of Math.\/} {\bf 131} (1990), 1--34.

[J] \name{J.\ A.\  Jenkins},
 On the existence of certain general extremal metrics,
{\it Ann.\ of Math.\/} {\bf 66} (1957), 440--453.

[J-M]  \name{L. Jorge} and \name{W.\ H.\  Meeks, III},
 The topology of complete minimal surfaces of finite total Gaussian
curvature,
{\it Topology}
 {\bf 22} (1983), 203--221.

[Kap] \name{N. Kapouleas},
 Complete embedded minimal surfaces of finite sotal curvature,
{\it J. Differential Geom.\/} {\bf 47}
(1997), 95--169.

[Kar] \name{H. Karcher},
 Construction of minimal surfaces,
 in {\it Surveys in Geometry}, University of Tokyo, 1989, 1--96.

[Ke] \name{S. Kerckhoff},
 The asymptotic geometry of Teichm\"uller space,
{\it Topology} {\bf 19}\break (1980), 23--41.

[K-S] \name{R. Kusner} and \name{N. Schmitt},
 The spinor representation of surfaces in space,
 University of Massachusetts, preprint;
http://front.math.ucdavis.edu/dg-ga/9610005.

[Laws] \name{H.\ B.\  Lawson, Jr.},
  {\it Lectures on Minimal Submanifolds}, 
 Publish or Perish Press, Wilmington, Del., 1980.

[Lop] \name{F.\ J.\  Lopez},
 The classification of complete minimal surfaces with total
curvature greater than $-12\pi$,
{\it Trans.\ Amer.\  Math.\ Soc.\/}
 {\bf 334} (1992),  49--74.

[MW] \name{F. Mart\'\i{}n} and \name{M. Weber},
 On properly embedded minimal surfaces with three ends,
{\it Duke Math.\ J.}
 {\bf 107} (2001), 533--559.

[Oht] \name{M. Ohtsuka},
{\it  Dirichlet Problem{\rm ,} Extremal Length{\rm ,} and Prime Ends},
\ Van Nostrand Reinhold, New York, 1970.

[Oss1] \name{R. Osserman},
 Global properties of minimal surfaces in $E^3$ and $E^n$,
{\it Ann.\  of Math.\/} 
 {\bf 80} (1964), 340--364.

[Oss2] \bibline,
{\it A Survey of Minimal Surfaces},
  2nd edition, Dover Publications Inc., New York, 1986.

[Pin] \name{U. Pinkall},
 Regular homotopy classes of immersed surfaces,
{\it Topology}
 {\bf 24} (1985), 421--434.

[Sch] \name{R. Schoen},
 Uniqueness, symmetry and embeddedness of minimal surfaces,
{\it J.\ Differential Geom.\/}  {\bf 18} (1983), 791--809.

[Str] \name{K. Strebel},
{\it Quadratic Differentials},
 Springer-Verlag, New York, 1984.

[Th] \name{E. Thayer},
 Complete minimal surfaces in Euclidean $3$-space,
 Univ.\ of\break Mass\-achu\-setts Thesis, 1994.

[Tr] \name{M. Traizet},
 An embedded minimal surface with no symmetries,
 Univ.\ de Tours, France,  preprint, 2002.

[W] \name{M. Weber},
 On the Horgan minimal non-surface,
{\it Calc.\ Var.\ Partial Differential Equations}
 {\bf 7} (1998), 373--379.

[WW] \name{M. Weber} and \name{M. Wolf},
 Minimal surfaces of least total curvature and moduli spaces of plane
polygonal arcs,
{\it Geom.\ Funct.\ Anal.\/} 
 {\bf 8} (1998), 1129--1170.

[WW2] \name{M. Weber} and \name{M. Wolf},
{\it Minimal Surfaces, Flat Structures and Teichm{\rm \"{\it u}}ller
Theory},
 in preparation.

[Woh1] \name{M. Wohlgemuth},
 Higher genus mini\-mal sur\-fa\-ces by grow\-ing han\-dles out of a cate\-noid,
{\it Manu\-scripta Math.\/}  {\bf 70}  (1991),  397--428.

[Woh2] \bibline,
 Minimal surfaces of higher genus with finite total curvature,
{\it Arch.\ Rational Mech.\ Anal.\/}
 {\bf 137} (1997), 1--25.

[Wo] \name{M. Wolf},
 On realizing measured foliations via quadratic
differentials of harmonic maps to ${\scriptstyle\Bbb R}$-trees,
{\it J.\ Anal.\ Math.\/} {\bf 68} (1996), 107--120.

\endreferences
\bye
\bigskip

\centerline{(Received June 23, 1998)}
\enddocument